\input ./style/arxiv-general.cfg
\documentclass[MSNbibl,number,citesort,seceqn,dvips]{arxbj}
\makeatletter
   \@ifpackageloaded{graphicx}{}{\usepackage{graphicx}}
\makeatother
\usepackage{upgreek}

\volume{22}
\issue{2}
\pubyear{2016}
\firstpage{1131}
\lastpage{1183}
\doi{10.3150/14-BEJ689}
\docsubty{FLA}

\makeatletter
\def\textitt#1{\emph{#1}}

\def\overset{\stackrel}
\newtheorem{algorithm}{Algorithm}[section]

\newproclaim{assumption}{Assumption}[section]

\newremark{rem}{Remark}[section]
\newtheorem{theorem}{Theorem}[section]
\newremark{remark}{Remark}[section]
\newremark{example}{Example}[section]
\newtheorem{lemma}{Lemma}[section]%
\makeatother

\begin{document}
\begin{frontmatter}

\title{Analysis of the Forward Search using some new results for
martingales and empirical processes}
\runtitle{Asymptotic analysis of the Forward Search}

%
\begin{aug}
\author[A]{\inits{S.}\fnms{S\o ren}~\snm{Johansen}\thanksref{A}\ead[label=e1]{soren.johansen@econ.ku.dk}}
\and
\author[B]{\inits{B.}\fnms{Bent}~\snm{Nielsen}\corref{}\thanksref{B}\ead[label=e3]{bent.nielsen@nuffield.ox.ac.uk}}
\address[A]{Department of Economics, University of Copenhagen \O ster
Farimagsgade 5, 1353 Copenhagen K, Denmark, and CREATES Aarhus
University. \printead{e1}}
\address[B]{Nuffield College and Department of Economics, University of
Oxford \& Institute for New Economic Thinking at the Oxford Martin
School, UK. \printead{e3}}

\runauthor{S. Johansen and B. Nielsen}

\end{aug}
%

%
\received{\smonth{2} \syear{2013}}
%
\revised{\smonth{11} \syear{2014}}


\begin{abstract}
The Forward Search is an iterative algorithm for avoiding outliers in a
regression analysis suggested by
Hadi and Simonoff
(\textit{J. Amer. Statist. Assoc.} \textbf{88} (1993) 1264--1272),
see also
Atkinson and Riani
(\textit{Robust Diagnostic Regression Analysis} (2000) Springer).
The algorithm constructs subsets of ``good'' observations so that the
size of
the subsets increases as the algorithm progresses. It results in a sequence
of regression estimators and forward residuals. Outliers are detected by
monitoring the sequence of forward residuals. We show that the
sequences of
regression estimators and forward residuals converge to Gaussian processes.
The proof involves a new iterated martingale inequality, a theory for
a new
class of weighted and marked empirical processes, the corresponding quantile
process theory, and a fixed point argument to describe the iterative aspect
of the procedure.
\end{abstract}

%
\end{frontmatter}

\section{Introduction}\label{sintro}

\subsection{The Forward Search algorithm}

The Forward Search algorithm was suggested for the multivariate
location model
by Hadi \cite{Hadi} and for multiple regression by Hadi and
Simonoff \cite{HadiSimonoff} and
developed further by Atkinson \cite{Atkinson1994} and
Atkinson and Riani \cite{AtkinsonRiani2000},
see also Atkinson, Riani and
Cerioli \cite
{ARC2010a,ARC2010b}. It is an algorithm for avoiding outliers in a regression
analysis by recursively constructing subsets of ``good'' observations. The
algorithm starts with a robust estimate of the regression parameters
based on
all observations, and constructs the set of observations with the smallest
$m_{0}$ absolute residuals. It continues by estimating the parameters
by least
squares based on the $m_{0}$ observations selected. From this estimate, the
absolute residuals of all observations are computed and ordered. The
$(m_{0}+1)^{\prime}$st largest absolute residual is the forward
residual and
it is used to monitor the algorithm. The set of $m_{0}+1$ observations with
the smallest absolute residuals is the starting point for the next iteration.
The results of the analysis are plots of the recursively estimated forward
residuals and robust parameter estimates. This paper provides an asymptotic
theory for these forward plots when applied to multiple regression
under the
assumption of no outliers.

The Forward Search is used as a diagnostic tool in regression analysis. The
idea is that most observations are ``good'' in the sense that they
conform with
a regression model with symmetric, if not normal, errors. Some observations
may not conform with the model -- they are the outliers. When building a
statistical model, the user can apply the Forward Search in combination with
considerations about the substantive context to decide which
observations are
``good'' and how to treat the ``outliers'' in the analysis. In order to use the
algorithm, we need to understand its properties when all observations are
``good'' with symmetric or even normal errors. Currently this understanding
comes from simulations reported in, for instance, the above mentioned
papers. In
the present paper, we analyse the algorithm using asymptotic tools. In the
future, we hope to analyse the algorithm in the presence of outliers
that may
or may not be of a symmetric nature.

\subsection{Purpose of paper and results}

In this paper, the forward plots are analysed for a multiple regression model.
The model for the ``good'' observations has symmetric zero mean errors with
unknown scale, while the regressors can be stationary as well as
stochastically and deterministically trending. The plots of forward residuals
and estimators are embedded as stochastic processes in $D[0,1]$, and
their asymptotic properties are derived using new results on empirical
processes and martingales. The results can be applied to construct pointwise
and simultaneous confidence bands for the forward plots.

The first result is that the process of forward residuals behaves
asymptotically as if the parameters were known. That is, as the process of
ordered absolute errors from an i.i.d. sample from the error distribution.
Such empirical quantile processes are studied by analysing the empirical
distribution function as an empirical process. In order to show that the
estimation uncertainty is negligible, we introduce a class of weighted and
marked empirical processes, where the weights represent functions of the
regressors and the marks are functions of the regression error. A technical
difficulty is, that because the empirical processes are constructed from
estimated residuals, the argument of the empirical process is stochastically
varying. We develop the theory of such processes, applying and generalizing
the results of Koul and Ossiander \cite{KoulOssiander}.

In the second result, the process of forward residuals is scaled by recursive
estimates of the unknown standard error. The limiting process is
Gaussian and
the variance function is found.

In the study of weighted and marked empirical processes, the well-known method
of replacing the discontinuous processes by their smooth compensators is
applied. The difference is a martingale. To justify this replacement,
some new
iterated exponential martingale inequalities for the variation of the maximum
of finitely many martingales are developed by an iterative application
of an
exponential inequality of Bercu and Touati \cite{BT}.

\subsection{History and background}

The Forward Search starts with a robust estimator. Examples of robust
regression estimators are the least median squares estimator and the least
trimmed squares estimator of Rousseeuw \cite{Rousseuw}. These estimators are
known to
have good breakdown properties, see Rousseeuw and Leroy~\cite{RL1},
Section~3.4 and an
asymptotic theory for the least trimmed squares regression estimator is
provided by V\'{\i}\v{s}ek \cite{Visekb,Viseka,Visekc}. We will allow initial estimators
$\hat{\beta}^{(m_{0})}$ converging at a rate slower than the usual
$n^{1/2}%
$-rate, for the stationary case, as, for example, the least median squares
estimator, which is $n^{1/3}$-consistent in location-scale models.

Broadly speaking, we require three asymptotic tools. First, a theory for
weighted and marked empirical processes to describe the least squares
statistics. Second, an analysis of the corresponding quantile processes to
describe the forward residuals. Third, a fixed point result to describe the
iteration involved.

In empirical process theory, the weights represent functions of the regressors
and the marks are functions of the regression error. The results generalise
those of Johansen and Nielsen \cite{JN2009} who did not allow stochastic variation in
the quantiles and those of Koul and Ossiander \cite{KoulOssiander} who did not allow marks.
The proof combines a chaining argument with iterations of an exponential
inequality for martingales by Bercu and Touati \cite{BT}.

The quantile process theory draws on the exposition of Cs\"{o}rg\H{o} \cite{Cs}.
It is found that in the case of a known variance, the forward residuals
satisfy a Bahadur representation, so that, asymptotically, the forward
residuals have the same distribution as the order statistics of the absolute
regression errors. When the variance is estimated, an additional term appears
in the asymptotic distribution.

The last ingredient is a fixed point result to describe the iterative result.
A single step of the algorithm has been discussed for the
location-scale case
by Johansen and Nielsen \cite{JN2010}. Starting with Bickel \cite{Bickel}, see also Simpson, Ruppert
and
Carroll \cite
{SimpsonRuppertCarrol}, there are a number of asymptotic results for
one-step L- and M-estimators. These are predominantly concerned with objective
functions that have continuous derivatives, thereby excluding hard rejection
as for the one-step Huber-skip function. The Forward Search gives a sequence
of one-step estimators. Because the estimators are based on least
squares in a
sample selected by truncating the residuals, each estimator is a one-step
Huber-skip estimator. Such estimators have been studied by Ruppert and Carroll \cite
{RuppertCarrol}, Johansen and Nielsen \cite{JN2009,JN2010}, Theorem~3.3,
Welsh and Ronchetti \cite{Ronchetti},
and Hawkins and Olive \cite{Olive}.

There appears to be less work on iteration of one-step estimators. The
case of
smooth weights was considered by Dollinger and Staudte \cite{Dollinger}, but the case of
0--1 weights does not appear to have been studied until recently. Cavaliere and Georgiev \cite
{CG} analysed a sequence of Huber-skip estimators for a first order
autoregression with infinite variance errors, while Johansen and Nielsen~\cite{JN2013},
Theorem~3.3, analysed sequences of one-step Huber-skip estimators with
a fixed
critical value. Here we need a critical value which changes with $m$, the
chosen number of observations, so we need a generalisation of the fixed point
result of the latter paper.

Outline of the paper: The model and the Forward Search algorithm are defined
in Section~\ref{smodel}. The main asymptotic results are given in
Section~\ref{smainresults}. The weighted and marked empirical process
results are
given in Section~\ref{sempproc}, while the iterated exponential martingale
inequalities are presented in Section~\ref{smartingale} with proofs
following in
Appendices \ref{AppendixA}--\ref{sFproof}. The proofs of the main results
follow in Appendix~\ref{sGproof}. Finally, Appendix~\ref{st-order}
gives a
result on order statistics of $\mathsf{t}$-distributed variables.

\section{Model and Forward Search algorithm}
\label{smodel}

The multiple regression model is presented, and the Forward Search algorithm
is defined including the forward residual and forward deletion residual.

\subsection{Model}

We assume that $(y_{i},x_{i})$, $i=1,\ldots,n$ satisfy the multiple regression
equation with regressors of dimension $\dim x$
%
\begin{equation}
y_{i}=x_{i}^{\prime}\beta+\varepsilon_{i},\qquad i=1,
\ldots,n. \label{model}%
\end{equation}
The errors, $\varepsilon_{i}$, are assumed independent and identically
distributed with mean zero and variance~$\sigma^{2}$, and $\varepsilon
_{i}/\sigma$ has known density $\mathsf{f}$ and distribution function
$\mathsf{F}(c)=\mathsf{P}(\varepsilon_{i}\leq\sigma c)$. In practice,
the distribution $\mathsf{F}$ will often be standard normal.

The Forward Search is an algorithm based on ordering absolute residuals and
calculation of least squares estimators from the selected observations. Both
these choices implicitly assume a symmetric density. Because, unless symmetry
is assumed, truncating the errors symmetrically gives in general an error
distribution with mean different from zero and hence biased least squares
estimators, at least for the location parameter.

The distribution function of the absolute errors $|\varepsilon
_{i}|/\sigma$ of
a symmetric density is $\mathsf{G}(c)=\mathsf{P}(|\varepsilon_{1}|\leq
\sigma
c)=2\mathsf{F}(c)-1$ with density $\mathsf{g}(c)=2\mathsf{f}(c)$. We define
the quantiles of the absolute errors as%
%
\begin{equation}
c_{\psi}=\mathsf{G}^{-1}(\psi)=\mathsf{F}^{-1}\bigl
\{(1+\psi)/2\bigr\}, \qquad \psi \in{}[0,1[, \label{FGinv}%
\end{equation}
and the truncated moments
%
\begin{equation}
\tau_{\psi}=\int_{-c_{\psi}}^{c_{\psi}}u^{2}
\mathsf{f}(u)\,\mathrm{d}u \quad\mbox{and}\quad\varkappa_{\psi}=\int_{-c_{\psi}}^{c_{\psi}}u^{4}
\mathsf{f}(u)\,\mathrm{d}u. \label{tauandkappa}%
\end{equation}
Then the conditional variance of $\varepsilon_{1}/\sigma$ given
$\{|\varepsilon_{1}|\leq\sigma c\}$ is
%
\begin{equation}
\varsigma_{\psi}^{2}=\tau_{\psi}/\psi.
\label{biascorr}%
\end{equation}
This will serve as a bias correction for the variance estimator based
on the
truncated sample. Using l'H\^{o}pital's rule, it is seen that
%
\begin{equation}
\varsigma_{0}^{2}=0, \qquad \frac{c_{0}^{2}}{\varsigma_{0}^{2}}=3.
\label
{tau0}%
\end{equation}
If $\mathsf{f}=\varphi$ is Gaussian, then $\varsigma_{\psi
}^{2}=1-2c_{\psi
}\varphi(c_{\psi})/\psi$.

\subsection{Forward Search algorithm}

The Forward Search algorithm is designed to avoid outliers in a linear
multiple regression. The first step is given by the choice of a robust
estimator, $\hat{\beta}^{(m_{0})}$, of the regression parameter,
and the choice of the size $m_{0}$ of the initial set of ``good''
observations. The
algorithm generates a sequence of sets of ``good'' observations and least
squares regression estimators based on these. The $(m+1)^{\prime}$st
step of
the algorithm is given as follows.

\begin{algorithm}[(Forward Search)]
\label{alestimators}
\textup{1.} Given an estimator $\hat
{\beta}^{(m)}$ compute absolute residuals $\hat{\xi}_{i}^{(m)}=|y_{i}%
-x_{i}^{\prime}\hat{\beta}^{(m)}|$, $i=1,\ldots,n$.\vspace*{-6pt}
\begin{longlist}[2.]
\item[2.] Find the
$(m+1)^{\prime}$st smallest order statistic $\hat{z}^{(m)}=\hat{\xi}%
_{(m+1)}^{(m)}$.

\item[3.] Find the set of $(m+1)$ observations with smallest
residuals $S^{(m+1)}=(i\dvt \hat{\xi}_{i}^{(m)}\leq\hat{z}^{(m)})$.

\item[4.]
Compute the new least squares estimators on $S^{(m+1)}$
%
\begin{eqnarray}
\hat{\beta}^{(m+1)} & =&\biggl(%
{ \sum
_{i\in S^{(m+1)}}} 
x_{i}x_{i}^{\prime}
\biggr)^{-1}\biggl(%
{ \sum
_{i\in S^{(m+1)}}} 
x_{i}y_{i}
\biggr),\label{beta^m}
\\
\bigl(\hat{\sigma}^{(m+1)}\bigr)^{2} & =&\frac{1}{m+1}%
{
\sum_{i\in S^{(m+1)}}} 
\bigl(y_{i}-x_{i}^{\prime}
\hat{\beta}^{(m+1)}\bigr)^{2}. \label{sigma^m}%
\end{eqnarray}
\end{longlist}
\end{algorithm}

Note, that $\hat{\beta}^{(n)}$ and $(\hat{\sigma}^{(n)})^{2}$ are the full
sample least squares estimators, and that for $n\rightarrow\infty
,m/n\rightarrow\psi$, see Theorem~\ref{tshat},
\[
\bigl(\hat{\sigma}^{(n)}\bigr)^{2}\overset{\mathsf{P}} {
\rightarrow}\sigma^{2}\tau _{\psi
}/\psi.
\]
We therefore introduce also the (asymptotically) bias corrected variance
estimator using $\varsigma_{m/n}^{2}=\tau_{m/n}/(m/n)$, see (\ref{biascorr}),
so that%
%
\begin{equation}
\bigl(\hat{\sigma}_{\mathrm{corr}}^{(m)}\bigr)^{2}=
\frac{(\hat{\sigma
}^{(m)})^{2}}{\varsigma
_{m/n}^{2}}\overset{\mathsf{P}} {\rightarrow}\sigma^{2}.
\label
{shatcorrect}%
\end{equation}
Applying the algorithm for $m=m_{0},\ldots,n-1$, results in sequences
of order
statistics $\hat{z}^{(m)}=\hat{\xi}_{(m+1)}^{(m)}$, least squares estimators
$(\hat{\beta}^{(m)},(\hat{\sigma}^{(m)})^{2})$, along with the scaled forward
residuals
\[
\frac{\hat{z}^{(m)}}{\hat{\sigma}^{(m)}}=\frac{\hat{\xi}_{(m+1)}^{(m)}}%
{\hat{\sigma}^{(m)}}.
\]
Atkinson and Riani \cite{AtkinsonRiani2000} propose to use the minimum deletion residual
\[
\hat{d}^{(m)}=\min_{i\notin S^{(m)}}\hat{\xi}_{i}^{(m)},
\]
instead of the forward residuals. Thus, the deletion residual is based
on the
smallest residual with respect to $\hat{\beta}^{(m)}$ among those observations
that were not included in $S^{(m)}$ which in turn is based on $\hat
{\beta
}^{(m-1)}$, and the forward residual is the largest absolute residual in
$S^{(m+1)}$ which is based on~$\hat{\beta}^{(m)}$.

The plots of $\hat{\beta}^{(m)}$, $\hat{z}^{(m)}/\hat{\sigma}^{(m)}$, and
$\hat{d}^{(m)}/\hat{\sigma}^{(m)}$ against $m$ are called forward
plots, see
Atkinson and Riani \cite{AtkinsonRiani2000}, pages~12--13. The primary objective of this
paper is to
derive the asymptotic distribution of these plots.

When the method was proposed by Hadi and Simonoff \cite{HadiSimonoff}, they also suggested
scaling the residual by a leverage factor and replace the scaled residuals
$\hat{\xi}_{i}^{(m)}/\hat{\sigma}^{(m)}$ above by
\[
\frac{\hat{\xi}_{i}^{(m)}}{\hat{\sigma}^{(m)}\sqrt{1-h_{i}^{(m)}}}\qquad\mbox {for }i\in S^{(m)},\qquad \frac{\hat{\xi}_{i}^{(m)}}{\hat{\sigma}^{(m)}\sqrt
{1+h_{i}^{(m)}}}\qquad\mbox{for
}i\notin S^{(m)},
\]
where $h_{i}^{(m)}=x_{i}^{\prime}(\sum_{j\in S^{(m)}}x_{j}x_{j}^{\prime}
)^{-1}x_{i}$ is the leverage factor. Johansen and Nielsen \cite{JN2009} prove that
such a leverage factor does not change the asymptotic distribution for the
one-step Huber-skip estimator, and the methods presented there can be
used to
prove a similar result for the Forward Search. Another small sample
correction would be to replace $m+1$ with $m+1- \operatorname{dim} x$ in
(\ref{sigma^m}), but we are mainly concerned with asymptotic properties
in this paper.

\section{The main results}
\label{smainresults}

Johansen and Nielsen \cite{JN2010}, Theorems 5.1--5.3, analysed a single step of the
Forward Search applied in a location-scale setting. Those results show that
the one-step version of the scaled residuals $\hat{z}^{(m)}/\hat{\sigma}
^{(m)}$ has an asymptotic representation involving an empirical process
and a
term arising from the estimation error for the variance. The subsequent
analysis shows how this result generalises to a fully iterated Forward
Search. This section first gives the assumptions, then the
results, and
finally presents some simulations. The derivatives of $\mathsf{f}$ are denoted
$\mathsf{\dot{f}}$ and $\mathsf{\ddot{f}}$ and for more complicated
expressions by $\mathrm{d}/\mathrm{d}x$.

\subsection{Assumptions}

In the following, a series of sufficient assumptions are listed for the
asymptotic theory of the Forward Search. When using the Forward
Search, the
density $\mathsf{f}$ is assumed known. The leading case is the normal density,
$\varphi$, but the results are also discussed for the \textsf{t}-density.

\begin{assumption}
\label{asFS}Let $\mathcal{F}_{i}$ be an increasing sequence of $\sigma$
fields such that $\varepsilon_{i-1}$ and $x_{i}$ are $\mathcal{F}_{i-1}%
$-measurable and $\varepsilon_{i}$ is independent of $\mathcal
{F}_{i-1}$ with
symmetric, continuously differentiable density $\mathsf{f}$ which is positive
on the support $\mathsf{F}^{-1}(0)<c<\mathsf{F}^{-1}(1)$ which contains $0$.
For some $0\leq\kappa<\eta\leq1/4$ choose an $r\geq2$ so that $2^{r-1}%
\geq1+(1/4+\kappa-\eta)(1+\dim x)$. Let $q_{0}=1+\max\{2^{r+1},2/(\eta
-\kappa)\}$. Suppose:
\begin{enumerate}[(iii)]
\item[(i)] density satisfies:
\begin{enumerate}[(a)]
\item[(a)] tail monotonicity: $c^{q}\mathsf{f}(c)$, $|c^{q-1}\mathsf{\dot{f}}(c)|$
are decreasing for large $c$ and some $q>q_{0}$;

\item[(b)] quantile process condition: $\gamma=\sup_{c>0}\mathsf{F}%
(c)\{1-\mathsf{F}(c)\}|\mathsf{\dot{f}}(c)|/\{\mathsf{f}(c)\}^{2}<\infty
$;

\item[(c)] unimodality: $\mathsf{\dot{f}}(c)\leq0$ for $c>0$ and $\lim_{c\rightarrow0}\mathsf{\ddot{f}}(c)<0$;

\item[(d)] tail condition: $\{1-\mathsf{F}(c)\}/\{c\mathsf{f}(c)\}=\mathrm{O}(1)$
for $c\rightarrow\infty$;
\end{enumerate}
\item[(ii)] regressors $x_{i}$ are $\mathcal
{F}%
_{i-1}$-measurable and a non-stochastic normalisation matrix $N$ exists so
that
\begin{enumerate}[(a)]
\item[(a)] $\Sigma_{n}=N^{\prime}\sum_{i=1}^{n}x_{i}x_{i}^{\prime}%
N\overset{\mathsf{D}}{\rightarrow}\Sigma\overset{a.s.}{>}0$;

\item[(b)] $\max_{1\leq i\leq n}|n^{1/2-\kappa}N^{\prime}x_{i}|=\mathrm{O}%
_{\mathsf{P}}(1)$;

\item[(c)] $n^{-1}\mathsf{E}
{ \sum_{i=1}^{n}}
|n^{1/2}N^{\prime}x_{i}|^{q_{0}}=\mathrm{O}(1)$;
\end{enumerate}
\item[(iii)] initial
estimator: $N^{-1}(\hat{\beta}^{(m_{0})}-\beta)=\mathrm{O}_{\mathsf{P}%
}(n^{1/4-\eta})$ for some $\eta>0$.
\end{enumerate}
\end{assumption}

\begin{remark}
\textitt{The constant} $q_{0}$ involves the term $\eta-\kappa$ in two ways.
Here $\kappa$ is needed to control $N^{\prime}x_{i}$. If the regressors are
bounded, we can choose $\kappa=0$. This is also the case if the
regressors are
deterministic or of random walk type, see Example~\ref{exasreg}
below. If
$\kappa=0$ and the initial estimator is convergent at the standard rate,
$\eta=1/4$, then $q_{0}$ reduces to $q_{0}=9$ and moments of order $8+$ are
sufficient. Depending on the trade-off between $\kappa,\eta$ and $\dim x$,
moments of order $8+$ may suffice.
\end{remark}

\begin{remark}
\textitt{Assumption}~\ref{asFS}(i) is satisfied for the normal distribution,
see Example~\ref{exnormal} below. For other distributions, the regularity
conditions involve a trade-off between four features: $\eta$,~which indicates
the rate of the initial estimator, $\kappa$, which indicates the order of
magnitude of the maximum of the normalised regressors, and $\dim x$, the
dimension of the regressor. From these quantities a number $r$ is defined,
which controls the number of moments and the smoothness required for the
density $\mathsf{f.}$ The number $r$ is increasing in $\kappa$ and
$\dim x$
and decreasing in $\eta$. The number of required moments, $1+2^{r+1}$, is
larger than 8 in order to control the estimation error for the variance.
Condition (i)(a) is more severe than normally seen in empirical process theory
due to the marks $\varepsilon_{i}^{p}$. Condition~(i)(b) is used in Theorem~\ref{tdhat0}, which builds on Cs\"{o}rg\H{o} \cite{Cs}. Condition~(i)(c) is
needed to ensure that the iterative element of the Forward Search is a
contraction. The unimodality could be relaxed by assuming the
conclusion of
Lemma~\ref{trho}. Condition (i)(d) for Mill's ratio is milder than the
condition employed for kernel density estimation by
Cs\"{o}rg\H{o} \cite{Cs}, page~139.
\end{remark}

\begin{remark}
\textitt{Assumption}~\ref{asFS}(ii). Condition (ii)(a) is standard in
regression analysis and allows for stationary, random walk, and
deterministically trending regressors. Some specific examples are given in
Example~\ref{exasreg} below.
\end{remark}

As part of the proof, a class of weighted and marked empirical
processes are
analysed in Section~\ref{sempproc} and at that point somewhat weaker assumptions
are introduced, see Assumption~\ref{as}.

\begin{example}
\label{exnormal}\textitt{Assumption}~\ref{asFS}(i) \textitt{for the
reference distribution} $\mathsf{f.}$
\begin{longlist}[(a)]
\item[(a)] \textitt{Standard normal
distribution}, $\mathsf{f}=\varphi$. Condition (i) is
satisfied: (i)(a) holds since $c^{q}\varphi(c)=-c^{q-1}\dot
{\varphi
}(c)$ is decreasing for large $c$ for any $q$. (i)(b) holds with $\gamma=1$,
noting $\dot{\varphi}(c)=-c\varphi(c)$ and the Mill's ratio result
$\{(4+c^{2})^{1/2}-c\}/2<\{1-\Phi(c)\}/\varphi(c)<1/c$, see Sampford \cite{Sampford}.
(i)(d)~holds since $\{1-\Phi(c)\}/\{c\varphi(c)\}<1/c^{2}\rightarrow0$ as
$c\rightarrow\infty$.

\item[(b)] \textitt{Scaled distribution}. Consider a
density $\mathsf{f}_{\delta}(c)$ that has variance $\delta^{2}$ but otherwise
satisfies condition (i). Then $\mathsf{f}(c)=\delta\mathsf{f}_{\delta
}(c\delta)$ has unit variance, distribution function $\mathsf{F}%
(c)=\mathsf{F}_{\delta}(c\delta)$ and satisfies condition (i) with
the same
$\gamma$ in part (b).

\item[(c)] \textitt{Scaled $t$-distribution}. The
$t$-distribution with $d>2^{r+1}$ degrees of freedom has
density $\mathsf{f}_{d}(c)=C_{d}(1+c^{2}/d)^{-(d+1)/2}$ with $C_{d}%
=\Gamma\{(d+1)/2\}/\{(d\uppi)^{1/2}\Gamma(d/2)\}$ and variance $\delta_{d}
^{2}=d/(d-2)$. The reference density can be chosen as $\mathsf{f}%
(c)=\mathsf{f}_{d}(c\delta_{d})\delta_{d}$. Due to part (b), it
suffices to
check condition (i) for $\mathsf{f}_{d}$. It holds that $\mathsf{\dot
{f}%
}_{d}(c)=-\gamma h(c)\mathsf{f}_{d}(c)$ where $\gamma=(d+1)/d$ and
$h(c)=c/(1+c^{2}/d)$ so that $\frac{\mathrm{d}}{\mathrm{d}c}\log\mathsf{f}_{d}(c)=-\gamma h(c)$.
Condition (i)(a): for some constants $C$, it holds that $c^{q}\mathsf{f}
_{d}(c)\sim Cc^{q-d-1}$ and $c^{q-1}|\mathsf{\dot{f}}_{d}(c)|\sim Cc^{q-d-3}$
since $h(c)\sim c^{-1}$. Thus $c^{q}\mathsf{f}_{d}(c)$ and $c^{q-1}%
|\mathsf{\dot{f}}_{d}(c)|$ are both declining for large $c$, for $q$
chosen so
that $d+1>q>q_{0}$. (i)(b) holds with the stated $\gamma$ since
$1-c^{-2}%
d/(d+2)<h(c)\{1-\mathsf{F}_{d}(c)\}/\mathsf{f}_{d}(c)<1$, see Soms \cite{Soms}, equation~(3.2).
(i)(c) is well-known to hold. (i)(d) holds since
$\{1-\mathsf{F}_{d}(c)\}/\{c\mathsf{f}_{d}(c)\}<1/\{ch(c)\}\rightarrow
1/d$ as
$c\rightarrow\infty$.
\end{longlist}
\end{example}

\begin{example}
\label{exasreg}\textitt{Assumption}~\ref{asFS}(ii) \textitt{for the
regressors} $x_{i}$.
\begin{longlist}[(a)]
\item[(a)] \textitt{Stationary and autoregressive
regressors}. In this case $x_{i}$ and $\varepsilon_{i}$ have moments of the
same order and $N=n^{-1/2}I_{\dim x}$. (ii)(c) holds if $\mathsf{E}%
|x_{i}|^{q_{0}}<\infty$. (ii)(b) holds due to the Boole and Markov
inequalities if $\eta>\kappa>1/q_{0}$.

\item[(b)] \textitt{Deterministic
regressors} such as $x_{i}=(1,i)^{\prime}$. Let $N=\operatorname{diag}%
(n^{-1/2},n^{-3/2})$. Then $n^{1/2} N^{\prime}x_{i}=(1,i/n)^{\prime}$. Thus
condition (ii) follows with $\kappa=0$.

\item[(c)] \textitt{Random walk
regressors} such as $x_{i}=\sum_{s=1}^{i-1}\varepsilon_{s}$. Let $N=n^{-1}$.
Then $n^{-1/2}x_{\operatorname{int}(n\psi)}$ converges to a Brownian motion by
Donsker's invariance principle, see Billingsley \cite{Billingsley}, Theorem~14.1.
Conditions (ii)(a), (ii)(b)
follows from the continuous mapping theorem with $\kappa=0$. As $x_{i}$ is
defined in terms of $\varepsilon_{i}$ which has moments of order
$q_{0}$, so
has $x_{i}$ and (ii)(c) follows.
\end{longlist}
\end{example}

\begin{example}
\textitt{Assumption}~\ref{asFS}(iii) \textitt{for the initial estimator.}
The focus of this paper is the situation with no outliers.
Thus, a
wide range of $n^{1/2}$-consistent standard estimators or even $n^{1/3}%
$-consistent median based estimators can be used. Therefore, Assumption~\ref{asFS}(iii) only becomes binding when analysing cases with outliers.
\end{example}

\subsection{The results}

The forward plot of, for instance, $\hat{z}^{(m)}$ is a process on
$m=m_{0},\ldots,n-1$. It is useful to embed it in the space $D[0,1]$ of right
continuous process on $[0,1]$ with limits from the left, endowed with the
uniform norm since all limiting processes will be continuous. Thus, define
%
\begin{equation}
\hat{z}_{\psi}=\cases{
\hat{z}^{(m)}, & \quad$\mbox{for }m=\operatorname{int}(n\psi)\mbox{ and
}m_{0}%
/n\leq\psi\leq1,$
\vspace*{2pt}\cr
0, & \quad $\mbox{otherwise.}$}
\label{embedding}%
\end{equation}
Embed in a similar way $\hat{\beta}^{(m)}$, $\hat{\sigma}^{(m)}$ as
$\hat{\beta}_{\psi}$, $\hat{\sigma}_{\psi}$.

The main results are described in terms of three processes
%
\begin{eqnarray}
\mathbb{G}_{n}(c_{\psi}) & =&n^{-1/2}%
{
\sum_{i=1}^{n}} 
\{1_{(|\varepsilon_{i}/\sigma|\leq c_{\psi})}-\psi\},\label{GG}
\\
\mathbb{L}_{n}(c_{\psi}) & =&\tau_{\psi}^{-1}n^{-1/2}%
{
\sum_{i=1}^{n}} 
\bigl[\bigl\{(
\varepsilon_{i}/\sigma)^{2}-c_{\psi}^{2}
\bigr\}1_{(|\varepsilon
_{i}/\sigma|\leq
c_{\psi})}-\bigl(\tau_{\psi}-c_{\psi}^{2}
\psi\bigr)\bigr],\label{HH}
\\
\mathbb{K}_{n}(c_{\psi}) & =&%
{
\sum_{i=1}^{n}} 
N^{\prime}x_{i}
\varepsilon_{i}1_{(|\varepsilon_{i}/\sigma|\leq c_{\psi})}. \label{KK}%
\end{eqnarray}
The first two are asymptotically Gaussian processes and the same holds
for the
third if the regressors are stationary, see Theorem~\ref{tGhat0}.

The main results give asymptotic representations of the forward residuals
$\hat{z}_{\psi}/\sigma$ scaled with known scale, of the bias corrected
variance, $\hat{\sigma}_{\psi,\mathrm{corr}}^{2}$, and of the forward residuals
$\hat{z}_{\psi}/\hat{\sigma}_{\psi,\mathrm{corr}}$ scaled with the bias corrected
variance estimator. Next, it is shown that the forward residuals, $\hat
{z}_{\psi}$, and the deletion residuals, $\hat{d}_{\psi}$, have the same
asymptotic representation after an initial burn-in period. Finally, an
asymptotic representation is given for the forward plot of regression
estimators, $\hat{\beta}_{\psi}$. Proofs of these results are given in
Appendix \ref{sGproof}.

\begin{theorem}
\label{txi}Suppose Assumption~\ref{asFS} holds. Let $\psi_{0}>0$ and
$\omega<\eta-\kappa\leq1/4$. Then%
%
\begin{equation}
\sup_{\psi_{0}\leq\psi\leq n/(n+1)}\bigl|2\mathsf{f}(c_{\psi})n^{1/2}
\bigl(\sigma ^{-1}\hat{z}_{\psi}-c_{\psi}\bigr)+
\mathbb{G}_{n}(c_{\psi})\bigr|=\mathrm{o}%
_{\mathsf{P}}
\bigl(n^{-\omega}\bigr). \label{txibahadur}%
\end{equation}
Moreover, if $\hat{c}_{m/n}$ are the order statistics of $\xi_{i}%
/\sigma=|\varepsilon_{i}|/\sigma$, then%
%
\begin{equation}
\sup_{\psi_{0}\leq\psi\leq n/(n+1)}\bigl|\mathsf{f}(c_{\psi})n^{1/2}\bigl(
\sigma ^{-1}\hat{z}_{\psi}-\hat{c}_{\psi}\bigr)\bigr|=
\mathrm{o}_{\mathsf{P}}\bigl(n^{-\omega}\bigr). \label{txicpsi}%
\end{equation}
\end{theorem}

If $\beta$ and $\sigma$ were known, the residuals are the innovations,
$\varepsilon_{i}$, and the ordering of the absolute residuals $\xi_{i}%
=|y_{i}-\beta^{\prime}x_{i}|=|\varepsilon_{i}|$ can be done once, so that
$\sigma^{-1}\hat{z}_{m}=\sigma^{-1}\xi_{(m+1)}=\hat{c}_{(m+1)/n}$, and the
left-hand side of (\ref{txicpsi}) is trivially zero. In this situation,
(\ref{txibahadur}) reduces to the Bahadur \cite{Bahadur} representation for the
order statistics of the errors $\xi_{i}$, see also Theorem~\ref
{tdhat0} in
the \hyperref[app]{Appendix}. Theorem~\ref{txi} therefore has the interpretation that the
forward residuals $\hat{z}_{m}=\hat{\xi}_{(m+1)}^{(m)}$ behave asymptotically
as the order statistics of the absolute innovations $\xi
_{i}=|\varepsilon
_{i}|$.

\begin{theorem}
\label{tshat}Let $\psi_{0}>0$. Under Assumption~\ref{asFS}, the
asymptotically biased corrected variance estimator has the
representation
\[
\sup_{\psi_{0}\leq\psi\leq n/(n+1)}\bigl|n^{1/2}\bigl(\sigma^{-2}\hat{
\sigma}%
_{\psi,\mathrm{corr}}^{2}-1\bigr)-\mathbb{L}_{n}(c_{\psi})\bigr|=
\mathrm{o}_{\mathsf{P}}(1).
\]
\end{theorem}

\begin{remark}
In Theorems \ref{txi} and \ref{tshat}, the supremum is taken over a smaller
interval for $\psi$ than the unit interval. A left end point larger
than 0 is
needed to ensure consistency. The results potentially hold with a right end
point equal to 1. Proving this would, however, add significantly to the length
of the proof without practical benefit, since the last forward residual is
based on the set $S^{(n-1)}$ with $n-1$ selected observations.
\end{remark}

\begin{remark}
\label{rc2B}The least squares estimator for the variance is $\hat
{\sigma
}_{1,\mathrm{corr}}^{2}=\hat{\sigma}_{1}^{2}$, noting that $\tau_{1}=1$ and
$\varsigma_{1}=1$. Least squares theory shows that $n^{1/2}(\hat{\sigma}
_{1}^{2}/\sigma^{2}-1)=n^{-1/2}\sum_{i=1}^{n}(\varepsilon_{i}^{2}/\sigma
^{2}-1)+\mathrm{o}_{\mathsf{P}}(1)$. To see that Theorem~\ref{tshat} matches
this result, note that the leading term of the least squares
approximation is
$\lim_{\psi\rightarrow1}\tau_{\psi}^{-1}n^{-1/2}%
{ \sum_{i=1}^{n}}
\{(\varepsilon_{i}/\sigma)^{2}1_{(|\varepsilon_{i}/\sigma|\leq c_{\psi
})}%
-\tau_{\psi}\}$. It is therefore necessary that the other term in
$\mathbb{L}_{n}(\psi)$ satisfies
\[
\lim_{\psi\rightarrow1}\tau_{\psi}^{-1}c_{\psi}^{2}n^{-1/2}%
{
\sum_{i=1}^{n}} 
\{1_{(|\varepsilon_{i}/\sigma|\leq c_{\psi})}-\psi\}=\lim_{\psi
\rightarrow
1}c_{\psi}^{2}
\mathbb{G}_{n}(\psi)=\mathrm{o}_{\mathsf{P}}(1).
\]
Because $\varepsilon_{i}$ has more than 8 moments, $c_{\psi}^{2}%
=\mathrm{o}\{(1-\psi)^{-1/4}\}$, see also item 4 of the proof of Lemma~\ref{lshat}. Combine this with Theorems \ref{tOReilly}(a), \ref
{tBB} to
see that $\lim_{\psi\rightarrow1}c_{\psi}^{2}\mathbb{G}_{n}(c_{\psi
})=\mathrm{o}_{\mathsf{P}}(1)$.
\end{remark}

Combining Theorems \ref{txi} and \ref{tshat} gives an asymptotic
representation of the forward residuals with a bias corrected scale.

\begin{theorem}
\label{txisigma}Let $c_{\psi}=\mathsf{G}^{-1}(\psi)$ and $\psi
_{0}>0$. Under
Assumption~\ref{asFS}, the bias corrected scaled forward residual has the
expansion
\[
\sup_{\psi_{0}\leq\psi\leq n/(n+1)}\biggl|2\mathsf{f}(c_{\psi})n^{1/2}
\biggl(\frac
{\hat
{z}_{\psi}}{\hat{\sigma}_{\psi,\mathrm{corr}}}-c_{\psi}\biggr)+\mathbb{G}_{n}(c_{\psi
})+c_{\psi}
\mathsf{f}(c_{\psi})\mathbb{L}_{n}(c_{\psi})\biggr|=
\mathrm{o}%
_{\mathsf{P}}(1).
\]
\end{theorem}

The above results generalise those of Johansen and Nielsen \cite{JN2010}, Theorems 5.1, 5.3,
for a single forward step for location-scale models. It is
interesting to note that the results do not depend on the type of regressors
of the model. This is due to Lemma~\ref{tdhat0}, which for $g=1$ and
$p=0$ shows that
the empirical distribution of the absolute residuals, due to symmetry
of the
density, has an expansion which similarly does not depend on the regressors.

An exception occurs for the empirical process of the residuals
themselves, see the
expansion in Theorem~\ref{tFplin} for $b=\hat{b}$. The expansion in
general depends on the
regressors through the bias term $\hat{b}^{\prime}n^{1/2}N^{\prime
}\bar{%
x}=n^{1/2}(\hat{\beta}-\beta)^{\prime}\bar{x}$, see Lee and Wei \cite{LeeWei}, Theorem~3.2. If, however, the regressors contain a constant, we write $
\beta^{\prime}x_{i}=\mu+\gamma^{\prime}z_{i}$. The first order
condition for estimating $\mu$ is $\bar{y}=\hat{\mu}+\hat{\gamma
}^{\prime}%
\bar{z}$, and inserting $\bar{y}=\bar{\varepsilon}+\mu+\gamma^{\prime
}%
\bar{z}$ we find that $(\hat{\beta}-\beta)^{\prime}\bar{x}=(\hat{\mu
}-\mu
)+(\hat{\gamma}-\gamma)^{\prime}\bar{z}=\bar{\varepsilon}$. This shows,
that including a constant, the bias term does not depend on the other
regressors, $z_{i}$, see Engler and Nielsen \cite{Engler}, Theorem~2.1.

In finite samples the forward residuals and the deletion residuals can be
different, see, for instance, Johansen and Nielsen \cite{JN2010}, Section~2.2. The next
result implies that $\hat{d}^{(m)}$ and $\hat{z}^{(m)}$ have
the same
asymptotic distribution.

\begin{theorem}
\label{tdeletion}It follows from the definitions that $\hat
{d}^{(m)}%
\leq\hat{z}^{(m)}$. Let $m_{0}=\operatorname{int}(n\psi_{0})$ where $\psi
_{0}>0$, and let Assumption~\ref{asFS} hold. Then for all $\psi_{1}$ such
that $\psi_{0}<\psi_{1}<1$
\[
\sup_{\psi_{1}\leq\psi\leq n/(n+1)}\bigl|\mathsf{f}(c_{\psi})n^{1/2}\bigl(
\hat {z}%
^{(m)}-\hat{d}^{(m)}\bigr)\bigr|=
\mathrm{o}_{\mathsf{P}}(1).
\]
\end{theorem}

The last result is for the forward plot of the estimator error $N^{-1}%
(\hat{\beta}^{(m)}-\beta)$, which can be analysed in two stages.
First, it
is established that $N^{-1}(\hat{\beta}^{(m)}-\beta)$ satisfies a
recursion of the form%
%
\begin{equation}
N^{-1}\bigl(\hat{\beta}^{(m+1)}-\beta\bigr)=
\rho_{m/n}N^{-1}\bigl(\hat{\beta}%
^{(m)}-
\beta\bigr)+(\psi\Sigma_{n})^{-1}\mathbb{K}_{n}(c_{\psi})+e_{m/n}%
\bigl\{N^{-1}\bigl(\hat{\beta}^{(m)}-\beta\bigr)\bigr\},
\label{betaar1}%
\end{equation}
where $\rho_{\psi}=2c_{\psi}\mathsf{f}(c_{\psi})/\psi$ is an ``autoregressive
coefficient'' and $e_{\psi}$ is a vanishing remainder term. This result
generalises the result for the location model in Johansen and Nielsen \cite{JN2010}, Theorem~5.2. The unimodality required in Assumption~\ref{asFS}(i)(c)
implies that $\rho_{\psi}$ is bounded away from unity for $\psi\geq\psi_{0}$.
The recursion (\ref{betaar1}) can then be iterated by generalising the
argument in Johansen and Nielsen \cite{JN2013} for the iterated one-step Huber-skip
estimator for a fixed $\psi$. The following result arises.

\begin{theorem}
\label{tbhatexpand}Suppose Assumption~\ref{asFS} holds. Let $m_{0}%
=\operatorname{int}(n\psi_{0})$ where $\psi_{0}>0$. Then, for all $\psi_{1}$,
$\psi_{0}<\psi_{1}<1$, the forward plot of the estimator has the
expansion
\[
\sup_{\psi_{1}\leq\psi\leq1}\biggl|N^{-1}(\hat{\beta}_{\psi}-
\beta)-\frac
{1}{\psi-2c_{\psi}\mathsf{f}(c_{\psi})}\Sigma_{n}^{-1}\mathbb
{K}_{n}(c_{\psi
})\biggr|=\mathrm{o}_{\mathsf{P}}(1).
\]
\end{theorem}

\subsection{Applications of the result for the forward residuals}

The statements of Theorems \ref{txi}, \ref{txisigma}, \ref
{tdeletion} for
the forward residuals and Theorem~\ref{tshat} do not depend on the
type of
regressor. Thus, to apply these theorems it suffices to analyse the
asymptotically Gaussian processes $\mathbb{G}_{n}$ and $\mathbb{L}_{n}$ for
the chosen reference distribution.

\begin{theorem}
\label{tGhat0} Suppose Assumption~\ref{as} holds. Then $\mathbb
{G}_{n}$ and
$\mathbb{L}_{n}$ converge on $D[0,1]$ to zero mean Gaussian processes,
$\mathbb{G}$, $\mathbb{L}$. Their variances are given by%
%
\begin{eqnarray}
\operatorname{\mathsf{Var}}\bigl\{\mathbb{G}(c_{\psi})\bigr\} & =&\psi(1-
\psi),\label{VarG}
\\
\operatorname{\mathsf{Var}}\bigl\{\mathbb{L}(c_{\psi})\bigr\} & =&\frac{1}{\tau_{\psi}^{2}}%
\bigl\{\varkappa_{\psi}-\tau_{\psi}^{2}+c_{\psi}^{2}(1-
\psi) \bigl(c_{\psi}^{2}\psi -2\tau_{\psi}\bigr)\bigr
\},
\\
\operatorname{\mathsf{Cov}}\bigl\{\mathbb{G}(c_{\psi}),\mathbb{L}(c_{\psi})
\bigr\} & =&\frac{1}%
{\tau_{\psi}}\bigl(\tau_{\psi}-c_{\psi}^{2}
\psi\bigr) (1-\psi)<0, \label{CovGL}%
\end{eqnarray}
where the truncated moments $\tau_{\psi}$ and $\varkappa_{\psi}$ are
given in
(\ref{tauandkappa}).
\end{theorem}

The following pointwise results arise for $\psi_{0}\leq\psi\leq\psi
_{1}$, for
some $\psi_{0}>0$ and $\psi_{1}<1$,%
%
\begin{equation}
n^{1/2}\biggl(\frac{\hat{z}_{\psi}}{\hat{\sigma}_{\psi}}-\frac{c_{\psi}}%
{\varsigma_{\psi}}\biggr)=n^{1/2}
\frac{1}{\varsigma_{\psi}}\biggl(\frac{\hat
{z}_{\psi
}-c_{\psi}\hat{\sigma}_{\psi,\mathrm{corr}}}{\hat{\sigma}_{\psi
,\mathrm{corr}}}\biggr)=n^{1/2}
\frac
{1}{\varsigma_{\psi}}\biggl(\frac{\hat{z}_{\psi}}{\hat{\sigma}_{\psi
,\mathrm{corr}}}-c_{\psi
}\biggr)\overset{
\mathsf{D}} {\rightarrow}\mathsf{N}(0,\omega_{\psi}),
\label{normalscaledforwardresiduals}%
\end{equation}
where $\omega_{\psi}$ has contributions from $\hat{z}_{\psi}$, from
$\hat{\sigma}_{\psi,\mathrm{corr}}$, and from their covariance so that
\[
\omega_{\psi}=\frac{1}{\{2\mathsf{f}(c_{\psi})\}^{2}} \bigl[
\operatorname{\mathsf {Var}}\bigl\{%
\mathbb{G}(c_{\psi})\bigr\}+2c_{\psi}\mathsf{f}(c_{\psi})
\operatorname{\mathsf{Cov}}%
\bigl\{\mathbb{G}(c_{\psi}),\mathbb{L}(c_{\psi})
\bigr\}+c_{\psi}^{2}\mathsf{f}%
^{2}(c_{\psi})
\operatorname{\mathsf{Var}}\bigl\{\mathbb{L}(c_{\psi})\bigr\} \bigr].
\]

The above results shed light on some previously suggested distributional
approximations for the deletion residuals. The approximation of Atkinson and Riani \cite{AtkinsonRiani2006},
Theorem~2, has an asymptotic variance that matches
that of the
process $\mathbb{G}$, while omitting the estimation error for the
scale. Riani and Atkinson \cite{RianiAtkinson2007} presented an approximation to the
distribution of the
deletion residuals that comes from order statistics of certain $t$-distributed
variables. Due to Theorem~\ref{torderedt} in Appendix~\ref{st-order}, that
approximation also has an asymptotic variance matching that
of the process~$\mathbb{G}_{n}$.

\begin{example}\label{exnormalasymp}\textitt{Some particular reference distributions.}
\begin{longlist}[(a)]
\item[(a)] \textitt{Standard normal distribution}. If $\mathsf
{f}=\varphi$,
then $c_{\psi}=\Phi^{-1}\{(1+\psi)/2\}$ and
\begin{eqnarray*}
\tau_{\psi} & =&2\int_{0}^{c_{\psi}}x^{2}
\varphi(x)\,\mathrm{d}x= 2\bigl\{\Phi (x)-x\varphi(x)\bigr\}\Big\vert
_{0}^{c_{\psi}}=\psi-2c_{\psi}\varphi
(c_{\psi
}),
\\
\varkappa_{\psi} & =&2\int_{0}^{c_{\psi}}x^{4}
\varphi(x)\,\mathrm{d}x= 2\bigl\{3\Phi(x)-\bigl(x^{3}+3x\bigr)\varphi(x)
\bigr\}\Big\vert _{0}^{c_{\psi}}=3\psi -2\bigl(c_{\psi
}^{3}+3c_{\psi}
\bigr)\varphi(c_{\psi}).
\end{eqnarray*}
\item[(b)]\textitt{Scaled $t$-distribution} with $d$ degrees of
freedom of Example~\ref{exnormal}(c) has density $\mathsf
{f}(c)=\delta
_{d}\mathsf{f}_{d}(c\delta_{d})$ where $\mathsf{f}_{d}$ is the
$t$-density with
$d$ degrees of freedom and variance $\delta_{d}^{2}=d/(d-2)$. Then
$c_{\psi
}=\delta_{d}^{-1}\mathsf{F}_{d}^{-1}\{(1+\psi)/2\}$ and $\psi=2\mathsf
{F}%
_{d}(c_{\psi}\delta_{d})-1$, and
%
\begin{eqnarray}
\tau_{\psi} & = &( d-1 ) \bigl\{ 2\mathsf{F}_{d-2}(c_{\psi
})-1
\bigr\} - ( d-2 ) \bigl\{ 2\mathsf{F}_{d}(c_{\psi}\delta
_{d})-1 \bigr\},\label{taut}
\\
\label{kappat}%
\varkappa_{\psi} & =& ( d-2 ) ^{2} \biggl[ \frac{(d-1)(d-3)}%
{(d-2)(d-4)}
\biggl\{ 2\mathsf{F}_{d-4} \biggl( \frac{c_{\psi}}{\delta_{d-2}
} \biggr) -1 \biggr
\}
\nonumber
\\[-8pt]
\\[-8pt]
\nonumber
&&\hspace*{38pt}{}  -2\frac{d-1}{d-2} \bigl\{ 2\mathsf{F}%
_{d-2} (
c_{\psi} ) -1 \bigr\} + \bigl\{ 2\mathsf{F}_{d} (
c_{\psi}\delta_{d} ) -1 \bigr\} \biggr].
\end{eqnarray}
\end{longlist}
Note that for $c_{\psi}\rightarrow\infty$, the distribution functions approach
unity so that%
%
\begin{equation}
\tau_{\psi}\rightarrow1,\qquad \varkappa_{\psi}\rightarrow3
\frac{d-2}{d-4}, \label{taulimitt}%
\end{equation}
which are the variance and the kurtosis of the scaled $t$-distribution.%


\begin{figure}[b]

\includegraphics{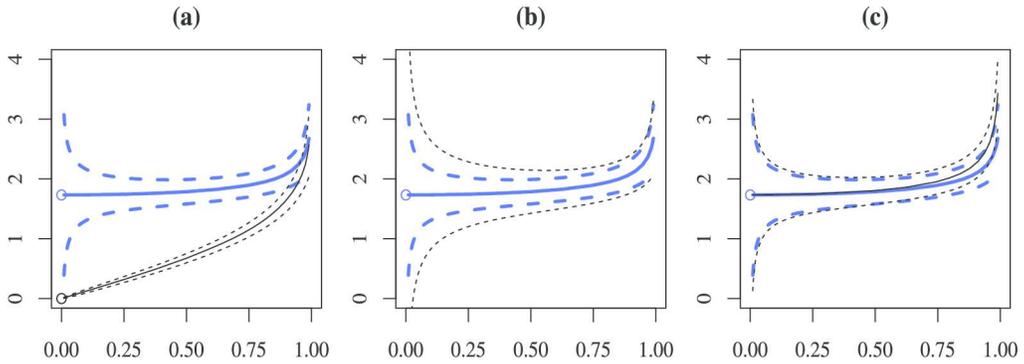}

\caption{Compares the asymptotic distribution of $\hat{z}_{\psi}/\hat
{\sigma
}_{\psi}$ for a normal reference distribution (thick line) with \textup{(a)}
$\hat
{z}_{\psi}/\hat{\sigma}_{\psi,\mathrm{corr}}$ using the corrected scale estimator,
\textup{(b)} $\hat{z}_{\psi}/(\sigma\varsigma_{\psi})$ using the known scale, and
\textup{(c)}~$\hat{z}_{\psi}/\hat{\sigma}_{\psi}$ for a $\mathsf{t}(5)$ reference
distribution. The solid lines indicate the mean, the dashed lines
indicate the
5\% and 95\% asymptotic quantiles for $n=128$.}%
\label{fbias}%
\end{figure}
\end{example}

Figure~\ref{fbias} compares the asymptotic distribution of $\hat
{z}_{\psi
}/\hat{\sigma}_{\psi}$ for a normal reference distribution with (a)
$\hat
{z}_{\psi}/\hat{\sigma}_{\psi,\mathrm{corr}}$ using the corrected scale estimator,
(b) $\hat{z}_{\psi}/(\sigma\varsigma_{\psi})$ using the known scale, and
(c) $\hat{z}_{\psi}/\hat{\sigma}_{\psi}$ for a $\mathsf{t}(5)$ reference
distribution. The solid lines are the point-wise means, while the
dashed lines
are asymptotic 5\% and 95\% quantiles computed for $n=128$. This value
of $n$
is chosen for comparability with the data example in Riani and Atkinson \cite{RianiAtkinson2007}, Figure~1. It is seen that the asymptotic mean $c_{\psi
}/\varsigma_{\psi}$ for
$\hat{z}_{\psi}/\hat{\sigma}_{\psi}$ approaches $\sqrt{3}$ for $\psi
\rightarrow0$, see (\ref{tau0}). Further, the 5\% and 95\% quantiles for
$\hat{z}_{\psi}/\hat{\sigma}_{\psi}$ and $\hat{z}_{\psi}/(\sigma
\varsigma_{\psi})$ diverge for $\psi\rightarrow0$, which is a
consequence of
the division by $\varsigma_{\psi}$ since $\varsigma_{0}=0$, see (\ref{tau0}).
The quantiles also diverge for $\psi\rightarrow1$ which is an extreme
value effect.

In panel (a), the forward residuals $\hat{z}_{\psi}/\hat{\sigma}_{\psi
}$ are
compared to the bias-corrected forward residuals $\hat{z}_{\psi}/\hat
{\sigma
}_{\psi,\mathrm{cor}}$. These representations are equivalent, but the former may be
preferable from a visual viewpoint.

Panel (b) compares situations with estimated and known variance. It is seen
that estimating the variance contributes to reducing the uncertainty. This
phenomenon is also seen for empirical processes of estimated residuals, see
Engler and Nielsen \cite{Engler}, equation (2.10).

Finally, panel (c) compares the result for $\mathsf{f}=\phi$ with the results
for $\mathsf{f}=\mathsf{t}(5)$. With 5 degrees of freedom, Assumption~\ref{asFS} is not met. For higher degrees of freedom, the results will
be in
between the $t_{5}$ and the normal results. A striking feature of this panel
is the excellent agreement between the curves when $\psi$ is not too large.
For larger $\psi$, the long tails of the $t$-distribution have an increasing
effect.\footnote{Graphics were done using R 2.13, see R Development Core Team \cite{R}.}

\subsection{Application of the result for the forward estimators}

In an application of Theorem~\ref{tbhatexpand} for the forward estimators,
the distribution of the kernel $\Sigma_{n}^{-1}\mathbb{K}_{n}(c_{\psi})$
depends on the type of regressors. Building on the analysis in
Johansen and Nielsen~\cite{JN2009},
Sections~1.4, 1.5, we present a result for the stationary
case. For situations with deterministic trends or unit roots, see those
papers. In the case of stationary and autoregressive regressors, we take
$N=n^{-1/2}$ and the normalised matrix of squared regressors, $\Sigma
_{n}=n^{-1}\sum_{i=1}^{n}x_{i}x_{i}^{\prime}$, described in
Assumption~\ref{asFS}(ii)(a), has a deterministic limit.

\begin{theorem}
\label{tK}Suppose Assumption~\ref{as} holds and that $x_{i}$ is stationary
and autoregressive with finite variance. Then $\Sigma_{n}\overset
{\mathsf{P}%
}{\rightarrow}\Sigma>0$ and $\mathbb{K}_{n}$ converges on $D[0,1]$ to a zero
mean Gaussian process $\mathbb{K}$ with variance given as%
%
\begin{equation}
\operatorname{\mathsf{Var}}\bigl\{\mathbb{K}(c_{\psi})\bigr\}=\tau_{\psi}
\sigma^{2}\Sigma.
\end{equation}
\end{theorem}

Theorem~\ref{tK} implies that
\[
n^{1/2}(\hat{\beta}_{\psi}-\beta)\overset{\mathsf{D}} {
\rightarrow }\mathsf{N}\biggl[0,\frac{\tau_{\psi}\sigma^{2}}{\{\psi-2c_{\psi}\mathsf
{f}(c_{\psi
})\}^{2}}\Sigma^{-1}\biggr],
\]
which generalises Johansen and Nielsen \cite{JN2010}, Corollaries 5.2, 5.3. The limiting
distribution matches that of the least trimmed squares estimator with trimming
$\psi$, see V\'{\i}\v{s}ek \cite{Visekc}, Theorem~1.

\section{A class of auxiliary weighted and marked empirical
processes}
\label{sempproc}

It is useful to consider an auxiliary class of weighted and marked empirical
distribution functions for errors $\varepsilon_{i}$ as opposed to absolute
errors $|\varepsilon_{i}|$. The analysis of this class generalises that of
Koul and Ossiander \cite{KoulOssiander} in two respects. First, the standardised estimation
error $b$ is permitted to diverge at a rate of $n^{1/4-\eta}$ rather than
being bounded. Second, non-bounded marks of the type~$\varepsilon
_{i}^{p}$, see also Section~\ref{ssG},
are allowed. These results are therefore of independent interest. This class
of weighted and marked empirical distribution functions is defined for
$b\in\mathbb{R}^{\dim x}$ and $c\in\mathbb{R}$ by%
%
\begin{equation}
\widehat{\mathsf{F}}_{n}^{g,p}(b,c)=\frac{1}{n}%
{
\sum_{i=1}^{n}} 
g_{in}
\varepsilon_{i}^{p}1_{(\varepsilon_{i}\leq\sigma c+x_{in}^{\prime}
b)},\label{Fhat}%
\end{equation}
with normalised regressors $x_{in}=N^{\prime}x_{i}$, weights $g_{in}$
which are measurable with respect to $(\varepsilon_{i-1},\ldots
,\varepsilon
_{1},x_{i},\ldots,x_{1})$, and marks $\varepsilon_{i}^{p}$. By proving results
that hold uniformly in $b$, we can handle the Forward Search. This
allows an
analysis of the order statistics of the residuals at a given step $m$
of the
Forward Search, since the order statistics depend on the previous estimation
error $\hat{b}=N^{-1}(\hat{\beta}^{(m)}-\beta)$, but are scale invariant.
In turn, we can apply the results for the estimation errors $N^{-1}(\hat
{\beta}^{(m+1)}-\beta)$ and $n^{1/2}(\hat{\sigma}_{\mathrm{corr}}^{(m+1)}-\sigma)$.

\subsection{Assumptions}\label{ssas}

We will keep track of the assumptions in a more explicit way than
above. In
the analysis of the one-sided empirical processes, the density $\mathsf
{f}$ is
not necessarily symmetric.

\begin{assumption}
\label{as}Let $\mathcal{F}_{i}$ be an increasing sequence of $\sigma$ fields
so that $\varepsilon_{i-1},x_{i},g_{in}$ are $\mathcal{F}_{i-1}$-measurable
and $\varepsilon_{i}$ is independent of $\mathcal{F}_{i-1}$ with continuously
differentiable density $\mathsf{f}$ which is positive on the
support $\mathsf{F}^{-1}(0)<c<\mathsf{F}^{-1}(1)$ which
contains $0$. Let $p,r,\eta,\kappa,\nu$ be given such that $p,r\in\mathbb{N}_{0}$,
$0\leq\kappa<\eta\leq1/4$ and $\nu\leq1$. Suppose:
\begin{enumerate}[(iii)]
\item[(i)] density
satisfies:
\begin{enumerate}[(a)]
\item[(a)] \textitt{moments}: $\int_{-\infty}^{\infty}|\varepsilon
|^{2^{r}p/\nu
}\mathsf{f}(u)\,\mathrm{d}u<\infty$;

\item[(b)] \textitt{boundedness}: $\{(1+|c|^{\max(0,2^{r}p-1)})\mathsf{f}%
(c)+(1+|c|^{2^{r}p})|\mathsf{\dot{f}}(c)|\}<\infty$;

\item[(c)] smoothness: a $C_{\mathsf{H}}\in\mathbb{N}$ exist such that for all
$a>0$
\[
\frac{\sup_{c\geq a}(1+c^{2^{r}p})\mathsf{f}(c)}{\inf_{0\leq c\leq
a}(1+c^{2^{r}p})\mathsf{f}(c)}\leq C_{\mathsf{H}},\qquad \frac{\sup_{c\leq
-a}(1+|c|^{2^{r}p})\mathsf{f}(c)}{\inf_{-a\leq u\leq0}(1+|c|^{2^{r}%
p})\mathsf{f}(c)}\leq
C_{\mathsf{H}};
\]
\end{enumerate}
\item[(ii)] regressors $x_{i}$ satisfy $\max_{1\leq i\leq n}|n^{1/2-\kappa
}N^{\prime}x_{i}|=\mathrm{O}_{\mathsf{P}}(1)$ for some non-stochastic
normalisation matrix $N$;

\item[(iii)] weights $g_{in}$ are matrix valued
and satisfy:
\begin{enumerate}[(a)]
\item[(a)] $n^{-1}\mathsf{E}%
{ \sum_{i=1}^{n}}
|g_{in}|^{2^{r}}(1+|n^{1/2}N^{\prime}x_{i}|)=\mathrm{O}(1)$;

\item[(b)] $n^{-1}%
{ \sum_{i=1}^{n}}
|g_{in}|(1+|n^{1/2}N^{\prime}x_{i}|^{2})=\mathrm{O}_{\mathsf{P}}(1)$.\vadjust{\goodbreak}
\end{enumerate}
\end{enumerate}
\end{assumption}

\begin{remark}
\label{ras}Discussion of Assumption~\ref{as}.
\begin{longlist}[(a)]
\item[(a)] \textitt{The case
of no marks} $p=0$. This is the situation discussed in Koul and Ossiander \cite
{KoulOssiander}. The primary role of $r$ is to control the tail
behaviour of the
density. When $p=0$ then $2^{r}p=0$ for all $r\in\mathbb{N}_{0}$, so
$r$ can
be chosen as $r=0$ and the assumptions simplify considerably.

\item[(b)]
\textitt{The moment condition in Assumption}~\ref{as}(i)(a) is used for some
$\nu<1$ for the tightness result in Theorem~\ref{tFptight}. Otherwise,
$\nu=1$ suffices.

\item[(c)] \textitt{The smoothness of density in Assumption}~\ref{as}(i)(c) is satisfied if $\mathsf{h}_{r}(c)=(1+\epsilon^{2^{r}%
p})\mathsf{f}(\epsilon)$ is monotone for $|c|>d_{1}$ for some $d_{1}\geq0$.
Indeed, choose $d_{2}\geq d_{1}$ so that $\sup_{c\geq d_{2}}\mathsf{h}%
_{r}(c)=\inf_{0\leq c\leq d_{2}}\mathsf{h}_{r}(c)=\mathsf
{h}_{r}(d_{2})$. Then
choose
\[
C_{\mathsf{H}}>\sup_{0\leq c\leq d_{2}}\mathsf{h}_{r}(c)\big/\inf
_{0\leq
c\leq
d_{2}}\mathsf{h}_{r}(c).
\]
A similar argument applies for $c<0$. Note, that the smoothness condition
implies that the density has connected support.

\item[(d)]
\textitt{Sufficient condition for Assumption}~\ref{as}(i). If $\mathsf
{f}$ is
symmetric and differentiable with $c^{q}\mathsf{f}(c)$, $c^{q-1}%
|\mathsf{\dot{f}}(c)|$ both decreasing for large $c$ for some $q>1+2^{r}p$,
then Assumption~\ref{as}(i) holds. Indeed, (i)(a) holds, since when
$c^{q}\mathsf{f}(c)$ is decreasing, then $c^{2^{r}p/\nu}\mathsf{f}(c)$ is
integrable for some $\nu<1$. Further, (i)(b) holds, since, first, the
continuity and decreasingness of $c^{q}\mathsf{f}(c)$ and hence of
$\mathsf{f}(c)$ implies $(1+|c|^{1+2^{r}p})\mathsf{f}(c)$ is bounded, and,
second, since $\mathsf{\dot{f}}(c)<0$ for large $c$ so that $|c^{q-1}%
\mathsf{\dot{f}}(c)|$ decreases, then $(1+|c|^{2^{r}p})|\mathsf{\dot{f}}(c)|$
is bounded. Finally, (i)(c) holds due to remark (c) above.
\end{longlist}
\end{remark}

\subsection{The empirical process results}

The weighted and marked empirical distribution function $\widehat
{\mathsf{F}%
}_{n}^{g,p}(b,c)$ defined in (\ref{Fhat}) is analysed through martingale
arguments. Thus, introduce the sum of conditional expectations%
%
\begin{equation}
\overline{\mathsf{F}}_{n}^{g,p}(b,c)=\frac{1}{n}%
{
\sum_{i=1}^{n}} 
g_{in}
\mathsf{E}_{i-1}\bigl\{\varepsilon_{i}^{p}1_{(\varepsilon_{i}\leq\sigma
c+x_{in}^{\prime}b)}
\bigr\}, \label{Fbar}%
\end{equation}
and the weighted and marked empirical process%
%
\begin{equation}
\mathbb{F}_{n}^{g,p}(b,c)=n^{1/2}\bigl\{\widehat{
\mathsf{F}}_{n}^{g,p}%
(b,c)-\overline{
\mathsf{F}}_{n}^{g,p}(b,c)\bigr\}. \label{Ftilde}%
\end{equation}
Three results follows. These are proved in Appendix~\ref{sFproof}. The first
result shows that the dependence of $\mathbb{F}_{n}^{g,p}$ on the estimation
error $b$ is negligible.

\begin{theorem}
\label{tFpest}Let $c_{\psi}=\mathsf{F}^{-1}(\psi)$.
Suppose Assumption~\ref{as}\textup{(i), (ii), (iii)(a)} holds with $\nu=1$, some $\eta>0$ and an $r$ such that
$2^{r-1}\geq1+(1/4+\kappa-\eta)(1+\dim x)$. Then, for any $B>0$ and
$n\rightarrow\infty$
\[
\sup_{0\leq\psi\leq1}\sup_{|b|\leq n^{1/4-\eta}B}\bigl|\mathbb{F}_{n}%
^{g,p}(b,c_{\psi})-
\mathbb{F}_{n}^{g,p}(0,c_{\psi})\bigr|=
\mathrm{o}_{\mathsf
{P}%
}(1).
\]
\end{theorem}

For the standard empirical process with weights $g_{in}=1$ and marks
$\varepsilon_{i}^{p}=1$, the order of the remainder term can be
improved as
follows. Note that when $p=0$, then $r$ will be irrelevant in Assumption~\ref{as}(i), see also Remark~\ref{ras}(a).

\begin{theorem}
\label{tFpestomega}Let $c_{\psi}=\mathsf{F}^{-1}(\psi)$.
Under Assumption~\ref{as}\textup{(i)},
\textup{(ii)}, \textup{(iii)(a)} with $\nu=1$, $p=0$, $r=2$ and some $\eta>0$ it holds
that for any $B>0$, any $\omega<\eta-\kappa\leq1/4$ and $n\rightarrow
\infty$,
\[
\sup_{0\leq\psi\leq1}\sup_{|b|,|d|\leq n^{1/4+\kappa-\eta}B}\bigl|\mathbb{F}%
_{n}^{1,0}
\bigl(b,c_{\psi}+n^{\kappa-1/2}d\bigr)-\mathbb{F}_{n}^{1,0}(0,c_{\psi
})\bigr|=
\mathrm{o}_{\mathsf{P}}\bigl(n^{-1/8-\omega/2}\bigr).
\]
\end{theorem}

The next results presents a linearization of $\overline{\mathsf{F}}_{n}%
^{g,p}(b,c)$.

\begin{theorem}
\label{tFplin}Let $c_{\psi}=\mathsf{F}^{-1}(\psi)$. Suppose
Assumption~\ref{as}\textup{(i)(b)}, \textup{(iii)(b)} holds with $r=0$ and some $\eta>0$. Then, for all $B>0$
and $n\rightarrow\infty$
\[
\sup_{0\leq\psi\leq1}\sup_{|b|\leq n^{1/4-\eta}B}\Biggl|n^{1/2}\bigl
\{\overline {\mathsf{F}}_{n}^{g,p}(b,c_{\psi})-
\overline{\mathsf {F}}_{n}^{g,p}(0,c_{\psi
})\bigr\}-
\sigma^{p-1}c_{\psi}^{p}\mathsf{f}(c_{\psi})n^{-1}%
{
\sum_{i=1}^{n}} 
g_{in}n^{1/2}x_{in}^{\prime}b\Biggr|
\]
is $\mathrm{O}_{\mathsf{P}}(n^{-2\eta})$.
\end{theorem}

Finally, we argue that the weighted and marked empirical process
$\mathbb{F}_{n}^{g,p}(0,c_{\psi})$ in (\ref{Ftilde}) is tight when viewed as a sequence
in $n$ of processes on $D[0,1]$. Following Billingsley \cite{Billingsley}, Theorem~13.2, we
need to check two conditions. First, it holds by construction that
$\mathbb{F}_{n}^{g,p}(0,c_{0})=0$. Second, the next results shows that the
modulus of continuity is small.

\begin{theorem}
\label{tFptight}Let $c_{\psi}=\mathsf{F}^{-1}(\psi)$.
Under Assumption~\ref{as}\textup{(i)(a)}, \textup{(iii)(a)} with $r=2$ and some $\nu<1$ it holds that, for all
$\epsilon>0$,
\[
\lim_{\phi\downarrow0}\mathop{\lim\sup}_{n\rightarrow\infty}
\mathsf{P} \Bigl\{ \sup_{0\leq\psi\leq\psi^{\dag}\leq1\dvt \psi^{\dag}-\psi
\leq\phi}\bigl|\mathbb{F}_{n}^{g,p}(0,c_{\psi^{\dag}})-
\mathbb{F}_{n}%
^{g,p}(0,c_{\psi})\bigr|>\epsilon
\Bigr\} \rightarrow0.
\]
\end{theorem}

The proofs of these results are given in Appendix~\ref{sFproof}.

\section{Iterated exponential martingale inequalities}
\label{smartingale}

Chaining arguments will be used to handle tightness properties of the
empirical processes. This reduces the tightness proof to a problem of finding
the tail probability for the maximum of a certain family of
martingales. We
first give a general result on a bound of a finite number of martingales,
which we prove by iterating a martingale inequality by Bercu and Touati \cite{BT}.
Subsequently, two special cases are analysed: where the number of
elements in the martingale family is increasing and where it is fixed.

\begin{theorem}
\label{lMGineq}For $\ell$, $1\leq\ell\leq L$, let $z_{\ell,i}$ be
$\mathcal{F}_{i}$-adapted and $\mathsf{E}z_{\ell,i}^{2^{\bar
{r}}}<\infty$
for some $\bar{r}\in\mathbb{N}$. Let $D_{r}=\max_{1\leq\ell\leq L}%
{ \sum_{i=1}^{n}}
\mathsf{E}_{i-1}z_{\ell,i}^{2^{r}}$ for $1\leq r\leq\bar{r}$.
Then, for
all $\kappa_{0},\kappa_{1},\ldots,\kappa_{\bar{r}}>0$,
\[
\mathsf{P} \Biggl\{ \max_{1\leq\ell\leq L}\Biggl|%
{
\sum_{i=1}^{n}} 
(z_{\ell,i}-
\mathsf{E}_{i-1}z_{\ell,i})\Biggr|>\kappa_{0} \Biggr\} \leq
L\frac{\mathsf{E}D_{\bar{r}}}{\kappa_{\bar{r}}}+%
{ \sum
_{r=1}^{\bar{r}}} 
\frac{\mathsf{E}D_{r}}{\kappa_{r}}+2L%
{
\sum_{r=0}^{\bar{r}-1}} 
\exp\biggl(-
\frac{\kappa_{r}^{2}}{14\kappa_{r+1}}\biggr).
\]
\end{theorem}

The proof is given in Appendix \ref{AppendixA}.

\begin{theorem}
\label{lsupMrasymp}For $\ell$, $1\leq\ell\leq L$, let $z_{\ell,i}$ be
$\mathcal{F}_{i}$-adapted and $\mathsf{E}z_{\ell,i}^{2^{\bar
{r}}}<\infty$
for some $\bar{r}\in\mathbb{N}$. Let $D_{r}=\max_{1\leq\ell\leq L}%
{ \sum_{i=1}^{n}}
\mathsf{E}_{i-1}z_{\ell,i}^{2^{r}}$ for $1\leq r\leq\bar{r}$. Suppose,
for some $\varsigma\geq0$, $\lambda>0$, that $L=\mathrm{O}(n^{\lambda
})$ and
$\mathsf{E}D_{r}=\mathrm{O}(n^{\varsigma})$ for $r\leq\bar{r}$.
Then, if
$\upsilon>0$ is chosen such that:
\begin{longlist}[(ii)]
\item[(i)] $\varsigma<2\upsilon$,

\item[(ii)] $\varsigma+\lambda<\upsilon2^{\bar{r}}$,
\end{longlist}
it holds
that for
all $\kappa>0$ and $n\rightarrow\infty$,
\[
\lim_{n\rightarrow\infty}\mathsf{P} \Biggl\{ \max_{1\leq\ell\leq L}\Biggl|%
{
\sum_{i=1}^{n}} 
(z_{\ell,i}-
\mathsf{E}_{i-1}z_{\ell,i})\Biggr|>\kappa n^{\upsilon} \Biggr\} =0.
\]
\end{theorem}

\begin{pf}
Apply Theorem~\ref{lMGineq} with
$\kappa_{q}=(\kappa n^{\upsilon})^{2^{q}}(28\lambda\log n)^{1-2^{q}}$
for any
$\kappa>0$ so that $\kappa_{0}=\kappa n^{\upsilon}$ and $\kappa_{q}^{2}%
/\kappa_{q+1}=28\lambda\log n$ and exploit conditions (i), (ii) to see
that the
probability of interest satisfies
\[
\mathcal{P}_{n}\mathcal{=}\mathrm{O} \Biggl\{ n^{\lambda}
\frac
{n^{\varsigma
}(\log n)^{2^{\bar{r}}-1}}{n^{\upsilon2^{\bar{r}}}}+%
{ \sum
_{r=1}^{\bar{r}}} 
\frac{n^{\varsigma}(\log n)^{2^{r}-1}}{n^{\upsilon2^{r}}}+2n^{\lambda
}
\bar{r}n^{-2\lambda} \Biggr\} =\mathrm{o}(1),
\]
as desired since $\varsigma+\lambda<\upsilon2^{\bar{r}}$ and
$\varsigma<2\upsilon\leq\upsilon2^{r}$ for $r\geq1$.
\end{pf}

\begin{theorem}
\label{lsupM4finite}For $\ell$, $1\leq\ell\leq L$, let $z_{\ell,i}$ be
$\mathcal{F}_{i}$-adapted and $\mathsf{E}z_{\ell,i}^{4}<\infty$.
Suppose
$\mathsf{E}\max_{1\leq\ell\leq L}%
{ \sum_{i=1}^{n}}
\mathsf{E}_{i-1}z_{\ell,i}^{2^{q}}\leq Dn$ for $q=1,2$ and some $D>0$. Then,
for all $\theta,\kappa>0$,
\[
\mathsf{P} \Biggl\{ \max_{1\leq\ell\leq L}\Biggl|%
{
\sum_{i=1}^{n}} 
(z_{\ell,i}-
\mathsf{E}_{i-1}z_{\ell,i})\Biggr|>\kappa n^{1/2} \Biggr\} \leq
\frac{(L+1)\theta^{3}D}{\kappa n}+\frac{\theta D}{\kappa}+4L\exp\biggl(-\frac
{\kappa\theta}{14}\biggr).
\]
\end{theorem}

\begin{pf}
Apply Theorem~\ref{lMGineq} with
$\kappa_{q}=\kappa n^{2^{q-1}}\theta^{1-2^{q}}$ for any $\kappa,\theta
>0$ so
that $\kappa_{0}=\kappa n^{1/2}$ and $\kappa_{q}^{2}/\kappa_{q+1}=\kappa
\theta$, while $\bar{r}=2$, to get the bound
\[
\mathcal{P}\leq\frac{(L+1)\theta^{3}}{\kappa n^{2}}\mathsf{E}\max_{1\leq
\ell\leq L}%
{
\sum_{i=1}^{n}} 
\mathsf{E}_{i-1}z_{\ell,i}^{4}+\frac{\theta}{\kappa n}
\mathsf{E}\max_{1\leq\ell\leq L}%
{ \sum
_{i=1}^{n}} 
\mathsf{E}_{i-1}z_{\ell,i}^{2}+4L\exp\biggl(-
\frac{\kappa\theta}{14}\biggr).
\]
Exploit the moment conditions to get the desired result.
\end{pf}

\section{Conclusion}

The intention of the Forward Search is to determine the number of
outliers by
looking at the forward plot of the forward residuals. The main results
for the
Forward Search, given in Section~\ref{smainresults}, describe the asymptotic
distribution of that process in a situation where there are no
outliers. We
can therefore add pointwise confidence bands to the forward plot, using
Theorem~\ref{txisigma}. These give an impression of the pointwise variation
we would expect for the forward plot, if there were in fact no
outliers. In
practice we would want to make a simultaneous decision based on the entire
graph. A theory is developed in Johansen and Nielsen \cite{JN2014b} and implemented
in the R-package ForwardSearch, see Nielsen \cite{N2014}.

We suspect that the iterated martingale inequalities will be useful in a
variety of situations. For instance, in ongoing research, we are
finding that
the inequalities are helpful in establishing consistency and asymptotic
distribution results for general M-estimators, see Johansen and Nielsen~\cite{JN2014a}.

The results and techniques in this paper could potentially also be used to
shed light on other iterative 1-step methods in robust statistics such as
those discussed in Bickel \cite{Bickel}, Simpson, Ruppert and
Carroll \cite{SimpsonRuppertCarrol}, and
Hawkins and Olive \cite{Olive}. Another example would be to establish an asymptotic
theory for the Forward Search applied to multivariate location and scatter,
see Cerioli, Farcomini and
Riani \cite{CeroliFarcomeniRiani} for a discussion of consistency
as well as Riani, Atkinson and Cerioli~\cite{RAC2009}.
Finally, we mention Bellini~\cite{Bellini} for an application of the Forward
Search to the cointegrated vector autoregressive model.

\begin{appendix}\label{app}

\section{Proofs of martingale inequalities}\label{AppendixA}
\begin{pf*}{Proof of Theorem~\ref{lMGineq}}
1. \textit{Notation}. For $0\leq
r\leq\bar{r}$ define $A_{\ell,r}=%
{ \sum_{i=1}^{n}}
(z_{\ell,i}^{2^{r}}-\mathsf{E}_{i-1}z_{\ell,i}^{2^{r}})$ and
\[
\mathcal{P}_{r}(\kappa_{r})=\mathsf{P}\Bigl(\max
_{1\leq\ell\leq L}A_{\ell,r}%
>\kappa_{r}\Bigr),\qquad
\mathcal{Q}_{r}(\kappa_{r})=\mathsf{P}\Bigl(\max
_{1\leq
\ell\leq
L}|A_{\ell,r}|>\kappa_{r}\Bigr),
\]
where $\mathcal{Q}_{0}(\kappa_{0})$ is the probability of interest, while
$\mathcal{P}_{r}(\kappa_{r})\leq\mathcal{Q}_{r}(\kappa_{r})$.

2. \textitt{The terms} $\mathcal{Q}_{r}(\kappa_{r})$ \textitt{for} $0\leq
r<\bar{r}$. We first prove that, for any $\kappa_{r},\kappa
_{r+1}>0$,%
%
\begin{equation}
\mathcal{Q}_{r}(\kappa_{r})\leq2L\exp\biggl(-
\frac{\kappa_{r}^{2}}{14\kappa
_{r+1}%
}\biggr)+\mathcal{P}_{r+1}(\kappa_{r+1})+
\frac{\mathsf{E}D_{r+1}}{\kappa_{r+1}}. \label{plMGineqPm}%
\end{equation}
The idea is now to apply the following inequality for sets $\mathcal{A}%
,\mathcal{B}$
\[
\mathsf{P}(\mathcal{A})=\mathsf{P}(\mathcal{A}\cap\mathcal{B})+\mathsf
{P}%
\bigl(\mathcal{A}\cap\mathcal{B}^{c}\bigr)\leq
\mathsf{P}(\mathcal{A}\cap\mathcal {B}%
)+\mathsf{P}\bigl(
\mathcal{B}^{c}\bigr).
\]
In the first term, $\mathcal{A}$ relates to the tails of a martingale and
$\mathcal{B}$ to the central part of the distribution of the quadratic
variation. Thus, the first term can be controlled by a martingale inequality.
In the second term, $\mathcal{B}^{c}$ relates to the tail of the quadratic
variation. The sum of the predictable and the total quadratic variation of
$A_{\ell,r}$ is $B_{\ell,r}=%
{ \sum_{i=1}^{n}}
B_{\ell,r,i}$ where $B_{\ell,r,i}=(z_{\ell,i}^{2^{r}}-\mathsf{E}_{i-1}%
z_{\ell,i}^{2^{r}})^{2}+\mathsf{E}_{i-1}(z_{\ell,i}^{2^{r}}-\mathsf{E}%
_{i-1}z_{\ell,i}^{2^{r}})^{2}$. We then get%
%
\begin{equation}
\mathcal{Q}_{r}(\kappa_{r})\leq\mathsf{P}\Bigl\{\Bigl(\max
_{1\leq\ell\leq
L}|A_{\ell
,r}|>\kappa_{r}\Bigr)\cap\Bigl(
\max_{1\leq\ell\leq L}B_{\ell,r}\leq7\kappa _{r+1}\Bigr)
\Bigr\}+\mathsf{P}\Bigl(\max_{1\leq\ell\leq L}B_{\ell,r}>7
\kappa_{r+1}\Bigr). \label{plMGineqPm1}%
\end{equation}

Consider the first term in (\ref{plMGineqPm1}), $\mathcal{S}_{1,r}$ say.
By Boole's inequality this satisfies
\[
\mathcal{S}_{1,r}\leq%
{ \sum
_{\ell=1}^{L}} 
\mathsf{P}\Bigl
\{\bigl(|A_{\ell,r}|>\kappa_{r}\bigr)\cap\Bigl(\max_{1\leq\ell\leq L}B_{\ell
,r}%
\leq7\kappa_{r+1}\Bigr)\Bigr\}.
\]
Noting that $(\max_{1\leq\ell\leq L}B_{\ell,r}\leq7\kappa_{r+1})\subset
(B_{\ell,r}\leq7\kappa_{r+1})$ gives the further bound
\[
\mathcal{S}_{1,r}\leq%
{ \sum
_{\ell=1}^{L}} 
\mathsf{P}\bigl
\{\bigl(|A_{\ell,r}|>\kappa_{r}\bigr)\cap(B_{\ell,r}<7
\kappa_{r+1})\bigr\}.
\]
Because $A_{\ell,r}$ is a martingale, the exponential inequality of
Bercu and Touati \cite{BT}, Theorem~2.1, shows
\[
\mathsf{P}\bigl\{\bigl(|A_{\ell,r}|>\kappa_{r}\bigr)
\cap(B_{\ell,r}<7\kappa_{r+1})\bigr\} \leq 2\exp\bigl\{-
\kappa_{r}^{2}/(14\kappa_{r+1})\bigr\}.
\]
Taken $L$ times, this gives the first term in (\ref{plMGineqPm}).

Consider the second term in (\ref{plMGineqPm1}), $\mathcal
{S}_{2,r}$ say.
Ignore the indices on $B_{\ell,r,i}, E_{i-1}$ and $z_{\ell
,i}^{2^{r}}$, and
apply the inequality $(z-\mathsf{E}z)^{2}\leq2(z^{2}+\mathsf{E}^{2}z)$ along
with $\mathsf{E}^{2}z\leq\mathsf{E}z^{2}$ and $\mathsf{E}(z-\mathsf{E}%
z)^{2}\leq\mathsf{E}z^{2}$ to get that $B=(z-\mathsf{E}z)^{2}+\mathsf{E}
(z-\mathsf{E}z)^{2}\leq2z^{2}+3\mathsf{E}z^{2}=2(z^{2}-\mathsf{E}%
z^{2})+5\mathsf{E}z^{2}$. Thus,
\[
\mathcal{S}_{2,r}\leq\mathsf{P}\Biggl\{\max_{1\leq\ell\leq L}%
{
\sum_{i=1}^{n}} 
\bigl(z_{\ell,i}^{2^{r+1}}-\mathsf{E}_{i-1}z_{\ell,i}^{2^{r+1}}
\bigr)>\kappa _{r+1}\Biggr\}+\mathsf{P}\Biggl(\max_{1\leq\ell\leq L}%
{
\sum_{i=1}^{n}} 
\mathsf{E}_{i-1}z_{\ell,i}^{2^{r+1}}>\kappa_{r+1}
\Biggr).
\]
Use the notation from above and then the Markov inequality to get
\[
\mathcal{S}_{2,r}\leq\mathcal{P}_{r+1}(\kappa_{r+1})+
\mathsf{P}(D_{r+1}%
>\kappa_{r+1})\leq
\mathcal{P}_{r+1}(\kappa_{r+1})+\frac{1}{\kappa_{r+1}%
}
\mathsf{E}D_{r+1},
\]
which are the last two terms of (\ref{plMGineqPm}).

3. \textitt{The term} $\mathcal{P}_{\bar{r}}(\kappa_{\bar{r}})$.
Apply the inequality $|z|-\mathsf{E}_{i-1}|z|\leq|z|$ and then Boole's and
Markov's inequalities to get
\[
\mathcal{P}_{\bar{r}}(\kappa_{\bar{r}})\leq\mathsf{P}\Biggl(\max
_{1\leq\ell\leq L}%
{ \sum
_{i=1}^{n}} 
z_{\ell,i}^{2^{\bar{r}}}>
\kappa_{\bar{r}}\Biggr)\leq L\max_{1\leq\ell\leq
L}\mathsf{P}
\Biggl(%
{ \sum_{i=1}^{n}}
z_{\ell,i}^{2^{\bar{r}}}>\kappa_{\bar{r}}\Biggr)\leq
\frac{L}{\kappa
_{\bar{r}}}\max_{1\leq\ell\leq L}\mathsf{E}%
{
\sum_{i=1}^{n}} 
z_{\ell,i}^{2^{\bar{r}}}.
\]
Apply iterated expectations and interchange maximum and expectation to
get%
%
\begin{equation}
\mathcal{P}_{\bar{r}}(\kappa_{\bar{r}})\leq\frac{L}{\kappa
_{\bar{r}}}\max
_{1\leq\ell\leq L}\mathsf{E}%
{ \sum
_{i=1}^{n}} 
\mathsf{E}_{i-1}z_{\ell,i}^{2^{\bar{r}}}\leq\frac{L}{\kappa_{\bar
{r}}%
}
\mathsf{E}\max_{1\leq\ell\leq L}%
{ \sum
_{i=1}^{n}} 
\mathsf{E}_{i-1}z_{\ell,i}^{2^{\bar{r}}}=\frac{L}{\kappa_{\bar{r}}%
}
\mathsf{E}D_{\bar{r}}. \label{plMGineqPrbar}%
\end{equation}
4. \textitt{Combine expressions}. Since $\mathcal{P}_{r+1}(\kappa_{r+1}%
)\leq\mathcal{Q}_{r+1}(\kappa_{r+1})$ then write (\ref{plMGineqPm})
as%
%
\begin{eqnarray}
\mathcal{Q}_{r}(\kappa_{r}) & \leq&2L\exp\biggl(-
\frac{\kappa_{r}^{2}}%
{14\kappa_{r+1}}\biggr)+\mathcal{Q}_{r+1}(\kappa_{r+1})+
\frac{\mathsf
{E}D_{r+1}%
}{\kappa_{r+1}} \qquad\mbox{for }r=0,\ldots,\bar{r}-2,\label
{plMGineqQr}
\\
\mathcal{Q}_{r}(\kappa_{r}) & \leq&2L\exp\biggl(-
\frac{\kappa_{r}^{2}}%
{14\kappa_{r+1}}\biggr)+\mathcal{P}_{r+1}(\kappa_{r+1})+
\frac{\mathsf
{E}D_{r+1}%
}{\kappa_{r+1}} \qquad\mbox{for }r=\bar{r}-1. \label{plMGineqQrbar}%
\end{eqnarray}
Then sum from $r=0$ to $\bar{r}-2$ to get
\[
\mathcal{Q}_{0}(\kappa_{0})=\mathcal{Q}_{\bar{r}-1}(
\kappa_{\bar{r}-1})+ 
{ \sum
_{r=0}^{\bar{r}-2}} 
\bigl\{\mathcal{Q}_{r}(
\kappa_{r})-\mathcal{Q}_{r+1}(\kappa_{r+1})\bigr\}
\]
and insert the bounds (\ref{plMGineqQr}), (\ref{plMGineqQrbar}) and
$\mathcal{P}_{\bar{r}}(\kappa_{\bar{r}})\leq\kappa_{\bar
{r}%
}^{-1}L\mathsf{E}D_{\bar{r}}$ from (\ref{plMGineqPrbar}).
\end{pf*}

\section{A metric on $\mathbb{R}$ and some inequalities}\label{schaining}

A metric is set up that will be used for the chaining argument. Then a number
of inequalities are shown, mostly related to this metric. Throughout
the rest
of this appendix, we denote by $C$ a constant which need not be the same in
different expressions.

Introduce the function
%
\begin{equation}
J_{i,p}(x,y)=(\varepsilon_{i}/\sigma)^{p}
\{1_{(\varepsilon_{i}\leq\sigma
y)}-1_{(\varepsilon_{i}\leq\sigma x)}\}, \label{J}%
\end{equation}
where $p\in\mathbb{N}_{0}$ and $\varepsilon_{i}/\sigma$ has density
$\mathsf{f}$. We will be interested in powers of $J_{i,p}(x,y)$ of order
$2^{r}$ where $r\in\mathbb{N}$ was chosen in Assumption~\ref{as}(i). Note
that $2^{r}p$ is even for $p\in\mathbb{N}_{0}$ and $r\in\mathbb{N}$ so that
$\varepsilon_{i}^{2^{r}p}$ is non-negative. Thus, define the increasing
function
\[
\mathsf{H}_{r}(x)=\int_{-\infty}^{x}
\bigl(1+u^{2^{r}p}\bigr)\mathsf{f}(u)\,\mathrm{d}u,
\]
with derivative $\mathsf{\dot{H}}_{r}(x)=(1+x^{2^{r}p})\mathsf{f}(x)$, along
with the constant
\[
H_{r}=\mathsf{H}_{r}(\infty)=\int_{-\infty}^{\infty}
\bigl(1+u^{2^{r}p}%
\bigr)\mathsf{f}(u)\,\mathrm{d}u<\infty.
\]
It follows that, for $x\leq y$ and $0\leq s\leq r$,%
%
\begin{equation}
0\leq\bigl|\mathsf{E}\bigl\{J_{i,p}(x,y)\bigr\}^{2^{s}}\bigr|\leq
\mathsf{E}\bigl\{\bigl|J_{i,p}%
(x,y)\bigr|^{2^{s}}\bigr\}\leq
\mathsf{H}_{r}(y)-\mathsf{H}_{r}(x), \label{EJ}%
\end{equation}
noting that, for $q\geq p\geq0$ and $\varepsilon\in\mathbb{R}$,
$|\varepsilon
^{p}|<1+|\varepsilon|^{q}$. We denote $\mathsf{H}_{r}(y)-\mathsf{H}_{r}(x)$
the $H_{r}$-distance between $y$ and $x$.

For the chaining, partition the range of $\mathsf{H}_{r}(c)$ into $K$
intervals of equal size $H_{r}/K$. That is, partition the support into $K$
intervals defined by the endpoints
%
\begin{equation}
-\infty=c_{0}<c_{1}<\cdots<c_{K-1}<c_{K}=
\infty, \label{xpartition}%
\end{equation}
and for $1\leq k\leq K$,
\[
\mathsf{E}\bigl[\bigl\{J_{i,p}(c_{k-1},c_{k})
\bigr\}^{2^{r}}\bigr]\leq\mathsf{H}_{r}%
(c_{k})-
\mathsf{H}_{r}(c_{k-1})=\frac{H_{r}}{K}.
\]
Let $c_{-k}=c_{0}$ for $k\in\mathbb{N}$.

The number of intervals $K$ will be chosen so large that $c_{-},c_{+}$ exist
which are (weakly) separated from zero by grid points in the sense that
$c_{k_{-}-1}\leq c_{-}\leq c_{k_{-}}\leq0$ and $0\leq c_{k_{+}-1}\leq
c_{+}\leq c_{k_{+}}$ and so that%
%
\begin{equation}
\mathsf{\dot{H}}_{r}(c_{-})=\mathsf{\dot
{H}}_{r}(c_{+})=H_{r}/\bigl(C_{\mathsf{H}%
}K^{1/2}
\bigr). \label{Hcutoff}%
\end{equation}
This can be done for sufficiently large $K$ since $\mathsf{f}$ is continuous
and since the function $\mathsf{\dot{H}}_{r}(c)=(1+c^{2^{r}p})\mathsf{f}(c)$
is integrable by Assumption~\ref{as}(i)(a).

The first inequality concerns the $H_{r}$-distance of additive perturbations
of the $]c_{k-1},c_{k}]$ intervals. It is used in the proof of the inequality
in Lemma~\ref{lineqJ}.

\begin{lemma}
\label{lineqFosci}Suppose Assumption~\ref{as}\textup{(i)} only holds for
$\nu=1$.
Then a constant $C>0$ exists so that for all $K$ satisfying \textup{(\ref
{Hcutoff})}
\[
\sup_{1\leq k\leq K}\sup_{|d|\leq K^{-1/2}}\bigl\{
\mathsf{H}_{r}(c_{k}%
+d)-\mathsf{H}_{r}(c_{k-1}+d)
\bigr\}\leq CH_{r}/K.
\]
\end{lemma}

\begin{pf}
1. \textitt{Definitions}. Consider positive
$c_{k}$ only, with a similar argument for negative $c_{k}$. Let
$\mathcal{H}=\mathsf{H}_{r}(c_{k}+d)-\mathsf{H}_{r}(c_{k-1}+d)$. Let
$\mathsf{\dot{H}}_{r}(c)=(1+c^{2^{r}p})\mathsf{f}(c)$ and
\[
\underline{\mathsf{\dot{H}}}_{r}(c)=\inf_{0\leq d\leq c}
\mathsf{\dot{H}} 
_{r}(d),\qquad \overline{\mathsf{
\dot{H}}}_{r}(c)=\sup_{d\geq c}\mathsf {\dot
{H}}_{r}(d),
\]
which are decreasing in $c$. Assumption~\ref{as}(i)(c) then implies%
%
\begin{equation}
C_{\mathsf{H}}^{-1}\overline{\mathsf{\dot{H}}}_{r}(c)\leq
\underline{\mathsf{\dot{H}}}_{r}(c)\leq\mathsf{\dot{H}}_{r}(c)
\leq \overline{\mathsf{\dot{H}}}_{r}(c)\leq C_{\mathsf{H}}
\underline{\mathsf {\dot {H}}}_{r}(c). \label{H'}%
\end{equation}
Since $\mathsf{\ddot{H}}_{r}(c)=2^{r}pc^{2^{r}p-1}\mathsf
{f}(c)+(1+c^{2^{r}%
p})\mathsf{\dot{f}}(c)$ then Assumption~\ref{as}(i)(b) gives%
%
\begin{equation}
\sup_{c\in\mathbb{R}}\bigl|\mathsf{\ddot{H}}_{r}(c)\bigr|<\infty.
\label{H''}%
\end{equation}
\begin{longlist}[1.]
\item[2.] \textitt{Apply the mean-value theorem} to get, for some $c_{\ell
}^{\ast}$ so
$c_{\ell-1}\leq c_{\ell}^{\ast}\leq c_{\ell}$, that%
%
\begin{equation}
H_{r}/K=\mathsf{H}_{r}(c_{\ell})-
\mathsf{H}_{r}(c_{\ell-1})=(c_{\ell}%
-c_{\ell-1})
\mathsf{\dot{H}}_{r}\bigl(c_{\ell}^{\ast}\bigr).
\label{mean-value}%
\end{equation}
Two inequalities for $\mathsf{\dot{H}}_{r}(c)$ arise from (\ref{H'}) and
condition (\ref{Hcutoff}). These are%
%
\begin{eqnarray}
\mathsf{\dot{H}}_{r}(c)&\leq&\overline{\mathsf{\dot{H}}}_{r}(c)
\leq \overline{\mathsf{\dot{H}}}_{r}(c_{+})\leq
C_{\mathsf{H}}\mathsf{\dot {H}}%
_{r}(c_{+})=H_{r}/K^{1/2}\qquad \mbox{for }c\geq c_{+}%
,\label{plineqFosciineqA}
\\
\label{plineqFosciineqB}%
\mathsf{\dot{H}}_{r}(c)&\geq&\underline{\mathsf{\dot{H}}}_{r}(c)
\geq \underline{\mathsf{\dot{H}}}_{r}(c_{+})\geq\overline{
\mathsf{\dot{H}}}%
_{r}(c_{+})/C_{\mathsf{H}}
\geq\mathsf{\dot{H}}_{r}(c_{+})/C_{\mathsf{H}}
\nonumber
\\[-8pt]
\\[-8pt]
\nonumber
&=&H_{r}/\bigl(C_{\mathsf{H}}^{2}K^{1/2}
\bigr) \qquad \mbox{for }0\leq c\leq c_{+}.
\end{eqnarray}
In parallel to (\ref{plineqFosciineqB}), which is derived for positive
$c$, it holds for negative $c$ that%
%
\begin{equation}
\mathsf{\dot{H}}_{r}(c)\geq H_{r}/\bigl(C_{\mathsf{H}}^{2}K^{1/2}
\bigr)\qquad \mbox{for }0\geq c\geq c_{-}. \label{plineqFosciineqD}%
\end{equation}

\item[3.] \textitt{Small arguments} $c_{-}\leq c_{k}^{\ast}\leq c_{+}$. Combine
(\ref{mean-value}), (\ref{plineqFosciineqB}) and
(\ref{plineqFosciineqD}) to get
%
\begin{equation}
c_{k}-c_{k-1}=H_{r}/\bigl\{K\mathsf{
\dot{H}}_{r}\bigl(c_{k}^{\ast}\bigr)\bigr\}\leq
C_{\mathsf{H}%
}^{2}/K^{1/2}. \label{Deltacbound}%
\end{equation}
Two second order Taylor expansions give
\begin{eqnarray*}
\mathsf{H}_{r}(c_{k}+d)-\mathsf{H}_{r}(c_{k})
& =&d\mathsf{\dot{H}}_{r}%
(c_{k})+
\bigl(d^{2}/2\bigr)\mathsf{\ddot{H}}_{r}\bigl(c_{k}^{\ast\ast}
\bigr),
\\
\mathsf{H}_{r}(c_{k-1}+d)-\mathsf{H}_{r}(c_{k-1})
& =&d\mathsf{\dot{H}}%
_{r}(c_{k-1})+
\bigl(d^{2}/2\bigr)\mathsf{\ddot{H}}_{r}\bigl(c_{k-1}^{\ast\ast}
\bigr),
\end{eqnarray*}
where $c_{k}^{\ast\ast},c_{k-1}^{\ast\ast}$ satisfy $\max(|c_{k}^{\ast
\ast
}-c_{k}|,|c_{k-1}^{\ast\ast}-c_{k-1}|)\leq|d|\leq K^{-1/2}$. The difference
is, when recalling the definition of $\mathcal{H}$ in item 1,
\[
\mathcal{H-\bigl\{}\mathsf{H}_{r}(c_{k})-
\mathsf{H}_{r}(c_{k-1})\bigr\}=d\bigl\{\mathsf {\dot
{H}}_{r}(c_{k})-\mathsf{\dot{H}}_{r}(c_{k-1})
\bigr\}+\bigl(d^{2}/2\bigr)\bigl\{\mathsf{\ddot {H}%
}_{r}
\bigl(c_{k}^{\ast\ast}\bigr)-\mathsf{\ddot{H}}_{r}
\bigl(c_{k-1}^{\ast\ast}\bigr)\bigr\}.
\]
The left-hand side is $\mathcal{H-}H_{r}/K$. The mean-value theorem
gives that
for some $\tilde{c}_{k}$, $c_{k-1}\leq\tilde{c}_{k}\leq c_{k}$, $\mathsf
{\dot
{H}}_{r}(c_{k})-\mathsf{\dot{H}}_{r}(c_{k-1})=(c_{k}-c_{k-1})\mathsf
{\ddot{H}%
}_{r}(\tilde{c}_{k})$. Insert this and rearrange to get
\[
0\leq\mathcal{H}=\frac{H_{r}}{K}+d(c_{k}-c_{k-1})
\mathsf{\ddot{H}}_{r}%
(\tilde{c}_{k})+
\frac{d^{2}}{2}\bigl\{\mathsf{\ddot{H}}_{r}\bigl(c_{k}^{\ast\ast
}
\bigr)-\mathsf{\ddot{H}}_{r}\bigl(c_{k-1}^{\ast\ast}
\bigr)\bigr\}.
\]
Using the bound $c_{k}-c_{k-1}\leq C_{\mathsf{H}}^{2}/K^{1/2}$ from
(\ref{Deltacbound}), and the bound $|d|\leq K^{-1/2}$, it follows that
$0\leq\mathcal{H}\leq C/K$, where $C=H_{r}+(C_{\mathsf{H}}^{2}+1)\sup_{c\in\mathbb{R}}|\mathsf{\ddot{H}}_{r}(c)|$ does not depend on $K$.

\item[4.] \textitt{Inequalities on tail grid point intervals.} Suppose
$c_{k}^{\ast
}\geq c_{+}$. This includes the situation where $c_{k}^{\ast}$ and
$c_{+}$ are
in the same grid interval. Expansion (\ref{mean-value}) and inequality
(\ref{plineqFosciineqA}) imply
\[
c_{k}-c_{k-1}=H_{r}/\bigl\{K\mathsf{
\dot{H}}_{r}\bigl(c_{k}^{\ast}\bigr)\bigr\}\geq
H_{r}%
/\bigl\{KH_{r}/K^{1/2}\bigr
\}=K^{-1/2}\geq|d|.
\]

\item[5.] \textitt{Large arguments} $c_{k}^{\ast}\geq c_{+}$ so either $k\geq k_{+}+2$
or $k=k_{+}+1$ with $c_{k-1}^{\ast}\geq c_{+}$. In this case
$c_{k-1}^{\ast
}\geq c_{+}$ so that item 4 shows that $c_{k+1}-c_{k}$, $c_{k}-c_{k-1}$ and
$c_{k-1}-c_{k-2}$ are all larger than $|d|$. Therefore,
\begin{eqnarray*}
c_{k}+d & \leq& c_{k}+|d|\leq c_{k}+c_{k+1}-c_{k}=c_{k+1},
\\
c_{k-1}+d & \geq& c_{k-1}-|d|\geq c_{k-1}-(c_{k-1}-c_{k-2})=c_{k-2}.
\end{eqnarray*}
It then holds that $0\leq\mathcal{H}\leq\mathsf{H}_{r}(c_{k+1})-\mathsf
{H}%
_{r}(c_{k-2})=3H_{r}/K$.

\item[6.] \textitt{Intermediate arguments} $c_{k}^{\ast}\geq c_{+}$ so that $k=k_{+}$.
Item 4 shows $c_{k}-c_{k-1}\geq|d|$ and $c_{k+1}-c_{k}\geq|d|$. 
Thus, for $d>0,$
$0\leq\mathcal{H\leq}\mathsf{H}_{r}(c_{k+1})-\mathsf{H}_{r}(c_{k-1}%
)=2H_{r}/K.$ For $d<0,$ write $\mathcal{H=H}_{1}+\mathcal{H}_{2}$ where
$\mathcal{H}_{1}=\mathsf{H}_{r}(c_{k}+d)-\mathsf{H}_{r}(c_{+})$ and
$\mathcal{H}_{2}=\mathsf{H}_{r}(c_{+})-\mathsf{H}_{r}(c_{k-1}+d).$ Again,
$0\leq\mathcal{H}_{1}\leq\mathsf{H}_{r}(c_{k+1})-\mathsf{H}_{r}(c_{k-1}%
)=2H_{r}/K.$ For $\mathcal{H}_{2}$ use the mean-value theorem to get
\begin{equation}
\mathcal{H}_{2}=(c_{+}-c_{k-1}-d)\mathsf{\dot{H}}_{r}(c_{+})-\tfrac{1}{2}%
(c_{+}-c_{k-1}-d)^{2}\mathsf{\ddot{H}}_{r}\bigl(c_{+}^{\ast}\bigr) \label{new}%
\end{equation}
for an intermediate point $c_{k-1}+d\leq c_{+}^{\ast}\leq c_{+}$. Now, argue
as in (\ref{Deltacbound}) in item 3 to get $c_{+}-c_{k-1}\leq C_{\mathsf{H}}^{2}/K^{1/2}.$
Since $|d|\leq K^{-1/2}$ while $\mathsf{\dot{H}}_{r}(c_{+})=H_{r}%
/(C_{\mathsf{H}}K^{1/2}),$ see (\ref{Hcutoff}), then the first term in (\ref{new}) is of
order $K^{-1}$ uniformly in $k.$ Similarly, the second term in (\ref{new}) is
of the same order since $\mathsf{\ddot{H}}_{r}$ is bounded by (\ref{H''}).

\item[7.] \textitt{Intermediate arguments} $c_{k}^{\ast}\geq c_{+}$ so that
$k=k_{+}+1$ with $c_{k-1}^{\ast}<c_{+}$ and $c_{k-1}+d\geq c_{+}$. Decompose
$0\leq\mathcal{H}\leq\mathcal{H}_{1}+\mathcal{H}_{2}$ where
\[
\mathcal{H}_{1}=\mathsf{H}_{r}(c_{k}+d)-
\mathsf{H}_{r}(c_{k-1}),\qquad \mathcal{H}_{2}=
\mathsf{H}_{r}(c_{k-1})-\mathsf{H}_{r}(c_{+}).
\]
Consider $\mathcal{H}_{1}$. Argue $\mathcal{H}_{1}\leq2H_{r}/K$ as in item
5.

Consider $\mathcal{H}_{2}$. Argue $c_{k-1}-c_{+}\geq|d|$ as
in item
4 and in turn $\mathcal{H}_{2}\leq H_{r}/K$ as in item 5.

\item[8.] \textitt{Intermediate arguments} $c_{k}^{\ast}\geq c_{+}$ so $k=k_{+}+1$
with $c_{k-1}^{\ast}<c_{+}$ and $c_{k-1}+d<c_{+}$. Decompose $0\leq
\mathcal{H}=\mathcal{H}_{1}+\mathcal{H}_{2}+\mathcal{H}_{3}$ where
$\mathcal{H}_{1}$ and $\mathcal{H}_{2}$ were defined and analyzed in
item 6,
while
\[
\mathcal{H}_{3}=\mathsf{H}_{r}(c_{+})-
\mathsf{H}_{r}(c_{k-1}+d).
\]
Since $c_{+}\leq c_{k_{+}}=c_{k-1}$ and $c_{k-1}+d<c_{+}$ then $c_{k-1}%
-c_{+}\leq|d|$. The mean-value theorem shows
\[
\mathcal{H}_{3}=\delta_{k,d}\mathsf{\dot{H}}_{r}(c_{+})+
\bigl(\delta_{k,d}%
^{2}/2\bigr)\mathsf{
\ddot{H}}_{r}\bigl(c^{\ast\ast}\bigr),
\]
where $\delta_{k,d}=c_{k-1}+d-c_{+}$ while $c^{\ast\ast}$ satisfies
$|c^{\ast\ast}-c_{+}|\leq|\delta_{k,d}|$. Here $|\delta_{k,d}|\leq
c_{k-1}-c_{+}+|d|\leq2|d|\leq2K^{-1/2}$. Because (\ref{Hcutoff}) shows
$\mathsf{\dot{H}}_{r}(c_{+})=H_{r}/(C_{\mathsf{H}}K^{1/2})$, while
$\mathsf{\ddot{H}}_{r}(c^{\ast\ast})$ is bounded by (\ref{H''}), it follows
that $\mathcal{H}_{3}\leq C/K$.
\end{longlist}
\end{pf}

The next lemma shows how small fluctuations in the arguments of the function
$J_{i,p}$ can be controlled in terms of $J_{i,p}$ functions defined on the
grid points. The results are used in the proofs of Theorems \ref{tFpest}, \ref{tFpestomega},
that are
concerned with estimation error $b$ in the empirical process $\mathbb{F}
_{n}^{q,p}(b,c)$. The proof uses Lemma~\ref{lineqFosci}.

\begin{lemma}
\label{lineqJ}Suppose Assumption~\ref{as}\textup{(i)} only holds for $\nu=1$. For
any $c\leq c_{K-1}$ we choose grid points, see \textup{(\ref{xpartition})},
$c_{k-1}<c\leq c_{k}(\leq c_{K-1})$. For $c>c_{K-1}$ we consider
$c_{K-1}<c<c_{K}(=\infty)$. Then an integer $k_{J}>0$ exists such that, for
all $K$ satisfying \textup{(\ref{Hcutoff})} and all $c,d,d_{m}\in\mathbb{R}$
for which
$|d|\leq K^{-1/2}$ and $|d-d_{m}|\leq K^{-1}$, integers $k^{\dag
},k^{\ddag}$
exist for which
\[
\bigl|J_{i,p}(c,c+d)-J_{i,p}(c_{k},c_{k}+d_{m})\bigr|
\leq\bigl|J_{i,p}(c_{k-k_{J}}%
,c_{k})\bigr|+\bigl|J_{i,p}(c_{k^{\dag}-k_{J}},c_{k^{\dag}})\bigr|
+\bigl|J_{i,p}(c_{k^{\ddag
}-k_{J}},c_{k^{\ddag}})\bigr|.
\]
\end{lemma}

\begin{pf}
1. \textitt{Decomposition}. Only the case $k<K$
is proved. The proof for $k=K$ is similar. Let $\sigma=1$ for notational
simplicity. Write
\[
\mathcal{J}=J_{i,p}(c,c+d)-J_{i,p}(c_{k},c_{k}+d_{m})=
\varepsilon_{i}%
^{p}(\mathcal{I}_{1}+
\mathcal{I}_{2}-\mathcal{I}_{3}),
\]
in terms of indicator functions $\mathcal{I}_{1}=1_{(c<\varepsilon
_{i}\leq
c_{k})}$, $\mathcal{I}_{2}=1_{(\varepsilon_{i}\leq
c_{k}+d)}-1_{(\varepsilon
_{i}\leq c_{k}+d_{m})}$ and $\mathcal{I}_{3}=1_{(c+d<\varepsilon_{i}\leq
c_{k}+d)}$. It follows that $|\mathcal{J}|\leq|\varepsilon_{i}^{p}%
|(\mathcal{I}_{1}+|\mathcal{I}_{2}|+\mathcal{I}_{3})$.

2. \textitt{Bound for} $\mathcal{I}_{1}$. Since $c_{k-1}<c\leq c_{k}$ then
$0\leq\mathcal{I}_{1}=1_{(c<\varepsilon_{i}\leq c_{k})}\leq1_{(c_{k-1}%
<\varepsilon_{i}\leq c_{k})}$.

3. \textitt{Bound for} $\mathcal{I}_{2}$. Write $d=d_{m}+(d-d_{m})$ where
$|d-d_{m}|\leq K^{-1}$. Let $c^{\dag}=c_{k}+d_{m}$. Then $|\mathcal{I}%
_{2}|\leq1_{(c^{\dag}-K^{-1}\leq\varepsilon_{i}\leq c^{\dag}+K^{-1})}$. Using
first this inequality and then the mean-value theorem, it follows that
\[
\mathcal{E}_{2}=\mathsf{E}\bigl(\bigl|\varepsilon_{i}^{p}
\mathcal{I}_{2}\bigr|\bigr)\leq \mathsf{H}_{r}
\bigl(c^{\dag}+K^{-1}\bigr)-\mathsf{H}_{r}
\bigl(c^{\dag}-K^{-1}\bigr)\leq 2H_{r}%
^{-1}
\sup_{c\in\mathbb{R}}\mathsf{\dot{H}}_{r}(c)H_{r}/K.
\]
Therefore, a $k^{\dag}$ exists for which $|\mathcal{I}_{2}|\leq
1_{(c_{k^{\dag
}-k_{J}}<\varepsilon_{i}\leq c_{k^{\dag}})}$, where $k_{J}\leq
2H_{r}^{-1}%
\sup_{c\in\mathbb{R}}\mathsf{\dot{H}}_{r}(c)+2$.

4. \textitt{Bound for} $\mathcal{I}_{3}$. Because $c_{k-1}<c\leq c_{k}$, then
$\mathcal{I}_{3}\leq1_{(c_{k-1}+d<\varepsilon_{i}\leq c_{k}+d)}$. Using first
this inequality and then Lemma~\ref{lineqFosci} and noting that
$|d|\leq
K^{-1/2}$, we find
\[
\mathcal{E}_{3}=\mathsf{E}\bigl(\bigl|\varepsilon_{i}^{p}\bigr|
\mathcal{I}_{3}\bigr)\leq \mathsf{H}_{r}(c_{k}+d)-
\mathsf{H}_{r}(c_{k-1}+d)\leq CH_{r}/K.
\]
Therefore, a $k^{\ddag}$ exists for which $|\mathcal{I}_{3}|\leq
1_{(c_{k^{\ddag}-k_{J}}<\varepsilon_{i}\leq c_{k^{\ddag}})}$ where
$k_{J}\leq
C+1$.
\end{pf}

The next inequality gives a tightness type result for the function
$\mathsf{H}_{r}$. This lemma is used in the proof of the tightness
result for
the empirical process $\mathbb{F}_{n}^{g,p}(0,c)$ in Theorem~\ref{tFptight}.

\begin{lemma}
\label{lFtoHscale}Let $c_{\psi}=\mathsf{F}^{-1}(\psi)$. For all densities
satisfying Assumption~\ref{as}\textup{(i)(a)} for some $\nu<1$, there exist
$C_{\nu
},\phi_{0}>0$ such that for all $0\leq\phi\leq\phi_{0}$ it follows that
\[
\max_{0\leq\psi\leq1-\phi}\bigl\{\mathsf{H}_{r}(c_{\psi+\phi})-
\mathsf{H}%
_{r}(c_{\psi})\bigr\}\leq
C_{\nu}\phi^{1-\nu}.
\]
\end{lemma}

\begin{pf}
Let $\psi_{0}=\mathsf{F}(0)$. Note that
$2^{r}p$ is even for $r\in\mathbb{N}$, $p\in\mathbb{N}_{0}$.
\begin{longlist}[1.]
\item[1.] \textitt{Let} $\psi\geq\psi_{0}$. Then $\mathsf{H}_{r}(c_{\psi+\phi
})-\mathsf{H}_{r}(c_{\psi})$ is increasing in $\psi$ since, with $\dot
{c}_{\psi}=1/\mathsf{f}(c_{\psi})$,
\[
\frac{\mathrm{d}}{\mathrm{d}\psi}\bigl\{\mathsf{H}_{r}(c_{\psi+\phi})-
\mathsf{H}_{r}(c_{\psi
})\bigr\}=\frac{\mathsf{\dot{H}}_{r}(c_{\psi+\phi})}{\mathsf{f}(c_{\psi+\phi}
)}-
\frac{\mathsf{\dot{H}}_{r}(c_{\psi})}{\mathsf{f}(c_{\psi})}=c_{\psi
+\phi
}^{p2^{r}}-c_{\psi}^{p2^{r}}>0.
\]
Thus, $\max_{\psi_{0}\leq\psi\leq1-\phi}\{\mathsf{H}_{r}(c_{\psi+\phi
})-\mathsf{H}_{r}(c_{\psi})\}\leq\mathsf{H}_{r}(\infty)-\mathsf{H}%
_{r}(c_{1-\phi})$. This bound satisfies
\[
\mathsf{H}_{r}(\infty)-\mathsf{H}_{r}(c_{1-\phi})=
\int_{c_{1-\phi
}}^{\infty
}\bigl(1+u^{p2^{r}}\bigr)
\mathsf{f}(u)\,\mathrm{d}u=\phi+\int_{c_{1-\phi}}^{\infty}u^{p2^{r}%
}
\mathsf{f}(u)\,\mathrm{d}u.
\]
Assumption~\ref{as}(i)(a) shows $\mathsf{E}\varepsilon^{p2^{r}/\nu}\leq
C$ for
some $C>0$ so that $1-\mathsf{F}(u)\leq Cu^{-p2^{r}/\nu}$ by the Chebychev
inequality. Hence, $u^{p2^{r}}\leq C^{\nu}\{1-\mathsf{F}(u)\}^{-\nu}$,
so that
\[
\mathsf{H}_{r}(\infty)-\mathsf{H}_{r}(c_{1-\phi})
\leq\phi+C^{\nu}%
\int_{c_{1-\phi}}^{\infty}
\bigl\{1-\mathsf{F}(u)\bigr\}^{-\nu}\mathsf{f}(u)\,\mathrm{d}u.
\]
Substituting $x=\mathsf{F}(u)$, so that $\mathrm{d}x=\mathsf{f}(u)\,\mathrm{d}u$ gives
\[
\mathsf{H}_{r}(\infty)-\mathsf{H}_{r}(c_{1-\phi})
\leq\phi+C^{\nu}\int_{1-\phi
}^{1}(1-x)^{-\nu}\,\mathrm{d}x=
\phi+\frac{C^{\nu}}{1-\nu}\phi^{1-\nu}.
\]

\item[2.] \textitt{Let} $\psi\leq\psi_{0}-\phi$. Apply a similar argument as in item
1, to show that $\mathsf{H}_{r}(c_{\psi+\phi})-\mathsf{H}_{r}(c_{\psi
})$ is
decreasing because $c_{\psi}<c_{\psi+\phi}\leq0$. Thus, $\mathsf{H}%
_{r}(c_{\phi})-\mathsf{H}_{r}(-\infty)$ satisfies the same bound.

\item[3.] \textitt{Let} $\psi_{0}-\phi\leq\psi\leq\psi_{0}$. Then
\[
\mathcal{H=}\max_{\psi_{0}-\phi\leq\psi\leq\psi_{0}}\bigl\{\mathsf{H}_{r}%
(c_{\psi+\phi})-
\mathsf{H}_{r}(c_{\psi})\bigr\}\leq\mathsf{H}_{r}(c_{\psi
_{0}+\phi
})-
\mathsf{H}_{r}(c_{\psi_{0}-\phi}).
\]
Using the mean-value theorem there exists a $\psi^{\ast}$, in the interval
$\psi_{0}-\phi\leq\psi^{\ast}\leq\psi_{0}+\phi$, for which
\[
\mathcal{H\leq}\frac{\mathsf{\dot{H}}_{r}(c_{\psi
^{\ast
}})}{\mathsf{f}(c_{\psi^{\ast}})}2\phi=2\bigl(1+c_{\psi^{\ast
}}^{2^{r}p}
\bigr)\phi \leq2\bigl\{1+\max\bigl(c_{\psi_{0}-\phi_{0}}^{2^{r}p},c_{\psi_{0}+\phi
_{0}}^{2^{r}%
p}
\bigr)\bigr\}\leq C\phi,
\]
for some $C>0$, because $\phi_{0}$ can be chosen so that the two
quantiles are finite.

\item[4.] \textitt{Combine results}. Note that $\phi\leq\phi^{1-\nu}$. Let
$C_{\nu
}=\max\{2H_{r},1+C^{\nu}/(1-\nu)\}$.\quad\qed
\end{longlist}
\noqed\end{pf}

\section{Proofs of auxiliary Theorems 
\texorpdfstring{\protect\ref{tFpest}}{4.1}--\texorpdfstring{\protect\ref{tFptight}}{4.4}}
\label{sFproof}

\begin{pf*}{Proof of Theorem~\ref{tFpest}}
Without loss of generality, let $\sigma=1$.
Let $\tilde{R}(b,c_{\psi})=\mathbb{F}_{n}^{g,p}(b,c_{\psi})-\mathbb{F}%
_{n}^{g,p}(0,c_{\psi})$ and $\mathcal{R}_{n}=\sup_{0\leq\psi\leq1}%
\sup_{|b|\leq n^{1/4-\eta}B}|\mathbb{F}_{n}^{g,p}(b,c_{\psi})-\mathbb{F}
_{n}^{g,p}(0,c_{\psi})|$.

1. \textitt{Partition the support}. For $\delta,n>0$, partition the axis as
laid out in (\ref{xpartition}) with $K=\operatorname{int}(H_{r}n^{1/2}/\delta)$
using Assumption~\ref{as}(i)(a) with $\nu=1$ only.

2. \textitt{Assign} $c_{\psi}$ \textitt{to the partitioned support}. Consider
$0\leq\psi\leq1$. Thus, for each $c_{\psi}$ there exists
$c_{k-1},c_{k}$ so $c_{k-1}<c_{\psi}\leq c_{k}$.

3. \textitt{Construct} $b$-\textitt{balls}. For a $\zeta>\kappa$, cover
the set
$|b|\leq n^{1/4-\eta}B$ with $M$ balls of radius $n^{-\zeta}$ with centers
$b_{m}$, that is $M=\mathrm{O}\{n^{(1/4-\eta+\zeta)\dim x}\}$. Thus,
for any
$b$ there exists a $b_{m}$ so that $|b-b_{m}|<n^{-\zeta}$.

4. \textitt{Apply chaining}. For $k<K$ where $c_{\psi}\leq c_{k}\leq c_{K-1}$,
we compare $c_{\psi}$ to the nearest right grid point, $c_{k}$, using
$\tilde{R}(b,c_{\psi})=\tilde{R}(b_{m},c_{k})+\{\tilde{R}(b,c_{\psi}%
)-\tilde{R}(b_{m},c_{k})\}$, whereas for $k=K$, we use the nearest left grid
point, $c_{K-1}$, and get $\tilde{R}(b,c_{\psi})=\tilde{R}(b_{m}%
,c_{K-1})+\{\tilde{R}(b,c_{\psi})-\tilde{R}(b_{m},c_{K-1})\}$. Therefore
$R_{n}\leq R_{n,1}+R_{n,2}$, where
\begin{eqnarray*}
\mathcal{R}_{n,1} & =&\max_{1\leq k<K}\max
_{1\leq m\leq M}\bigl|\tilde{R}%
(b_{m},c_{k})\bigr|,
\\
\mathcal{R}_{n,2} & =&\max_{1\leq k<K}\max
_{1\leq m\leq M}\sup_{c_{k-1}%
<c_{\psi}\leq c_{k}}\sup_{|b-b_{m}|<n^{-\zeta}}\bigl|
\tilde{R}(b,c_{\psi}%
)-\tilde{R}(b_{m},c_{k})\bigr|
\\
&&{} +\max_{1\leq m\leq M}\sup_{c_{K-1}<c_{\psi}}\sup
_{|b-b_{m}|<n^{-\zeta}
}\bigl|\tilde{R}(b,c_{\psi})-\tilde{R}(b_{m},c_{K-1})\bigr|.
\end{eqnarray*}
Thus, it suffices to show that $\mathsf{P}(\mathcal{R}_{n,j}>\gamma)$ vanishes
for $j=1,2$.

5. \textitt{The term} $\mathcal{R}_{n,1}$. Use Theorem~\ref{lsupMrasymp} to see that $\mathcal{R}_{n,1}=\mathrm{o}_{\mathsf{P}}(1)$.
To see this, let $\upsilon=1/2$ and let $g_{in}$ have coordinates
$g_{in}^{\ast}$. Then, for $z_{\ell i}=g_{in}^{\ast}J_{i,p}(c_{k},c_{k}%
+\sigma^{-1}x_{in}^{\prime}b_{m})$ we write the coordinates of $\tilde
{R}(b_{m},c_{k})$ as $n^{-1/2}\sum_{i=1}^{n}(z_{\ell i}-\mathsf{E}%
_{i-1}z_{\ell i})$, see definition in (\ref{J}), and where $\ell$ represents
the indices $k,m$. The conditions of Theorem~\ref{lsupMrasymp} need to
be verified.

\textitt{The parameter} $\lambda$. The set of indices $\ell$ has size
$L=\mathrm{O}(n^{\lambda})$ where $\lambda=1/2+(1/4-\eta+\zeta)\dim x$ since
$K=\mathrm{O}(n^{1/2})$ and $M=\mathrm{O}\{n^{(1/4-\eta+\zeta)\dim b}\}$.

\textitt{The parameter} $\varsigma$. Because $|1_{(\varepsilon_{i}\leq
c_{k}+x_{in}^{\prime}b_{m})}-1_{(\varepsilon_{i}\leq c_{k})}|\leq
1_{(c_{k}-|x_{in}||b_{m}|<\varepsilon_{i}\leq c_{k}+|x_{in}||b_{m}|)}$
we find
for $1\leq q\leq r$, that
\[
\mathsf{E}_{i-1}(J_{i,p})^{2^{q}}\leq
\mathsf{H}_{r}\bigl(c_{k}+|x_{in}%
||b_{m}|\bigr)-
\mathsf{H}_{r}\bigl(c_{k}-|x_{in}||b_{m}|\bigr)
\leq2|x_{in}||b_{m}|\sup_{v\in\mathbb{R}}\mathsf{
\dot{H}}_{r}(v),
\]
using the mean-value theorem. Because $|b_{m}|\leq n^{1/4-\eta}B$, while
$\sup_{v\in\mathbb{R}}\mathsf{\dot{H}}_{r}(v)<\infty$ by Assumption~\ref{as}(i)(b), we find
%
\begin{equation}
D_{q}=\max_{1\leq\ell\leq L}%
{
\sum_{i=1}^{n}} 
\mathsf{E}_{i-1}(z_{\ell i})^{2^{q}}\leq C_{1}
\Biggl(n^{-1}%
{ \sum
_{i=1}^{n}} 
\bigl|g_{in}^{\ast}\bigr|^{2^{q}}\bigl|n^{1/2}x_{in}\bigr|
\Biggr)n^{3/4-\eta}. \label
{plFptildebD}%
\end{equation}
Thus, $\mathsf{E}D_{q}=\mathrm{O}(n^{\varsigma})$ where $\varsigma
=3/4-\eta$
by Assumption~\ref{as}(iii)(a).

\textitt{Condition} (i) is that $\varsigma<2\upsilon$. This holds since
$0<\eta$ so that $\varsigma=3/4-\eta<1=2\upsilon$.

\textitt{Condition} (ii) is that $\varsigma+\lambda<\upsilon
2^{\bar{r}}$
with $\bar{r}=r$. If $\zeta>\kappa$ is chosen sufficiently small,
then
\[
\varsigma+\lambda=1+(1/4+\kappa-\eta) (1+\dim x)+(\zeta-\kappa)\dim x-\kappa<
\upsilon2^{r}=2^{r-1},
\]
provided $r$ is chosen so that $2^{r}-1\geq1+(1/4+\kappa-\eta)(1+\dim x)$.

6. \textitt{Decompose} $\mathcal{R}_{n,2}$. It will be argued that
$\mathcal{R}_{n,2}\leq3(\tilde{\mathcal{R}}_{n,2}+2\overline
{\mathcal{R}%
}_{n,2})+\mathrm{o}_{\mathsf{P}}(1)$, where%
%
\begin{eqnarray}
\tilde{\mathcal{R}}_{n,2} & =&\max_{1\leq k\leq K}n^{-1/2}%
{
\sum_{i=1}^{n}} 
|g_{in}|
\bigl\{\bigl|J_{i,p}(c_{k-k_{J}},c_{k})\bigr|-\mathsf
{E}_{i-1}\bigl|J_{i,p}(c_{k-k_{J}%
},c_{k})\bigr|\bigr
\},\label{plempproclocVtilde}
\\
\overline{\mathcal{R}}_{n,2} & =&\max_{1\leq k\leq K}n^{-1/2}%
{
\sum_{i=1}^{n}} 
|g_{in}|
\mathsf{E}_{i-1}\bigl|J_{i,p}(c_{k-k_{J}},c_{k})\bigr|.
\label{plempproclocVbar}%
\end{eqnarray}
To see this, let $c_{k}$ denote the nearest right grid point for
$c_{\psi}\leq
c_{K-1}$ while $c_{k}=c_{K-1}$ for $c_{\psi}>c_{K-1}$. Note first that
$\tilde{R}_{\mathsf{F}}^{p}(b,c_{\psi})-\tilde{R}_{\mathsf{F}}^{p}(b_{m}
,c_{k})$ involves the functions
\[
\mathcal{J}_{i}=J_{i,p}\bigl(c_{\psi},c_{\psi}+x_{in}^{\prime}b
\bigr)-J_{i,p}%
\bigl(c_{k},c_{k}+x_{in}^{\prime}b_{m}
\bigr).
\]
Assumption~\ref{as}(ii) gives that $\max_{1\leq i\leq n}|x_{in}%
|=\mathrm{O}_{\mathsf{P}}(n^{\kappa-1/2})$. Thus, for all $\epsilon>0$ a
$C_{x}>0$ exists so that the set $(\max_{1\leq i\leq n}|x_{in}|\leq
n^{\kappa-1/2}C_{x})$ has probability of at least $1-\epsilon$. On that set,
using $d=x_{in}^{\prime}b$ and $d_{m}=x_{in}^{\prime}b_{m}$,
$|d|=\mathrm{O}%
(n^{-1/4+\kappa-\eta})=\mathrm{o}(K^{-1/2})$ for $\eta-\kappa>0$ and
$|d-d_{m}|=\mathrm{O}(n^{-1/2+\kappa-\zeta})=\mathrm{o}(K^{-1})$ for
$\zeta-\kappa>0$. Thus, for sufficiently large $n,|d|<K^{-1/2}$ and
$|d-d_{m}|<K^{-1}$. Lemma~\ref{lineqJ} using Assumption~\ref{as}(i) then
shows that a $k_{J}$ exists so that, for all $c,d,d_{m}$, there exist
$k^{\dag},k^{\ddag}$ for which%
%
\begin{equation}
|\mathcal{J}_{i}|\leq\bigl|J_{i,p}(c_{k-k_{J}},c_{k})\bigr|+\bigl|J_{i,p}(c_{k^{\dag
}-k_{J}%
},c_{k^{\dag}})\bigr|+\bigl|J_{i,p}(c_{k^{\ddag}-k_{J}},c_{k^{\ddag}})\bigr|.
\label{plFptildebJi}%
\end{equation}
As a consequence it holds, as desired, that $\mathcal{R}_{n,2}\leq
3(\tilde{\mathcal{R}}_{n,2}+2\overline{\mathcal{R}}_{n,2})+\mathrm
{o}%
_{\mathsf{P}}(1)$.

7. \textitt{The term} $\tilde{\mathcal{R}}_{n,2}$ is $\mathrm{o}%
_{\mathsf{P}}(1)$ by Lemma~\ref{lsupMrasymp}. Let $\upsilon=1/2$. To see
this, note that $\tilde{\mathcal{R}}_{n,2}$ is the maximum of a
family of
martingales of the required form with $\ell=k$ so that $L=K$ and
$z_{\ell
i}=|g_{in}||J_{i,p}(c_{k-k_{J}},c_{k})|$ and it suffices to set $\bar{r}=2$.

Condition (i) holds with $\lambda=1/2$ since $K=\operatorname{int}(H_{r}%
n^{1/2}/\delta)$.

Condition (ii) holds with $\varsigma=1/2$ since $\mathsf{E}_{i-1}%
(J_{i,p})^{2^{\bar{r}}}\leq\mathsf{H}_{r}(c_{k})-\mathsf
{H}_{r}(c_{k-k_{J}%
})=k_{J}H_{r}/K$ for $r\geq\bar{r}=2$. Thus $\sum_{i=1}^{n}\mathsf{E}%
_{i-1}(J_{i,p})^{2^{\bar{r}}}=\mathrm{O}(n^{1-1/2})$, uniformly in $\ell,i$.

It holds that $\lambda+\varsigma=1$ which is less than $\upsilon2^{\bar
{r}%
}=2$.

8. \textitt{Bounding} $\mathcal{\bar{R}}_{n,2}$. Note $\mathsf{E}%
_{i-1}|J_{i,p}(c_{k-k_{J}},c_{k})|\leq k_{J}H_{r}/K\leq2k_{J}\delta n^{-1/2}$
uniformly in $i,k$ by the same argument as in item 7 and since
$K=\operatorname{int}%
(H_{r}n^{1/2}/\delta)$. It follows that $\mathcal{\bar{R}}_{n,2}\leq
2k_{J}\delta n^{-1}\sum_{i=1}^{n}|g_{in}|$. Here $n^{-1}\sum_{i=1}^{n}%
|g_{in}|=\mathrm{O}_{\mathsf{P}}(1)$ by Markov's inequality and Assumption~\ref{as}(iii)(a), so that $\mathcal{\bar{R}}_{n,2}=\mathrm{O}_{\mathsf
{P}%
}(\delta)$. Thus, choosing $\delta$ sufficiently small, $\mathcal{\bar
{R}%
}_{n,2}$ is small in probability.
\end{pf*}

\begin{pf*}{Proof of Theorem~\ref{tFpestomega}}
It suffices to show, for all
$\omega<\eta-\kappa$ where $\eta-\kappa\leq1/4$, that
\begin{eqnarray*}
\mathcal{S}_{1}&=&\sup_{0\leq\psi\leq1}\sup
_{|b|\leq n^{1/4+\kappa-\eta}B} 
\sup_{d\in\mathbb{R}}\bigl|
\mathbb{F}_{n}^{1,0}\bigl(b,c_{\psi}+n^{\kappa-1/2}%
d
\bigr)-\mathbb{F}_{n}^{1,0}\bigl(0,c_{\psi}+n^{\kappa-1/2}d
\bigr)\bigr|  =\mathrm{o}%
_{\mathsf{P}}\bigl(n^{-\omega}\bigr),
\\
\mathcal{S}_{2}&=&\sup_{0\leq\psi\leq1}\sup_{|d|\leq n^{1/4+\kappa-\eta}%
B}\bigl|
\mathbb{F}_{n}^{1,0}\bigl(0,c_{\psi}+n^{\kappa-1/2}d
\bigr)-\mathbb{F}_{n}%
^{1,0}(0,c_{\psi})\bigr|
=\mathrm{o}_{\mathsf{P}}\bigl(n^{-\omega}\bigr).
\end{eqnarray*}
For each term the proof of Theorem~\ref{tFpest} is used with minor
modifications. Since $p=0$ then $2^{r}p=0$ for all $r$, which
simplifies the
assumptions, see Remark~\ref{ras}(a). Moreover, when using Theorem~\ref{lsupMrasymp}, $z_{\ell,i}^{2^{r}}=z_{\ell,i}$ for all $r\geq1$.
Thus, it
suffices to check $\varsigma<2\upsilon$ and $\lambda<\infty$.

\textitt{A. The term} $\mathcal{S}_{1}$. The steps of the proof of Theorem~\ref{tFpest} are modified as follows.

1. Choose $K=\operatorname{int}(H_{r}n^{1/2+1/8+\omega/2}/\delta)$ where
$\omega<\eta-\kappa\leq1/4$.

2. For each $c_{\psi}+n^{\kappa-1/2}d$, there exist $c_{k-1},c_{k}$ depending
on $n$ so that $c_{k-1}<c_{\psi}+n^{\kappa-1/2}d\leq c_{k}$.

3. Choose $\zeta\geq\eta$ which implies $\zeta>\kappa$ since $\kappa
<\eta$.
The $b$-set is now $|b|\leq n^{1/4+\kappa-\eta}B$ so that the number of
$b$-balls is $M=\mathrm{O}\{n^{(1/4+\kappa-\eta+\zeta)\dim x}\}$.

4. Note that in the chaining argument, $c_{\psi}$ is replaced by
$c_{\psi
}+n^{\kappa-1/2}d$. This only affects $\mathcal{R}_{n,2}$.

5. The term $\mathcal{R}_{n,1}$ is $\mathrm{o}_{\mathsf
{P}}(n^{-1/8-\omega
/2})$. Use Theorem~\ref{lsupMrasymp} with $\upsilon=3/8-\omega
/2>1/2+\kappa
-\eta$. Define $z_{\ell i}$ as before. Since $p=0$, $g_{i}n=1$ then
$|J_{i,p}(x,y)|^{2^{r}}=|J_{i,p}(x,y)|$ and $|z_{\ell i}|=|z_{\ell i}{}%
^{2^{r}}|$ for any $r\in\mathbb{N}_{0}$. The inequality
(\ref{plFptildebD}) for $D_{q}$ holds as before, uniformly in
$q\in\mathbb{N}$ so $\varsigma=3/4+\kappa-\eta$. Thus, \textitt{condition}
(i) holds since $\varsigma=3/4+\kappa-\eta<3/4-\omega=2\upsilon$. Moreover,
$\lambda=1/2+\omega+(1/4+\kappa-\eta+\zeta)\dim x$ is finite so
\textitt{condition} (ii) holds for some $r$.

6. Lemma~\ref{lineqJ} is an analytic result holding in finite samples.
So the
argument is not affected by the dependence of $c_{k}$ on $n$ through
$c_{\psi
}+n^{\kappa-1/2}d$. In particular, (\ref{plFptildebJi}) holds as stated
and therefore the decomposition of $\mathcal{R}_{n,2}$ holds, noting
that $K$
is now chosen differently.

7. The term $\mathcal{\tilde{R}}_{n,2}$ is $\mathrm{o}_{\mathsf{P}}%
(n^{-1/4})$. Use Theorem~\ref{lsupMrasymp} with some $\upsilon
>3/16-\omega
/4$. Here $\lambda=5/8+\omega/2<\infty$ by the definition of $K$, while
$\varsigma=3/8-\omega/2$ since $\mathsf{E}_{i-1}(J_{i,p})^{4}=\mathsf{E}
_{i-1}(J_{i,p})\leq\mathsf{H}_{r}(c_{k})-\mathsf{H}_{r}(c_{k-k_{J}}%
)=k_{J}H_{r}/K$ so that $\sum_{i=1}^{n}\mathsf{E}_{i-1}(J_{i,p})^{4}%
=\mathrm{O}(n^{1-5/8-\omega/2})$, uniformly in $\ell,i$. Thus,
\textitt{condition} (i) holds with $\varsigma=3/8-\omega/2\leq2\upsilon$
while \textitt{condition} (ii) holds for some $r$.

8. Note $\mathsf{E}_{i-1}|J_{i,p}(c_{k-k_{J}},c_{k})|\leq2k_{J}\delta
n^{-5/8-\omega}$ uniformly in $i,k$ by the same argument as in item 7. Since
$g_{in}=1$ then $\mathcal{\bar{R}}_{n,2}=\mathrm{O}_{\mathsf{P}%
}(n^{-5/8-\omega})=\mathrm{o}_{\mathsf{P}}(n^{-1/4})$.

B. \textitt{The term} $\mathcal{S}_{2}$. Rewrite
\[
\mathcal{S}_{2}=\sup_{0\leq\psi\leq1}\sup_{|d|\leq n^{1/4+\kappa-\eta}%
B}\bigl|
\mathbb{F}_{n}^{1,0}\bigl(0,c_{\psi}+n^{\kappa-1/2}d
\bigr)-\mathbb{F}_{n}%
^{1,0}(0,c_{\psi})\bigr|.
\]
Choosing the regressor as $x_{in}^{\ast}=n^{\kappa-1/2}$, then $\mathbb
{F}%
_{n}^{1,0}(0,c_{\psi}+n^{\kappa-1/2}d)=\mathbb{F}_{n}^{1,0}(d,c_{\psi})$.
Apply the argument of part A.
\end{pf*}

\begin{pf*}{Proof of Theorem~\ref{tFplin}}
The expression of interest is
\[
R(b,c_{\psi})=n^{1/2}\bigl\{\overline{\mathsf{F}}_{n}^{g,p}(b,c_{\psi}%
)-
\overline{\mathsf{F}}_{n}^{g,p}(0,c_{\psi})\bigr\}-
\sigma^{p-1}c_{\psi}%
^{p}
\mathsf{f}(c_{\psi})n^{-1}%
{
\sum_{i=1}^{n}} 
g_{in}n^{1/2}x_{in}^{\prime}b.
\]
Recalling the definition of $\overline{\mathsf{F}}_{n}^{g,p}$ from
(\ref{Fbar}), this satisfies $R(b,c_{\psi})=n^{-1/2}\sum_{i=1}^{n}%
g_{in}\mathcal{S}_{i}(b,c_{\psi})$, where
\[
\mathcal{S}_{i}(b,c_{\psi})=\mathsf{E}_{i-1}\bigl[
\varepsilon_{i}^{p}%
\{1_{(\varepsilon_{i}\leq\sigma c_{\psi}+b^{\prime
}x_{in})}-1_{(\varepsilon
_{i}\leq\sigma c_{\psi})}
\}\bigr]-\sigma^{p-1}x_{in}^{\prime}bc_{\psi}%
^{p}
\mathsf{f}(c_{\psi}).
\]
A bound is needed for $\mathcal{S}_{i}(b,c_{\psi})$. Let $h_{in}=\sigma
^{-1}x_{in}^{\prime}b$ and $\mathsf{g}(c)=c^{p}\mathsf{f}\mathbf{(}c)$. Write
$\mathcal{S}_{i}(b,c_{\psi})$ as an integral and Taylor expand to
second order
to get
\[
\mathcal{S}_{i}(b,c_{\psi})=\int_{c_{\psi}}^{c_{\psi}+h_{in}}
\mathsf{g}%
(c)\,\mathrm{d}c-h_{in}\mathsf{g}(c_{\psi})=
\frac{1}{2}h_{in}^{2}\mathsf{\dot{g}}%
\bigl(c^{\ast}\bigr),
\]
for an intermediate point so that $|c^{\ast}-c_{\psi}|\leq|h_{in}|$. Exploit
the bound $|b|\leq n^{1/4-\eta}B$ to get
\[
\bigl|\mathcal{S}_{i}(b,c_{\psi})\bigr|\leq\frac{1}{2}\sigma
^{-2}|b|^{2}|x_{in}|^{2}%
\sup
_{c\in\mathbb{R}}\bigl|\mathsf{\dot{g}}\bigl(c^{\ast}
\bigr)\bigr|=|x_{in}|^{2}\sup_{c\in\mathbb{R}}\bigl|\mathsf{
\dot{g}}(c)\bigr|\mathrm{O}\bigl(n^{1/2-2\eta}\bigr).
\]
Thus, by the triangular inequality
\[
\bigl|R(b,c_{\psi})\bigr|\leq n^{-1/2}%
{
\sum_{i=1}^{n}} 
|g_{in}|\bigl|S_{i}(b,c_{\psi})\bigr|
\leq\mathrm{O}\bigl(n^{-2\eta}\bigr)n^{-1}%
{
\sum_{i=1}^{n}} 
|g_{in}|\bigl|n^{1/2}x_{in}\bigr|^{2}
\sup_{c\in\mathbb{R}}\bigl|\mathsf{\dot{g}}(c)\bigr|.
\]
Due to Assumption~\ref{as}(i)(b), (iii)(b), this expression is of order
$\mathrm{O}_{\mathsf{P}}(n^{-2\eta})$ uniformly in $\psi,b$.
\end{pf*}

\begin{pf*}{Proof of Theorem~\ref{tFptight}}
1. \textitt{Coefficients} $\sigma
,\epsilon,\phi,r$. Without loss of generality, let $\sigma=1$ and
$0<\phi<1$
and $\epsilon<1$. Take $0<\epsilon$ and $n$ as well as $0<\phi^{(1-\nu
)/4}%
\leq\epsilon^{2}$ as given. Throughout, $C>0$ denotes as usual a
constant not
depending on $\phi,n,\epsilon$, which may have a different value in different
expressions. Let $r=2$. Since $\psi^{\dag}-\psi\leq\phi$, Lemma~\ref{lFtoHscale} with Assumption~\ref{as}(i)(a) shows that $0<\nu<1$ and
$C_{v},\phi_{0}>0$ exist such that $\mathsf{H}_{r}(c_{\psi^{\dag}}%
)-\mathsf{H}_{r}(c_{\psi})\leq C\phi^{1-\nu}$ for $0\leq\phi\leq\phi
_{0}$. The
proof will use a dyadic argument. Given $\epsilon,\phi,n$ we will choose
numbers $\bar{m},\underline{m}$ and derive a bound to the
probability not
depending on $\bar{m},\underline{m}$.

2. \textitt{Fine grid.} Let $\bar{m}$ satisfy $2^{-\bar{m}}\leq
n^{-1/2}\epsilon\phi^{(1-\nu)/4}\leq2^{1-\bar{m}}$.

3. \textitt{Coarse grid.} Let $\underline{m}$ satisfy $2^{-\underline{m}%
-1}H_{r}<C\phi^{1-\nu}\leq2^{-\underline{m}}H_{r}$. For large $n$,
$\bar{m}>\underline{m}$.

4. \textitt{Partition support}. For each of $m=\underline{m},\ldots
,\bar
{m}$ partition axis as laid out in (\ref{xpartition}) with $K_{m}=2^{m}$
points. For each $m$, points $c_{k_{m},m}$ and $c_{k^{\dag},m}$ exist
so that
$\underline{c}_{m}=c_{k_{m}-1,m}<c_{\psi}\leq c_{k_{m},m}=\bar{c}_{m}$
and $\underline{c}_{m}^{\dag}=c_{k_{m}^{\dag}-1,m}<c_{\psi^{\dag}}\leq
c_{k_{m}^{\dag},m}=\bar{c}_{m}^{\dag}$. Then $\bar{c}_{m-1}%
=c_{k_{m-1},m-1}$ equals either $\bar{c}_{m}=c_{k_{m},m}$ or
$c_{k_{m}+1,m}$ so that $\bar{c}_{m-1}\geq\bar{c}_{m}$ and
$\mathsf{H}(\bar{c}_{m-1})-\mathsf{H}(\bar{c}_{m})$ is either zero
or $2^{-m}H_{r}$. There is at most one $\underline{m}$-grid point in the
interval $c_{\psi},c_{\psi^{\dag}}$.

5. \textitt{Decompose} $J_{i,p}(c_{\psi},c_{\psi^{\dag}})$, see
definition in
(\ref{J}). Split the $c_{\psi},c_{\psi^{\dag}}$ interval into three intervals
where the partitioning points are $\bar{c}_{\bar{m}}$ and
$\underline{c}_{\bar{m}}^{\dag}$ which are the fine grid points to the
right of $c_{\psi}$ and to the left of $c_{\psi^{\dag}}$, respectively. Note,
that if $c_{\psi},c_{\psi^{\dag}}$ are in the same $\bar{m}$-interval
then $\bar{c}_{\bar{m}}>\underline{c}_{\bar{m}}^{\dag}$
and if
they are in neighbouring $\bar{m}$-interval then $\bar{c}%
_{\bar{m}}=\underline{c}_{\bar{m}}^{\dag}$. Thus,
\[
J_{i,p}(c_{\psi},c_{\psi^{\dag}})=J_{i,p}(c_{\psi},
\bar {c}_{\bar{m}%
})+J_{i,p}\bigl(\underline{c}_{\bar{m}}^{\dag},c_{\psi^{\dag
}}
\bigr)-1_{(\bar
{c}_{\bar{m}}>\underline{c}_{\bar{m}}^{\dag
})}J_{i,p}(\underline{c}%
_{\bar{m}},
\bar{c}_{\bar{m}})+1_{(\bar{c}_{\bar
{m}%
}<\underline{c}_{\bar{m}}^{\dag})}J_{i,p}\bigl(
\bar{c}_{\bar
{m}%
},\underline{c}_{\bar{m}}^{\dag}\bigr).
\]
Consider the fourth term. An iterative argument can be made. Since
$\bar{c}_{\bar{m}}<\underline{c}_{\bar{m}}^{\dag}$, the coarser
$(\bar{m}-1)$-grid satisfies $\bar{c}_{\bar{m}}\leq
\bar
{c}_{\bar{m}-1}\leq\underline{c}_{\bar{m}-1}^{\dag}\leq
\underline{c}_{\bar{m}}^{\dag}$, so that
\[
J_{i,p}\bigl(\bar{c}_{\bar{m}},\underline{c}_{\bar{m}}^{\dag
}
\bigr)=J_{i,p}(\bar{c}_{\bar{m}},\bar{c}_{\bar{m}-1}%
)+J_{i,p}
\bigl(\bar{c}_{\bar{m}-1},\underline{c}_{\bar
{m}-1}^{\dag
}
\bigr)+J_{i,p}\bigl(\underline{c}_{\bar{m}-1}^{\dag},
\underline {c}_{\bar{m}%
}^{\dag}\bigr).
\]
If $\bar{c}_{\bar{m}-1}=\underline{c}_{\bar{m}-1}^{\dag
}$, then
$J_{i,p}(\bar{c}_{\bar{m}-1},\underline{c}_{\bar
{m}-1}^{\dag
})=0$ and the iteration stops, noting that for $m<\bar{m}-1$ the $m$-grid
points cross over so that $\bar{c}_{m}\geq\bar{c}_{\bar
{m}%
-1}=\underline{c}_{\bar{m}-1}^{\dag}\geq\underline{c}_{m}^{\dag}$. If
$\bar{c}_{\bar{m}-1}<\underline{c}_{\bar{m}-1}^{\dag}$, the
argument can be made again for $J_{i,p}(\bar{c}_{\bar{m}%
-1},\underline{c}_{\bar{m}-1}^{\dag})$. In the $m$th step, the iteration
continues if $\bar{c}_{m}<\underline{c}_{m}^{\dag}$, so that if
there are
no other $m$-grid points between $\bar{c}_{\bar{m}}$ and
$\underline{c}_{\bar{m}}^{\dag}$, the contribution from the $(m-1)$-step
is zero. Because there is at most one $\underline{m}$-point in the interval
$c_{\psi},c_{\psi^{\dag}}$, the $\underline{m}$-step will either give a zero
contribution or the grid points will have crossed over at an earlier stage.
Therefore, the fourth term satisfies
\[
1_{(\bar{c}_{\bar{m}}<\underline{c}_{\bar{m}}^{\dag})}%
J_{i,p}\bigl(\bar{c}_{\bar{m}},
\underline{c}_{\bar{m}}^{\dag
}\bigr)=%
{
\sum_{m=\underline{m}+1}^{\bar{m}}} 
1_{(\bar{c}_{m}<\underline{c}_{m}^{\dag})}
\bigl\{J_{i,p}(\bar{c}%
_{m},
\bar{c}_{m-1})+J_{i,p}\bigl(\underline{c}_{m-1}^{\dag},
\underline {c}%
_{m}^{\dag}\bigr)\bigr\}.
\]

6. \textitt{Decompose} $\mathcal{S}=n^{1/2}\{\mathbb{F}_{n}^{g,p}%
(0,c_{\psi^{\dag}})-\mathbb{F}_{n}^{g,p}(0,c_{\psi})\}$. Due to the
decomposition of $J_{i,p}(c_{\psi}, c_{\psi^{\dag}})$ in item 5, then
$|\mathcal{S}|\leq|Z_{1}|+|Z_{2}|+|Z_{3}|+|Z_{4}|+|Z_{5}|$, where
\begin{eqnarray*}
Z_{1} & =&\frac{1}{\sqrt{n}}%
{ \sum
_{i=1}^{n}} 
g_{in}
\bigl[J_{i,p}(c_{\psi},\bar{c}_{\bar{m}})-
\mathsf{E}_{i-1}%
\bigl\{J_{i,p}(c_{\psi},
\bar{c}_{\bar{m}})\bigr\}\bigr],
\\
Z_{2} & =&\frac{1}{\sqrt{n}}%
{ \sum
_{i=1}^{n}} 
g_{in}
\bigl[J_{i,p}\bigl(\underline{c}_{\bar{m}}^{\dag},c_{\psi^{\dag}}%
\bigr)-\mathsf{E}_{i-1}\bigl\{J_{i,p}\bigl(
\underline{c}_{\bar{m}}^{\dag},c_{\psi
^{\dag
}}\bigr)\bigr\}\bigr],
\\
Z_{3} & =&1_{(\bar{c}_{\bar{m}}>\underline{c}_{\bar
{m}}^{\dag
})}\frac{1}{\sqrt{n}}%
{
\sum_{i=1}^{n}} 
g_{in}
\bigl[J_{i,p}(\underline{c}_{\bar{m}},\bar{c}_{\bar{m}%
})-
\mathsf{E}_{i-1}\bigl\{J_{i,p}(\underline{c}_{\bar{m}},
\bar {c}_{\bar{m}})\bigr\}\bigr],
\\
Z_{4} & =&%
{ \sum
_{m=\underline{m}+1}^{\bar{m}}} 
1_{(\bar{c}_{m}<\underline{c}_{m}^{\dag})}
\frac{1}{\sqrt{n}}%
{ \sum
_{i=1}^{n}} 
g_{in}
\bigl[J_{i,p}(\bar{c}_{m},\bar{c}_{m-1})-
\mathsf{E}_{i-1}%
\bigl\{J_{i,p}(\bar{c}_{m},
\bar{c}_{m-1})\bigr\}\bigr],
\\
Z_{5} & =&%
{ \sum
_{m=\underline{m}+1}^{\bar{m}}} 
1_{(\bar{c}_{m}<\underline{c}_{m}^{\dag})}
\frac{1}{\sqrt{n}}%
{ \sum
_{i=1}^{n}} 
g_{in}
\bigl[J_{i,p}\bigl(\underline{c}_{m-1}^{\dag},
\underline{c}_{m}^{\dag}%
\bigr)-
\mathsf{E}_{i-1}\bigl\{J_{i,p}\bigl(\underline{c}_{m-1}^{\dag},
\underline{c}_{m} 
^{\dag}\bigr)\bigr\}\bigr].
\end{eqnarray*}

7. \textitt{The term} $Z_{1}$: \textitt{Finding martingale}. Bound
$|J_{i,p}(c_{\psi},\bar{c}_{\bar{m}})|\leq|J_{i,p}(\underline
{c}%
_{\bar{m}},\bar{c}_{\bar{m}})|$ where the points
$\underline{c}%
_{\bar{m}},\bar{c}_{\bar{m}}$ are two neighbouring
points on
the $\bar{m}$-grid, but their location depends on $\psi$. It
follows that
\begin{eqnarray*}
\sup_{0\leq\psi\leq\psi^{\dag}\leq1\dvt \psi^{\dag}-\psi\leq\phi}|Z_{1}|
 &\leq&\frac{1}{\sqrt{n}}%
{ \sum
_{i=1}^{n}} 
|g_{in}|\bigl
\{\bigl|J_{i,p}(\underline{c}_{\bar{m}},\bar{c}_{\bar
{m}%
})\bigr|+
\mathsf{E}_{i-1}\bigl|J_{i,p}(\underline{c}_{\bar{m}},
\bar {c}_{\bar{m}})\bigr|\bigr\}
\\
& \leq&\max_{1\leq\ell\leq2^{\bar{m}}}\frac{1}{\sqrt{n}}%
{
\sum_{i=1}^{n}} 
|g_{in}|
\bigl\{\bigl|J_{i,p}(c_{\ell-1,m},c_{\ell,m})\bigr|+
\mathsf{E}_{i-1}\bigl|J_{i,p}%
(c_{\ell-1,m},c_{\ell,m})\bigr|
\bigr\}.
\end{eqnarray*}
Thus, a martingale decomposition gives
\[
\sup_{0\leq\psi\leq\psi^{\dag}\leq1\dvt \psi^{\dag}-\psi\leq\phi}|Z_{1}|\leq \max_{1\leq\ell\leq2^{\bar{m}}}|
\tilde{V}_{1,\ell,\bar{m}%
}|+2\max_{1\leq\ell\leq2^{\bar{m}}}\overline{V}_{1,\ell,\bar{m}},
\]
where%
%
\begin{eqnarray}
\tilde{V}_{1,\ell,\bar{m}} & =&\frac{1}{\sqrt{n}}%
{
\sum_{i=1}^{n}} 
|g_{in}|
\bigl[\bigl|J_{i,p}(c_{\ell-1,\bar{m}},c_{\ell,\bar
{m}})\bigr|-
\mathsf{E}%
_{i-1}\bigl\{\bigl|J_{i,p}(c_{\ell-1,\bar{m}},c_{\ell,\bar{m}}%
)\bigr|
\bigr\}\bigr],\label{ptFptightZ1}
\\
\overline{V}_{1,\ell,\bar{m}} & =&\frac{1}{\sqrt{n}}%
{
\sum_{i=1}^{n}} 
|g_{in}|
\mathsf{E}_{i-1}\bigl\{\bigl|J_{i,p}(c_{\ell-1,\bar{m}},c_{\ell
,\bar{m}})\bigr|
\bigr\}.
\end{eqnarray}
8. \textitt{The term} $Z_{1}$: \textitt{The compensator} $\overline{V}$. Since
\[
\mathsf{E}_{i-1}\bigl\{\bigl|J_{i,p}(c_{\ell-1,\bar{m}},c_{\ell,\bar{m}%
})\bigr|^{2^{r}}
\bigr\}\leq\mathsf{H}_{r}(c_{\ell,\bar{m}})-\mathsf{H}_{r}%
(c_{\ell-1,\bar{m}})=2^{-\bar{m}}H_{r},
\]
Assumption~\ref{as}(i)(a), (iii)(a) implies%
%
\begin{equation}
\mathsf{E}\max_{1\leq\ell\leq2^{\bar{m}}}%
{ \sum
_{i=1}^{n}} 
|g_{in}|^{2^{r}}
\mathsf{E}_{i-1}\bigl\{\bigl|J_{i,p}(c_{\ell-1,\bar{m}}%
,c_{\ell,\bar{m}})\bigr|^{2^{r}}
\bigr\}\leq nC2^{-\bar{m}}H_{r}%
.\label{ptFptightMEbound}%
\end{equation}
Item 2 shows $2^{-\bar{m}}\leq n^{-1/2}\epsilon\phi^{(1-\nu)/4}$. Thus,
the Markov inequality implies
\[
\mathsf{P}\Bigl(\max_{1\leq\ell\leq2^{\bar{m}}}\overline{V}_{1,\ell
,\bar{m}}>
\epsilon\Bigr)\leq\frac{1}{\epsilon}\mathsf{E}\max_{1\leq\ell
\leq2^{\bar{m}}}
\overline{V}_{1,\ell,\bar{m}}\leq n^{1/2}\frac
{1}{\epsilon}C2^{-\bar{m}}H_{r}=C
\phi^{(1-\nu)/4}.
\]
9. \textitt{The term} $Z_{1}$: \textitt{The martingale} $\tilde{V}$. Apply
Theorem~\ref{lsupM4finite} with $z_{\ell,i}=g_{in}^{\ast}|J_{i,p}%
(c_{\ell,\bar{m}},c_{\ell+1,\bar{m}})|$ where $g_{in}^{\ast}$
is a
coordinate of $|g_{in}|$, and with $L=2^{\bar{m}}$ and $\kappa
=\epsilon$
while $D=C2^{-\bar{m}}$ by the inequality (\ref{ptFptightMEbound}),
to get
%
\begin{equation}
\mathsf{P}\Bigl(\max_{1\leq\ell\leq2^{\bar{m}}}|\tilde{V}_{1,\ell
,\bar{m}}|>
\epsilon\Bigr)\leq C2^{-\bar{m}}\frac{\theta}{\epsilon}%
+C
\frac{\theta^{3}}{n\epsilon}+C2^{\bar{m}}\exp\biggl(-\frac{\epsilon
\theta
}{14}
\biggr),\label{ptFptightV1bound}%
\end{equation}
where we can choose $\theta=14\epsilon^{-1}(\log2^{2\bar{m}}+\log
\phi^{-1})$. First term in (\ref{ptFptightV1bound}) satisfies
\[
C2^{-\bar{m}}\frac{\theta}{\epsilon}\leq C\frac{1}{\epsilon^{2}%
}2^{-\bar{m}/2}
\bigl\{\bar{m}2^{-\bar{m}/2}+\phi^{-(1-\nu
)/2}2^{-\bar{m}/2}
\phi^{(1-\nu)/2}\log\phi^{-1}\bigr\}\leq C\phi^{(1-\nu)/4},
\]
since the bounds in items 1, 3 imply $\epsilon^{-2}\leq\phi^{-(1-\nu
)/4}$ and
$2^{-\bar{m}/2}\leq2^{-\underline{m}/2}<C\phi^{(1-\nu)/2}$, while the
functions $m2^{-m/2}$ and $\phi^{(1-\nu)/2}\log\phi^{-1}$ are bounded for
$m\geq1$ and $0<\phi<1$. Second term in~(\ref{ptFptightV1bound}): Use
first the definition of $\theta$ with the inequality $(x+y)^{3}\leq
C(x^{3}+y^{3})$ and then that the bounds in items 1, 2 imply $\epsilon
^{-2}\leq\phi^{-(1-\nu)/4}$ and $n^{-1}\epsilon^{2}\leq C\phi^{-(1-\nu
)/2}2^{-2\bar{m}}$ so that
\[
\mathcal{P}_{1}=C\frac{\theta^{3}}{n\epsilon}\leq C\frac{1}{n\epsilon
^{4}%
}\biggl(
\bar{m}^{3}+\log^{3}\frac{1}{\phi}\biggr)\leq C
\phi^{-(1-\nu)/2}%
2^{-2\bar{m}}\phi^{-3(1-\nu)/4}\biggl(
\bar{m}^{3}+\log^{3}\frac
{1}{\phi}\biggr).
\]
Rewrite this bound and argue as for the first term, to get that
\begin{eqnarray*}
\mathcal{P}_{1}&\leq& C\bigl\{2^{-\bar{m}/2}\phi^{-(1-\nu)/2}\bigr\}^{3}\biggl\{
\bar{m}%
^{3}2^{-\bar{m}/2}+\phi^{-(1-\nu)/2}2^{-\bar{m}/2}
\phi^{(1-\nu
)/2}\log^{3}\frac{1}{\phi}\biggr\}
\phi^{(1-\nu)/4}\\
&\leq& C\phi^{(1-\nu)/4}.
\end{eqnarray*}
Third term in (\ref{ptFptightV1bound}) satisfies
\[
C2^{\bar{m}}\exp\biggl(-\frac{\epsilon\theta}{14}\biggr)=C2^{-\bar{m}}\phi
\leq C\phi^{(1-\nu)/4},
\]
since $\phi\leq\phi^{(1-\nu)/4}$ for $0<\phi<1$. In summary, $\mathsf{P}
(\max_{1\leq\ell\leq2^{\bar{m}}}|\tilde{V}_{1,\ell,\bar{m}
}|>\epsilon)\leq C\phi^{(1-\nu)/4}$.

10. \textitt{The terms} $Z_{2}$ \textitt{and} $Z_{3}$. Apply the same argument
as in items 7--9.

11. \textitt{The term} $Z_{4}$: \textitt{finding martingale.} Recall
that, for instance, $\bar{c}_{m}=c_{k_{m},m}$ while $\bar{c}_{m-1}$
either equals $c_{k_{m},m}$ or $c_{k_{m}+1,m}$, so that $\bar{c}%
_{m},\bar{c}_{m-1}$ are at most 1 step apart in the $m$-grid. Let
\[
M_{\ell,m,n}=\frac{1}{\sqrt{n}}%
{ \sum
_{i=1}^{n}} 
g_{in}
\bigl[J_{i,p}(c_{\ell,m},c_{\ell+1,m})-\mathsf{E}_{i-1}
\bigl\{ J_{i,p}(c_{\ell
,m},c_{\ell+1,m})\bigr\}\bigr].
\]
It then holds that
\[
|Z_{4}|\leq%
{ \sum
_{m=\underline{m}+1}^{\bar{m}}} 
|M_{k_{m},m,n}|
\leq%
{ \sum_{m=\underline{m}+1}^{\bar{m}}}
\max_{1\leq\ell\leq2^{m}}|M_{\ell,m,n}|.
\]
Note that $%
{ \sum_{m=\underline{m}+1}^{\bar{m}}}
2^{(\underline{m}-m)/4}\leq\sum_{j=1}^{\infty}2^{-j/4}=(2^{1/4}-1)^{-1}<6$,
and that the right-hand side does not depend on $\psi$. It therefore
holds
\[
\mathcal{P}_{4}=\mathsf{P}\Bigl(\sup_{0\leq\psi\leq\psi^{\dag}\leq1\dvt \psi
^{\dag
}-\psi\leq\phi}|Z_{4}|>
\epsilon\Bigr)\leq\mathsf{P}%
{ \bigcup
_{m=\underline{m}+1}^{\bar{m}}} 
 \biggl\{ \max
_{1\leq\ell\leq2^{m}}|M_{\ell,m,n}|>\frac{2^{(\underline{m}%
-m)/4}\epsilon}{6} \biggr\}.
\]
Using Boole's inequality, then
\[
\mathcal{P}_{4}\leq%
{ \sum
_{m=\underline{m}+1}^{\bar{m}}} 
\mathsf{P} \biggl\{ \max
_{1\leq\ell\leq2^{m}}|M_{\ell,m,n}|>\frac
{2^{(\underline{m}-m)/4}\epsilon}{6} \biggr\}.
\]

12. \textitt{The term} $Z_{4}$: \textitt{apply Lemma~\ref{lsupM4finite}} with
$z_{\ell,i}=g_{in}^{\ast}J_{i,p}(c_{\ell-1,m},c_{\ell,m})$ where $g_{in}
^{\ast}$ is a coordinate of $g_{in}$ and with $L=2^{m}$ while $\kappa
=2^{(\underline{m}-m)/4}\epsilon/6$ and $D=C2^{-m}$, due to the inequality
(\ref{ptFptightMEbound}) with $\bar{m}$ replaced by $m$. Thus%
%
\begin{equation}
\mathcal{P}_{4}\leq C%
{ \sum
_{m=\underline{m}+1}^{\bar{m}}} 
 \biggl\{ 2^{-m}
\frac{\theta_{m}}{2^{(\underline{m}-m)/4}\epsilon}+\frac
{\theta_{m}^{3}}{n2^{(\underline{m}-m)/4}\epsilon}+2^{m}\exp\biggl(-
\frac
{2^{(\underline{m}-m)/4}\epsilon\theta_{m}}{84}\biggr) \biggr\},\label{ptFptightV4bound}%
\end{equation}
where we choose $2^{(\underline{m}-m)/4}\epsilon\theta_{m}/84=\log
(4^{m-\underline{m}})+\log\phi^{-1}$. First term in
(\ref{ptFptightV4bound}) satisfies
\[
\mathcal{P}_{41}=%
{ \sum
_{m=\underline{m}+1}^{\bar{m}}} 
2^{-m}
\frac{\theta_{m}}{2^{(\underline{m}-m)/4}\epsilon}\leq C%
{ \sum
_{m=\underline{m}+1}^{\bar{m}}} 
\frac{1}{2^{(\underline{m}+m)/2}\epsilon^{2}}\biggl
\{(m-\underline{m})+\log \frac
{1}{\phi}\biggr\}.
\]
Note $2^{-(\underline{m}+m)/2}=2^{-(m-\underline{m})/2}2^{-\underline{m}}$.
Items 1, 3 imply $\epsilon^{-2}\leq\phi^{-(1-\nu)/4}$ and $2^{-\bar
{m}%
/2}\leq2^{-\underline{m}/2}<C\phi^{(1-\nu)/2}$. Next, use that
geometric sums
are finite and argue as in item 9 to see that
\[
\mathcal{P}_{41}\leq C%
{ \sum
_{m=\underline{m}+1}^{\bar{m}}} 
2^{-(m-\underline{m})/2}
\biggl\{2^{-\underline{m}}(m-\underline{m})+\phi ^{(1-\nu
)}\log
\frac{1}{\phi}\biggr\}\phi^{(1-\nu)/4}\leq C\phi^{(1-\nu)/4}.
\]
The second term in (\ref{ptFptightV4bound}) satisfies
\[
\mathcal{P}_{42}=%
{ \sum
_{m=\underline{m}+1}^{\bar{m}}} 
\frac{\theta_{m}^{3}}{n2^{(\underline{m}-m)/4}\epsilon}\leq
C%
{ \sum_{m=\underline{m}+1}^{\bar{m}}}
\frac{1}{n2^{(\underline{m}-m)}\epsilon^{4}}\biggl\{(m-\underline{m})^{3}+\log
^{3}\frac{1}{\phi}\biggr\}.
\]
Items 1, 2 imply $\epsilon^{-2}\leq\phi^{-(1-\nu)/4}$ and $n^{-1}\epsilon
^{2}\leq\phi^{-(1-\nu)/2}2^{2-2\bar{m}}$ so that
\[
\mathcal{P}_{42}\leq C%
{ \sum
_{m=\underline{m}+1}^{\bar{m}}} 
\phi^{-5(1-\nu)/4}2^{-2\bar{m}-\underline{m}+m}\biggl\{(m-\underline{m}%
)^{3}+
\log^{3}\frac{1}{\phi}\biggr\}.
\]
Rewrite $2^{-2\bar{m}-\underline{m}+m}=2^{3(m-\bar{m}%
)/2-(m-\underline{m})/2-\bar{m}/2-3\underline{m}/2}$ to get that
$\mathcal{P}_{42}$ is bounded by
\[
C%
{ \sum_{m=\underline{m}+1}^{\bar{m}}}
\bigl\{\phi^{-(1-\nu)/2}2^{-\underline{m}/2}\bigr\}^{3}2^{3(m-\bar{m}%
)/2}2^{-(m-\underline{m})/2}
\biggl\{2^{-\bar{m}/2}(m-\underline{m}%
)^{3}+2^{-\bar{m}/2}
\log^{3}\frac{1}{\phi}\biggr\}\phi^{(1-\nu)/4}.
\]
Argue as for first term using $2^{-\bar{m}/2}\leq2^{-\underline{m}%
/2}<C\phi^{(1-\nu)/2}$ from item 3 to get $\mathcal{P}_{42}\leq C\phi
^{(1-\nu)/4}$.

The third term in (\ref{ptFptightV4bound}) satisfies, noting
$2^{\underline{m}}\leq C\phi^{\nu-1}$
\[
\mathcal{P}_{43}=%
{ \sum
_{m=\underline{m}+1}^{\bar{m}}} 
2^{m}\exp
\biggl(-\frac{2^{(\underline{m}-m)/4}\epsilon\theta_{m}}{84}\biggr)=%
{ \sum
_{m=\underline{m}+1}^{\bar{m}}} 
2^{-(m-\underline{m})}2^{\underline{m}}
\phi.
\]
Noting that $2^{\underline{m}}\leq C\phi^{\nu-1}$ then $\mathcal
{P}_{43}\leq C%
{ \sum_{m=\underline{m}+1}^{\bar{m}}}
2^{-(m-\underline{m})}\phi^{\nu}=C\phi^{\nu}$.

13. \textitt{The terms} $Z_{5}$. Apply the same argument as for $Z_{4}$.

14. \textitt{Combine} the bounds from items 8, 9, 10, 12, 13 to get
\[
\mathsf{P}\Bigl(\sup_{0\leq\psi\leq\psi^{\dag}\leq1\dvt \psi^{\dag}-\psi\leq\phi
}|\mathcal{S}|>\epsilon\Bigr)
\leq%
{ \sum_{j=1}^{5}}
\mathsf{P}\Bigl(\sup_{0\leq\psi\leq\psi^{\dag}\leq1\dvt \psi^{\dag}-\psi\leq\phi} 
|Z_{j}|>
\epsilon\Bigr)\leq C\phi^{(1-\nu)/4}+2C\phi^{\nu}.
\]
For a given $\epsilon>0$ the only constraint to $\phi$ is that $0<\phi
^{(1-\nu)/4}\leq\epsilon^{2}$. Thus, the probability vanishes as
$\phi\downarrow0$.
\end{pf*}

\section{Proofs of main Theorems \texorpdfstring{\protect\ref{txi}}{3.1}--\texorpdfstring{\protect\ref{tK}}{3.7}}
\label{sGproof}

The main results for the Forward Search are proved in a series of steps.
Theorem~\ref{txi} shows that asymptotically the forward residuals
behave like
the quantiles of the absolute errors $|\varepsilon_{i}|$. It is therefore
useful to start by reviewing some known results from the theory of quantile
processes. Second, the Forward Search problem is reformulated in terms
of a
weighted and marked absolute empirical distribution function
$\widehat{\mathsf{G}}_{n}$. At this point, we work with absolute errors
and it
is natural to move from the general densities of Assumption~\ref{as} to the
symmetric densities of Assumption~\ref{asFS}. Third, this empirical
distribution function is analysed using the results from Section~\ref{sempproc}. Fourth, the corresponding quantile processes are analysed.
Fifth, a single step of the Forward Search is analysed using these results.
Sixth, the iteration of the Forward Search is analysed.

\subsection{Some known results from the theory of quantile
processes}\label{ssknown}

Introduce the empirical distribution function of the absolute errors,
$|\varepsilon_{i}|/\sigma$, that is%
%
\begin{equation}
\widehat{\mathsf{G}}_{n}(c)=\frac{1}{n}%
{
\sum_{i=1}^{n}} 
1_{(|\varepsilon_{i}|\leq\sigma c)}.
\label{Ghat0}%
\end{equation}
The first result gives the asymptotic distribution of the empirical
process
\[
\mathbb{G}_{n}(c_{\psi})=n^{1/2}\bigl\{\widehat{
\mathsf{G}}_{n}(c_{\psi
})-\psi\bigr\}.
\]

\begin{lemma}[(Billingsley \cite{Billingsley}, Theorem~14.3)]
\label{lGhat0} Let $\mathbb{B}$ be a
Brownian bridge so that $\mathbb{B}(\psi)$ is $\mathsf{N}\{0,\psi(1-\psi
)\}$-distributed. Then, it holds $\mathbb{G}_{n}\overset{\mathsf{D}%
}{\rightarrow}\mathbb{B}$ on $D[0,1]$.
\end{lemma}

The empirical quantiles of the absolute errors, $|\varepsilon
_{i}|/\sigma$,
are defined as%
%
\begin{equation}
\hat{c}_{\psi}=\widehat{\mathsf{G}}_{n}^{-1}(\psi)=
\inf\bigl\{c\dvt \widehat {\mathsf{G}%
}_{n}(c)\geq\psi\bigr\}.
\label{Ginv0}%
\end{equation}
Empirical quantiles and empirical distribution functions are linked as follows.

\begin{lemma}[(Cs\"{o}rg\H{o} \cite{Cs}, Corollaries 6.2.1, 6.2.2)]
\label{tdhat0} Suppose that
$\mathsf{f}$ is symmetric, differentiable, positive for $\mathsf{F}%
^{-1}(0)<c<\mathsf{F}^{-1}(1)$, decreasing for large $c$, and satisfying
$\gamma=\sup_{c>0}\mathsf{F}(c)\{1-\mathsf{F}(c)\}|\mathsf{\dot{f}%
}(c)|/\{\mathsf{f}(c)\}^{2}<\infty$.\vadjust{\goodbreak}

Then, for all $\zeta>0$,
\begin{longlist}[(a)]
\item[(a)]
$\sup_{0\leq\psi\leq1}|2\mathsf{f}(c_{\psi})n^{1/2}(\hat{c}_{\psi
}-c_{\psi
})+n^{1/2}\{\widehat{\mathsf{G}}_{n}(c_{\psi})-\psi\}|=\mathrm
{o}_{\mathsf{P}%
}(n^{\zeta-1/4})$;

\item[(b)] $\sup_{0\leq\psi\leq1}|2\mathsf
{f}(c_{\psi
})n^{1/2}(\hat{c}_{\psi}-c_{\psi})-n^{1/2}\{\mathsf{G}(\hat{c}_{\psi}%
)-\psi\}|=\mathrm{o}_{\mathsf{P}}(n^{\zeta-1/2})$;

\item[(c)] $\sup_{0\leq\psi\leq1}|n^{1/2}\{\mathsf{G(}\hat{c}_{\psi})-\psi\}+n^{1/2}%
\{\widehat{\mathsf{G}}_{n}(c_{\psi})-\psi\}|=\mathrm{o}_{\mathsf{P}}%
(n^{\zeta-1/4})$.
\end{longlist}
\end{lemma}

The result in Lemma~\ref{tdhat0}(a) shows that the empirical quantile
$\hat{c}_{\psi}$ satisfies, for $0<\psi<1$,
\[
n^{1/2}(\hat{c}_{\psi}-c_{\psi})=\frac{1}{2\mathsf{f}(c_{\psi})}n^{1/2}%
\bigl\{\psi-\widehat{\mathsf{G}}_{n}(c_{\psi})\bigr\}+
\mathrm{o}_{\mathsf{P}}(1).
\]
This is known as the Bahadur \cite{Bahadur} representation. Kiefer \cite{Kiefer}, equations (1.8), (1.9),
studied parts~(b), (c), which combine to (a). More details can be
found in
Cs\"{o}rg\H{o} \cite{Cs} who also gives almost sure, logarithmic rates.

Some weighted versions of the above results are also needed.

\begin{lemma}[(Shorack \cite{Shorack}, Cs\"{o}rg\H{o} \cite{Cs}, Theorem~5.1.1)]
\label{tOReilly} Let the
function $q_{\psi}$ be symmetric about $\psi=1/2$ (it suffices if
$q_{\psi}$
is bounded below by such a function), such that $q_{\psi}$ is
increasing and
continuous on $0\leq\psi\leq1/2$ and satisfies $q_{\psi}=\{\psi\log\log
(1/\psi)\}^{1/2}g_{\psi}$ for a function $g_{\psi}$ so $\lim_{\psi
\rightarrow0}g_{\psi}=\infty$. Then, a probability space exists on
which one
can define a Brownian bridge $\mathbb{B}_{n}$ for each $n$, so
that:

\begin{longlist}[(a)]
\item[(a)] $\sup_{0\leq\psi\leq1}|\{\mathbb{G}_{n}(c_{\psi})-\mathbb{B}_{n}%
(\psi)\}/q_{\psi}|=\mathrm{o}_{\mathsf{P}}(1)$;

\item[(b)] $\sup_{1/(n+1)\leq\psi\leq n/(n+1)}|\{\mathsf{f}(c_{\psi})n^{1/2}(\hat
{c}_{\psi
}-c_{\psi})-\mathbb{B}_{n}(\psi)\}/q_{\psi}|=\mathrm{o}_{\mathsf{P}}(1)$
provided the assumptions of Lemma~\textup{\ref{tdhat0}} hold.
\end{longlist}
\end{lemma}

In Lemma~\ref{tOReilly} a possible choice of $q_{w}$ is $\{\psi
(1-\psi)\}^{\alpha}$ for $\alpha<1/2$, which will be used in the proof of
Theorem~\ref{tshat}. Finally, a continuity property of the Brownian
bridge is needed.

\begin{lemma}[(Revuz and Yor \cite{RevuzYor}, Theorem~1.2.2)]
\label{tBB} A Brownian motion
$\mathbb{W}$ is locally H\"{o}lder continuous of order $\alpha$ for all
$\alpha<1/2$. That is,
\[
\sup_{0\leq\psi<\psi^{\dag}\leq1}\frac{|\mathbb{W}(\psi^{\dag})-\mathbb
{W}%
(\psi)|}{(\psi^{\dag}-\psi)^{\alpha}}\overset{a.s.} {<}\infty.
\]
Thus, for a Brownian bridge $\mathbb{B}$, $\lim_{\psi\rightarrow0}%
\mathbb{B}(\psi)/\psi^{\alpha}=0$ a.s.
\end{lemma}

\subsection{Absolute empirical process representation}\label{ssG}

Normalisations are needed for estimators and regressors. Depending on the
stochastic properties of the regressor $x_{i}$, choose a non-stochastic
normalisation matrix $N$ and define
\[
\hat{b}=N^{-1}(\hat{\beta}-\beta),\qquad x_{in}=N^{\prime}x_{i},
\]
so that $\sum_{i=1}^{n}x_{in}x_{in}^{\prime}$ converges, $n^{-1/2}\sum_{i=1}^{n}|x_{in}|$ is bounded, and $x_{i}^{\prime}(\hat{\beta}-\beta
)=x_{in}^{\prime}b$. If, for example, $(y_{i},x_{i})$ is stationary then
$N=n^{-1/2}I_{\dim x}$ so that $b=n^{1/2}(\hat{\beta}-\beta)$ and
$x_{in}=n^{-1/2}x_{i}$. If $x_{i}$ is a random walk then $N=n^{-1}$.

Introduce matrix-valued weights $g_{in}$ of the form $1$,
$n^{1/2}Nx_{i}$ or
$nNx_{i}x_{i}^{\prime}N$, so that the sum $n^{-1}\sum_{i=1}^{n}|g_{in}|$ is
bounded. In the stationary case, $g_{in}$ will be $1$, $x_{i}$ or
$x_{i}%
x_{i}^{\prime}$. When $x_{i}$ is a random walk, $g_{in}$ is $1$,
$n^{-1/2}x_{i}$ or $n^{-1}x_{i}x_{i}^{\prime}$.

Define the \textitt{weighted and marked absolute empirical distribution
functions}%
%
\begin{equation}
\widehat{\mathsf{G}}_{n}^{g,p}(b,c)=\frac{1}{n}%
{
\sum_{i=1}^{n}} 
g_{in}
\varepsilon_{i}^{p}1_{(|\varepsilon_{i}-x_{in}^{\prime}b|\leq
\sigma c)}, \label{Ghat}%
\end{equation}
for $b\in\mathbb{R}^{\dim x}$ and $c\geq0$. Here the weights are
$g_{in}$ and
the marks $\varepsilon_{i}^{p}$. Four combinations of weights and marks
are of
interest in the analysis of the Forward Search. The deletion residuals involve
$g_{in}=1$, $p=0$. The least squares estimator involves $g_{in}=n^{1/2}%
N^{\prime}x_{i}$, $p=1$ and $g_{in}=nN^{\prime}x_{i}x_{i}^{\prime}N$, $p=0$.
The variance estimator involves the terms mentioned as well as $g_{in}=1$,
$p=2$. When $p=0$, the marks are $\varepsilon_{i}^{0}=1$ so that
$\widehat{\mathsf{G}}_{n}^{g,0}$ is a weighted absolute empirical distribution
function, similar to that studied by Koul and Ossiander \cite{KoulOssiander}. When also
$b=0$, then $\widehat{\mathsf{G}}_{n}^{1,0}$ equals the empirical distribution
function $\widehat{\mathsf{G}}_{n}$ of (\ref{Ghat0}).

The Forward Search Algorithm \ref{alestimators} can now be cast as follows.
Step $(m+1)$ results in an order statistic
%
\begin{equation}
\hat{z}^{(m)}=\sigma\inf\biggl\{c\dvt \widehat{\mathsf{G}}_{n}^{1,0}
\bigl(\hat{b}%
^{(m)},c\bigr)\geq\frac{m+1}{n}\biggr\},
\label{xihat}%
\end{equation}
where $g_{in}=1$, $p=0$, so that%
%
\begin{equation}
\frac{m+1}{n}=\widehat{\mathsf{G}}_{n}^{1,0}\biggl(
\hat{b}^{(m)},\frac{\hat{z}
^{(m)}}{\sigma}\biggr)=\frac{1}{n}%
{
\sum_{i=1}^{n}} 
1_{(|\varepsilon_{i}-x_{in}^{\prime}\hat{b}^{(m)}|\leq\hat
{z}^{(m)})}=
\frac
{1}{n}%
{ \sum
_{i\in S^{(m+1)}}} 
1. \label{xihat2}%
\end{equation}
The least squares estimator has estimation error%
%
\begin{eqnarray}\label{betahatG}%
\hat{b}^{(m+1)}&=&N^{-1}\bigl(\hat{\beta}^{(m)}-\beta
\bigr)
\nonumber
\\[-8pt]
\\[-8pt]
\nonumber
&=&\biggl\{\widehat{\mathsf{G}}%
_{n}^{xx,0}
\biggl(\hat{b}^{(m)},\frac{\hat{z}^{(m)}}{\sigma}\biggr)\biggr\}^{-1}
\biggl\{n^{1/2}%
\widehat{\mathsf{G}}_{n}^{x,1}
\biggl(\hat{b}^{(m)},\frac{\hat{z}^{(m)}}{\sigma
}\biggr)\biggr\},
\end{eqnarray}
while the asymptotically bias corrected least squares variance estimator
satisfies%
%
\begin{eqnarray}\label
{sigmahatG}%
&&n^{1/2}\bigl\{\bigl(\hat{\sigma}_{\mathrm{cor}}^{(m+1)}
\bigr)^{2}-\sigma^{2}\bigr\}
\nonumber
\\[-8pt]
\\[-8pt]
\nonumber
&&\qquad=\frac{n^{1/2}}%
{\tau_{m/n}}\biggl[
\widehat{\mathsf{G}}_{n}^{1,2}\biggl(\hat{b}^{(m)},
\frac{\hat
{z}^{(m)}%
}{\sigma}\biggr)-\bigl\{\hat{b}^{(m+1)}\bigr\}^{\prime}
\widehat{\mathsf {G}}_{n}^{xx,0}\biggl(\hat
{b}^{(m)},\frac{\hat{z}^{(m)}}{\sigma}\biggr)\bigl\{\hat{b}^{(m+1)}\bigr\}
\biggr].
\end{eqnarray}

\subsection{The absolute empirical distribution}

The process $\widehat{\mathsf{G}}_{n}^{g,p}$ is now analysed using the
auxiliary Theorems \ref{tFpest}--\ref{tFptight} for the process
$\widehat{\mathsf{F}}_{n}^{g,p}$. Only the four combinations of $g_{in}%
,p$ are now considered as outlined in Section~\ref{ssG}. When checking
Assumption~\ref{as} it suffices to check the conditions for the hybrid case
where $g_{in}=nN^{\prime}x_{i}x_{i}^{\prime}N$ and
$p=2$. The process $\widehat{\mathsf{G}}_{n}^{g,p}$ can be
expressed in
terms of $\widehat{\mathsf{F}}_{n}^{g,p}$ by%
%
\begin{equation}
\widehat{\mathsf{G}}_{n}^{g,p}(b,c)=\widehat{
\mathsf{F}}_{n}^{g,p}%
(b,c)-\lim
_{c^{+}\downarrow c}\widehat{\mathsf{F}}_{n}^{g,p}
\bigl(b,-c^{+}\bigr). \label{FG}%
\end{equation}
The asymptotic arguments are made on the probability scale $\psi
=\mathsf{G}(c_{\psi})$. When $\mathsf{f}$ is symmetric, the probability scales
of $\mathsf{G}$ and $\mathsf{F}$ are related in a simple linear
fashion, see
(\ref{FGinv}), so that (\ref{FG}) translates into%
%
\begin{equation}
\widehat{\mathsf{G}}_{n}^{g,p}\bigl\{b,\mathsf{G}^{-1}(
\psi)\bigr\}=\widehat {\mathsf{F}%
}_{n}^{g,p}\biggl
\{b,\mathsf{F}^{-1}\biggl(\frac{1+\psi}{2}\biggr)\biggr\}-\lim
_{\psi^{+}%
\downarrow\psi}\widehat{\mathsf{F}}_{n}^{g,p}\biggl
\{b,\mathsf{F}^{-1}\biggl(\frac
{1-\psi^{+}}{2}\biggr)\biggr\}.
\label{FGsym}%
\end{equation}
Therefore, results for $\widehat{\mathsf{F}}_{n}$ transfer to
$\widehat{\mathsf{G}}_{n}$. The corresponding conditional mean process
is
%
\begin{equation}
\overline{\mathsf{G}}_{n}^{g,p}(b,c)=\frac{1}{n}%
{
\sum_{i=1}^{n}} 
g_{in}
\mathsf{E}_{i-1}\bigl\{\varepsilon_{i}^{p}1_{(|\varepsilon_{i}-x_{in}%
^{\prime}b|\leq\sigma c)}
\bigr\},\qquad p=0,1,2. \label{Gbar}%
\end{equation}
Form also the empirical process%
%
\begin{equation}
\mathbb{G}_{n}^{g,p}(b,c)=n^{1/2}\bigl\{\widehat{
\mathsf{G}}_{n}^{g,p}%
(b,c)-\overline{
\mathsf{G}}_{n}^{g,p}(b,c)\bigr\}. \label{Gtilde}%
\end{equation}

For later use note $\mathbb{G}_{n}^{1,0}(0,c)=\mathbb{G}_{n}(c)$. Note also
that $\mathsf{E}_{i-1}\{\varepsilon_{i}^{p}1_{(|\varepsilon_{i}|\leq
\sigma
c)}\}=0$ for odd $p$ since $\mathsf{f}$ is symmetric and $b=0$. Errors in
estimating the quantile are denoted $d=n^{1/2}(c_{\psi}^{b}-c_{\psi})$.
Estimation errors represented by $b,d$ vanish uniformly as shown in the next
result. Due to the two-sidedness of the absolute residuals and symmetry of
$\mathsf{f,}$ only one of the error terms $x_{in}^{\prime}b$ and $n^{-1/2}d$
enters the asymptotic expansion depending on the choice of $p$.

\begin{lemma}
\label{lG}For each $\psi$ let $c_{\psi}=\mathsf{G}^{-1}(\psi)$. Suppose
Assumption~\ref{asFS}\textup{(i)(a), (ii)(b), (ii)(c)} holds for some $0\leq\kappa<\eta
\leq1/4$,
but with $q_{0}=1+2^{r+1}$ only. Then, for all $B,\epsilon>0$ and all
$\omega<\eta-\kappa\leq1/4$,
\begin{longlist}[(b$^{\prime}$)]
\item[(a)] $\sup_{0\leq\psi\leq1}%
\sup_{|b|,|d|\leq n^{1/4-\eta}B}|n^{1/2}\{\overline{\mathsf{G}}_{n}%
^{g,p}(b,c_{\psi}+n^{\kappa-1/2}d)-\overline{\mathsf
{G}}_{n}^{g,p}(0,c_{\psi
})\}-\break 2\sigma^{p-1} c_{\psi}^{p} \mathsf{f}(c_{\psi})n^{-1/2}%
{ \sum_{i=1}^{n}}
g_{in}\{1_{(p\ \mathrm{odd})}x_{in}^{\prime}b+1_{(p\ \mathrm{even})}n^{\kappa-1/2}\sigma
d\}|=\mathrm{O}_{\mathsf{P}}\{n^{2(\kappa-\eta)}\}$;

\item[(b)] $\sup_{0\leq\psi\leq1}\sup_{|b|,|d|\leq n^{1/4-\eta}B}|\mathbb{G}_{n}%
^{g,p}(b,c_{\psi}+n^{\kappa-1/2}d)-\mathbb{G}_{n}^{g,p}(0,c_{\psi
})|=\mathrm{o}_{\mathsf{P}}(1)$;

\item[(b$^{\prime}$)] $\sup_{0\leq\psi
\leq
1}\sup_{|b|,|d|\leq n^{1/4-\eta}B}|\mathbb{G}_{n}^{1,0}(b,c_{\psi}%
+n^{\kappa-1/2}d)-\mathbb{G}_{n}^{1,0}(0,c_{\psi})|=\mathrm{o}_{\mathsf
{P}%
}(n^{-1/8-\omega/2})$;

\item[(c)] $\lim_{\phi\downarrow0}%
{\lim\sup}_{n\rightarrow\infty} \mathsf{P}\{\sup_{0\leq\psi\leq
\psi^{\dag}\leq1\dvt \psi^{\dag}-\psi\leq\phi}|\mathbb
{G}_{n}^{g,p}(0,c_{\psi
^{\dag}})-\mathbb{G}_{n}^{g,p}(0,c_{\psi})|>\epsilon\}\rightarrow0$.
\end{longlist}
\end{lemma}

\begin{pf}
(a) Assumption~\ref{asFS}(i)(a), (ii)(c) implies
Assumption~\ref{as}(i)(b), (iii)(b) with $r=0$, $p\leq2$ and
$g_{in}=1,n^{1/2}%
x_{in}$ or $nx_{in}x_{in}^{\prime}$, and hence the assumptions of Theorem~\ref{tFplin}. First, we want to apply this result to $\overline
{\mathsf{F}%
}_{n}^{g,p}(b,c_{\psi}+n^{\kappa-1/2}d)$. Thus, rewrite
\begin{eqnarray*}
\overline{\mathsf{F}}_{n}^{g,p}\bigl(b,c_{\psi}+n^{\kappa-1/2}d
\bigr) & =&n^{-1}%
{ \sum
_{i=1}^{n}} 
g_{in}
\mathsf{E}_{i-1}\varepsilon_{i}^{p}1_{\{\varepsilon
_{i}-x_{in}^{\prime
}b\leq\sigma(c_{\psi}+n^{\kappa-1/2}d)\}}
\\
& =&n^{-1}%
{ \sum
_{i=1}^{n}} 
g_{in}
\mathsf{E}_{i-1}\varepsilon_{i}^{p}1_{(\varepsilon_{i}-\bar
{x}%
_{in}^{\prime}\bar{b}\leq\sigma c_{\psi})},
\end{eqnarray*}
for $\bar{b}=(b^{\prime},n^{\kappa}d)^{\prime}$ and $\bar{x}%
_{in}=(x_{in}^{\prime},n^{-1/2}\sigma)^{\prime}$, where $|\bar{b}%
|\leq2n^{1/4+\kappa-\eta}B$ while $\bar{x}_{in}$ satisfies Assumption~\ref{as}(iii)(b) because $|\bar{x}_{in}|^{2}=|x_{in}|^{2}+n^{-1}%
\sigma^{2}$. Therefore we find, using that $\overline{\mathsf{G}}_{n}^{g,p}$
can be expressed in terms of $\overline{\mathsf{F}}_{n}^{g,p}$ as in
(\ref{FG}), that $\sigma^{1-p}n^{1/2}\{\overline{\mathsf{G}}_{n}%
^{g,p}(b,c_{\psi}+n^{\kappa-1/2}d)-\overline{\mathsf
{G}}_{n}^{g,p}(0,c_{\psi
})\}$ has correction term
\begin{eqnarray*}
&&c_{\psi}^{p}\mathsf{f}(c_{\psi})n^{-1}%
{
\sum_{i=1}^{n}} 
g_{in}n^{1/2}
\bigl(x_{in}^{\prime}b+n^{\kappa-1/2}\sigma d\bigr)
\\
&&\quad{}-(-c_{\psi})^{p}\mathsf{f}(-c_{\psi})n^{-1}%
{
\sum_{i=1}^{n}} 
g_{in}n^{1/2}
\bigl(x_{in}^{\prime}b-n^{\kappa-1/2}\sigma d\bigr)
\\
&&\qquad=c_{\psi}^{p}\mathsf{f}(c_{\psi})n^{-1/2}%
{
\sum_{i=1}^{n}} 
g_{in}
\bigl[\bigl\{1-(-1)^{p}\bigr\}x_{in}^{\prime}b+\bigl
\{1+(-1)^{p}\bigr\}n^{\kappa
-1/2}\sigma d\bigr],
\end{eqnarray*}
due to the symmetry of $\mathsf{f}$. This reduces as desired.

(b) Let $c_{\psi}^{\dag}=c_{\psi}+n^{\kappa-1/2}d$. Rewrite $\mathcal
{G=}%
\mathbb{G}_{n}^{g,p}(b,c_{\psi}^{\dag})-\mathbb{G}_{n}^{g,p}(0,c_{\psi
})$ as
$\mathcal{G}=\mathcal{G}_{1}+\mathcal{G}_{2}$, where
\[
\mathcal{G}_{1}=\mathbb{G}_{n}^{g,p}
\bigl(b,c_{\psi}^{\dag}\bigr)-\mathbb{G}_{n}%
^{g,p}
\bigl(0,c_{\psi}^{\dag}\bigr),\qquad \mathcal{G}_{2}=
\mathbb{G}_{n}^{g,p}%
\bigl(0,c_{\psi}^{\dag}
\bigr)-\mathbb{G}_{n}^{g,p}(0,c_{\psi}).
\]

The term $\mathcal{G}_{1}$ is $\mathrm{o}_{\mathsf{P}}(1)$ uniformly in
$|b|\leq n^{1/4-\eta}B$, $0\leq\psi\leq1$. To see this, expand $\mathbb
{G}%
_{n}^{g,p}$ in a similar fashion to (\ref{FG}). Apply Theorem~\ref{tFpest},
noting that Assumption~\ref{asFS}(i)(a), (ii)(b), (ii)(c) implies Assumption~\ref{as}(i), (ii), (iii)(a) with $p\leq2$, $g_{in}=1,n^{1/2}x_{in}$ or
$nx_{in}x_{in}^{\prime}$ and the chosen $r$.

The term $\mathcal{G}_{2}$. Apply Theorem~\ref{tFptight} noting that
Assumption~\ref{asFS}(i)(a), (ii)(c) implies Assumption~\ref{as}(i)(a), (iii)(a) with
$r=2$ and some $\nu<1$.
\begin{longlist}[(b$^{\prime}$)]
\item[(b$^{\prime}$)] Similar to (b), but using Theorem~\ref{tFpestomega}.

\item[(c)] Assumption~\ref{asFS}(i)(a), (ii)(c) implies Assumption~\ref{as}(i)(a), (iii)(a).
Apply Theorem~\ref{tFptight}.
\end{longlist}
\end{pf}

\subsection{A first analysis of the order statistics}

The Forward Search is defined in terms of order statistics $\hat{z}^{(m)}$,
see (\ref{xihat}). A process version gives quantiles%
%
\begin{equation}
\hat{c}_{\psi}^{b}=\inf\bigl\{c\dvt \widehat{
\mathsf{G}}_{n}^{1,0}(b,c)\geq\psi\bigr\}.
\label{chatb}%
\end{equation}
Setting $b=0$ gives $\hat{c}_{\psi}^{0}=\widehat{\mathsf
{G}}_{n}^{-1}(\psi)$
as defined in (\ref{Ginv0}) and studied in Lemma~\ref{tdhat0}.
Evaluating the
empirical distribution function at the quantile gives%
%
\begin{equation}
\widehat{\mathsf{G}}_{n}^{1,0}\bigl(b,\hat{c}_{\psi}^{b}
\bigr)=\frac{1}{n}\inf (x\in\mathbb{N}_{0}\dvt x\geq\psi n).
\label{Ghatb}%
\end{equation}
The first result gives an algebraic bound to the distance between $\hat
{c}_{\psi}^{b}$ and $\hat{c}_{\psi}^{0}$. Probabilistic bounds follow.

\begin{lemma}
\label{tdhatalg}For all $b,\psi$, the quantiles $\hat{c}_{\psi}^{b}$ and
$\hat{c}_{\psi}^{0}$ satisfy $\sigma|\hat{c}_{\psi}^{b}-\hat{c}_{\psi}%
^{0}|<2|b|\max_{1\leq i\leq n}|x_{in}|$.
\end{lemma}

\begin{pf}
1. \textitt{A property of}
$\widehat{\mathsf{G}}_{n}$. The quantile $\sigma\hat{c}_{\psi}^{0}$
is the left-continuous inverse of the right-continuous function
$\widehat{\mathsf{G}}_{n}^{1,0}(0,c)=\widehat{\mathsf{G}}_{n}(c)$ in
(\ref{Ginv0}). Thus,%
%
\begin{equation}
\widehat{\mathsf{G}}_{n}(y)<\widehat{\mathsf{G}}_{n}\bigl(
\hat{c}_{\psi}^{0}%
\bigr)\leq\widehat{
\mathsf{G}}_{n}(z) \quad\Rightarrow \quad y<\hat{c}_{\psi}
^{0}\leq z. \label{ptdhatalgineq1}%
\end{equation}

2. \textitt{A lower bound}. Let $x_{\max}=\max_{1\leq i\leq n}|x_{in}|$. Then
it follows that
\[
\mathcal{S}_{i}=\bigl[-\sigma\hat{c}_{\psi}^{b}+x_{in}^{\prime}b,
\sigma\hat {c}_{\psi}^{b}+x_{in}^{\prime}b
\bigr]\subset{}\bigl[-\sigma\hat{c}_{\psi}%
^{b}-x_{\max}|b|,
\sigma\hat{c}_{\psi}^{b}+x_{\max}|b|\bigr]=
\mathcal{S},
\]
so that for all $0\leq\psi\leq1$ and $z=\hat{c}_{\psi}^{b}+\sigma
^{-1}x_{\max
}|b|$,
\[
\widehat{\mathsf{G}}_{n}^{1,0}\bigl(b,\hat{c}_{\psi}^{b}
\bigr)\leq\frac{1}{n}%
{ \sum
_{i=1}^{n}} 
1_{(|\varepsilon_{i}|\leq\sigma z)}=\widehat{
\mathsf{G}}_{n}^{1,0}%
(0,z)=\widehat{
\mathsf{G}}_{n}(z).
\]
Using (\ref{Ghatb}) we find, for all $b,\psi$, that
\[
0=\widehat{\mathsf{G}}_{n}^{1,0}\bigl(b,
\hat{c}_{\psi}^{b}\bigr)-\widehat{\mathsf
{G}%
}_{n}^{1,0}\bigl(0,\hat{c}_{\psi}^{0}
\bigr)\leq\widehat{\mathsf{G}}_{n}%
(z)-\widehat{
\mathsf{G}}_{n}\bigl(\hat{c}_{\psi}^{0}\bigr),
\]
which implies that $\sigma z=\sigma\hat{c}_{\psi}^{b}+x_{\max}|b|\geq
\sigma\hat{c}_{\psi}^{0}$ by inequality (\ref{ptdhatalgineq1}).

3. \textitt{An upper bound}. For $y=\hat{c}_{\psi}^{b}-\sigma
^{-1}2x_{\max
}|b|$, we find
\[
\mathcal{S}_{i}=\bigl[-\sigma\hat{c}_{\psi}^{b}+x_{in}^{\prime}b,
\sigma\hat {c}_{\psi}^{b}+x_{in}^{\prime}b
\bigr]\supset{}[-\sigma y,\sigma y]=\mathcal{S}%
,
\]
noting that the smaller set is empty if $y<0$. It therefore follows that
\[
\widehat{\mathsf{G}}_{n}^{1,0}\bigl(b,\hat{c}_{\psi}^{b}
\bigr)\geq\frac{1}{n}%
{ \sum
_{i=1}^{n}} 
1_{(|\varepsilon_{i}|\leq\sigma y)}=\widehat{
\mathsf{G}}_{n}(y).
\]
Actually, this inequality must be strict. Indeed, at least one $i^{\dag}$
exists for which $\sigma\hat{c}_{\psi}^{b}=|\varepsilon_{i^{\dag
}}-x_{i^{\dag
}n}^{\prime}b|$. For this (these) $i^{\dag}$ it holds that $\varepsilon
_{i^{\dag}}\in\mathcal{S}_{i}$ but $\varepsilon_{i^{\dag}}\notin
\mathcal{S}$. Thus, $\widehat{\mathsf{G}}_{n}^{1,0}(b,\hat{c}_{\psi}%
^{b})>\widehat{\mathsf{G}}_{n}(y)$. Proceed as before to see that%
%
\begin{equation}
0=\widehat{\mathsf{G}}_{n}^{1,0}\bigl(b,
\hat{c}_{\psi}^{b}\bigr)-\widehat{\mathsf
{G}%
}_{n}^{1,0}\bigl(0,\hat{c}_{\psi}^{0}
\bigr)>\widehat{\mathsf{G}}_{n}%
(y)-\widehat{
\mathsf{G}}_{n}\bigl(\hat{c}_{\psi}^{0}\bigr),
\label
{ptdhatalglower}%
\end{equation}
which implies that $y=\hat{c}_{\psi}^{b}-\sigma^{-1}2x_{\max}|b|<\hat{c}
_{\psi}^{0}$ by inequality (\ref{ptdhatalgineq1}).
\end{pf}

The next result introduces a convergence rate for $\hat{c}_{\psi
}^{b}-\hat
{c}_{\psi}^{0}$.

\begin{lemma}
\label{tdhatf}Suppose Assumption~\ref{asFS}\textup{(i)(a), (ii)(b), (ii)(c)} holds, but with
$q_{0}=1+2^{r+1}$ only. Then, for all $\omega<\eta-\kappa$,
\[
\sup_{0\leq\psi\leq1}\sup_{|b|\leq n^{1/4-\eta}B}n^{1/2}%
\bigl|
\mathsf{f}\bigl(\hat{c}_{\psi}^{0}\bigr) \bigl(
\hat{c}_{\psi}^{b}-\hat{c}_{\psi}%
^{0}
\bigr)\bigr|=\mathrm{o}_{\mathsf{P}}\bigl(n^{-\omega}\bigr).
\]
\end{lemma}

\begin{pf}
If we combine Lemma~\ref{tdhatalg} with
Assumption~\ref{asFS}(ii)(b) we find that\break $\max_{1\leq i\leq n}%
|x_{in}|=\mathrm{O}_{\mathsf{P}}(n^{\kappa-1/2})$ to get that $\hat
{c}_{\psi
}^{b}-\hat{c}_{\psi}^{0}=\mathrm{O}_{\mathsf{P}}(n^{-1/4+\kappa-\eta})$ for
$|b|\leq n^{1/4-\eta}B$. Thus, for any $\epsilon>0$ a $C>0$ exists so
that the
set $\mathcal{C}_{n}=\{|n^{1/2-\kappa}(\hat{c}_{\psi}^{b}-\hat{c}_{\psi}
^{0})|\leq n^{1/4-\eta}C\}$ has probability $\mathsf{P}(\mathcal{C}%
_{n})>1-\epsilon$. On this set it holds, using (\ref{Ghatb}) and with
$d=n^{1/2-\kappa}(\hat{c}_{\psi}^{b}-\hat{c}_{\psi}^{0})$, that
\[
0=\widehat{\mathsf{G}}_{n}^{1,0}\bigl(b,
\hat{c}_{\psi}^{b}\bigr)-\widehat{\mathsf
{G}%
}_{n}^{1,0}\bigl(0,\hat{c}_{\psi}^{0}
\bigr)=\widehat{\mathsf{G}}_{n}^{1,0}\bigl(b,\hat
{c}_{\psi}^{0}+n^{\kappa-1/2}d\bigr)-\widehat{
\mathsf{G}}_{n}^{1,0}\bigl(0,\hat{c} 
_{\psi}^{0}
\bigr).
\]
Lemma~\ref{lG}(a), using Assumption~\ref{asFS}(i)(a), (ii)(c), shows that
\[
n^{1/2}\bigl\{\overline{\mathsf{G}}_{n}^{1,0}
\bigl(b,c_{\psi}+n^{\kappa-1/2}%
d\bigr)-\overline{
\mathsf{G}}_{n}^{1,0}(0,c_{\psi})\bigr\}-2
\sigma^{-1}\mathsf{f}%
(c_{\psi})n^{\kappa}
\sigma d=\mathrm{O}_{\mathsf{P}}\bigl(n^{2\kappa-2\eta
}\bigr)=
\mathrm{o}_{\mathsf{P}}\bigl(n^{-\omega}\bigr),
\]
uniformly in $0\leq\psi\leq1$ and $|b|,|d|\leq n^{1/4-\eta}B$, for all
$\omega<\eta-\kappa<2(\eta-\kappa)$. Lemma~\ref{lG}(b$^{\prime}$) using
Assumption~\ref{asFS}(i)(a), (ii)(b), (ii)(c) shows that, uniformly in $0\leq\psi
\leq1$
and $|b|,|d|\leq n^{1}/4-\eta B$,
\[
\mathbb{G}_{n}^{1,0}\bigl(b,c_{\psi}+n^{\kappa-1/2}d
\bigr)-\mathbb{G}_{n}^{1,0}%
(0,c_{\psi})=
\mathrm{o}_{\mathsf{P}}\bigl(n^{-\omega}\bigr),
\]
for all $\omega<\eta-\kappa$. Using the definition $\mathbb{G}_{n}%
^{1,0}=n^{1/2}(\widehat{\mathsf{G}}_{n}^{1,0}-\overline{\mathsf{G}}_{n}%
^{1,0})$,
\[
0=n^{1/2}\bigl\{\widehat{\mathsf{G}}_{n}^{1,0}
\bigl(b,\hat{c}_{\psi}^{0}+n^{\kappa
-1/2}d\bigr)-\widehat{
\mathsf{G}}_{n}^{1,0}\bigl(0,\hat{c}_{\psi}^{0}
\bigr)\bigr\}=2\mathsf {f}%
\bigl(\hat{c}_{\psi}^{0}
\bigr)n^{\kappa}d+\mathrm{o}_{\mathsf{P}}\bigl(n^{-\omega}\bigr).
\]
Inserting $d=n^{1/2-\kappa}(\hat{c}_{\psi}^{b}-\hat{c}_{\psi}^{0})$ we
get the
desired result.
\end{pf}

The next result provides a modification of
Cs\"{o}rg\H{o} \cite{Cs}, equation (2.8).

\begin{lemma}
\label{lf2f}Let $c_{\psi}=\mathsf{G}^{-1}(\psi)$. Suppose $\mathsf{f}$ is
symmetric and decreasing for large $c$ and that Assumption~\ref{asFS}\textup{(i)(b)}
holds, but with $q_{0}=1+2^{r+1}$ only. Let $|\psi^{\ast}-\psi|\leq
|\mathsf{G}(\hat{c}_{\psi}^{0})-\psi|$, then:
\begin{longlist}[(a)]
\item[(a)] $\sup_{0\leq
\psi
\leq1-c_{n}}|1-\mathsf{f}(c_{\psi})/\mathsf{f}(c_{\psi^{\ast}})|=\mathrm
{o}%
_{\mathsf{P}}(1)$, for any sequence $c_{n}\rightarrow0$ for which
$nc_{n}\rightarrow\infty$;

\item[(b)] $\sup_{0\leq\psi\leq n/(n+1)}%
|1-\mathsf{f}(c_{\psi})/\mathsf{f}(c_{\psi^{\ast}})|=\mathrm{O}_{\mathsf
{P}%
}(1)$.
\end{longlist}
\end{lemma}

\begin{pf}
(a) By (\ref{FGinv}) $\mathsf{G}^{-1}%
(\psi)=\mathsf{F}^{-1}(y)$ for $y=(1+\psi)/2$ varying in $1/2\leq
y\leq1-(2n+2)^{-1}$. Let $\gamma=\sup_{c\in\mathbb{R}}\mathsf{F}%
(c)\{1-\mathsf{F}(c)\}|\mathsf{\dot{f}}(c)|/\{\mathsf{f}(c)\}^{2}$
which is
finite by Assumption~\ref{asFS}(i)(b). It is first argued that for all
$\epsilon>0$ and $0<c<1$ and all $n$
%
\begin{equation}
\mathsf{P}\biggl\{\sup_{1/2+c\leq y\leq1-c}\biggl|\frac{\mathsf{f}\{\mathsf{F}^{-1}%
(y)\}}{\mathsf{f}\{\mathsf{F}^{-1}(y^{\ast})\}}-1\biggr|>\epsilon
\biggr\}\leq 4\bigl\{1+\operatorname{int}(\gamma)\bigr\} \bigl\{\exp(-nch_{1})+
\exp(-nch_{2}%
)\bigr\},\label{CsorgoTh151}%
\end{equation}
where, with $h(\lambda)=\lambda+\log(1/\lambda)-1$,
\begin{eqnarray*}
h_{1} & =&h\bigl[(1+\epsilon)^{\{1+\operatorname{int}(\gamma)\}/2}\bigr],
\\
h_{2} & =&h\bigl[1/(1+\epsilon)^{\{1+\operatorname{int}(\gamma)\}/2}\bigr].
\end{eqnarray*}
This is nearly the statement of Cs\"{o}rg\H{o} \cite{Cs}, Theorem~1.5.1, which,
however, has the denominator $\mathsf{f}(\theta_{y,n})$ instead of
$\mathsf{f}\{\widehat{\mathsf{F}}_{n}^{-1}(y^{\ast})\}$ where $\theta_{y,n}$
is a particular intermediate point between $\widehat{\mathsf{F}}_{n}^{-1}(y)$
and $\mathsf{F}^{-1}(y)$ rather than any intermediate point. Cs\"{o}rg\H{o}
states that the proof of this theorem is similar to that of his Theorem~1.4.3.
Equation (1.4.18.2) of that proof uses a bound only depending on
$\widehat{\mathsf{F}}_{n}^{-1}(y)$ and $\mathsf{F}^{-1}(y)$ and not on the
particular intermediate point $\theta_{y,n}$. This proves (\ref{CsorgoTh151}).

The inequality (\ref{CsorgoTh151}) implies that for any sequence
$c_{n}\rightarrow0$ for which $nc_{n}\rightarrow\infty$,
\[
\mathsf{P} \biggl\{ \sup_{1/2+c_{n}\leq y\leq1-c_{n}}\biggl|\frac{\mathsf{f}%
\{\mathsf{F}^{-1}(y)\}}{\mathsf{f}\{\widehat{\mathsf
{F}}_{n}^{-1}(y^{\ast}%
)\}}-1\biggr|>\epsilon
\biggr\} \rightarrow0.
\]
The reason is that $h(\lambda)>0$ for all $\lambda>0$ so $\lambda\neq1$.
Consider the tails.

\textitt{Left-hand tail}. Use that $c_{n}$ vanishes, that $\mathsf
{G}(\hat
{c}_{\psi}^{0})-\psi=\mathrm{O}_{\mathsf{P}}(n^{-1/2})$ by Lemmas
\ref{lGhat0}, \ref{tdhat0}, and that $\mathsf{f}$ is uniformly
continuous in
a neighbourhood of zero because $\mathsf{f}$ is bounded, positive and
continuous.

(b) \textitt{Right-hand tail}. It suffices to argue that
%
\begin{equation}
\lim_{\epsilon\rightarrow\infty}\mathop{\lim\sup }_{n\rightarrow\infty} \mathsf{P} \biggl\{
\sup_{1-c_{n}\leq y\leq1-(2n+2)^{-1}}\biggl|\frac
{\mathsf{f}\{\mathsf{F}^{-1}(y)\}}{\mathsf{f}\{\widehat{\mathsf{F}}_{n}%
^{-1}(y^{\ast})\}}-1\biggr|>\epsilon \biggr\} =0.
\label{plf2fCsorgo28}%
\end{equation}
Apply the inequality (\ref{CsorgoTh151}) with $c=(2n+2)^{-1}$ so that
$nc\sim1/2$. Then use that $h_{1},h_{2}\rightarrow\infty$ for $\epsilon
\rightarrow\infty$ since $h(\lambda)\rightarrow\infty$ for $\lambda
\rightarrow\infty$.
\end{pf}

The next result relates $\hat{c}_{\psi}^{0}$ to $c_{\psi}$.

\begin{lemma}
\label{tdhat0hat}Suppose Assumption~\ref{asFS}\textup{(i)(a), (i)(b)} holds with $q=1$
only. Then
\[
\sup_{0\leq\psi\leq1}\bigl|\bigl(\hat{c}_{\psi}^{0}
\bigr)^{k}\mathsf{f}\bigl(\hat{c}_{\psi}%
^{0}
\bigr)-(c_{\psi})^{k}\mathsf{f}(c_{\psi})\bigr|=
\mathrm{o}_{\mathsf{P}}%
(1) \qquad\mbox{for }k=0,1.
\]
\end{lemma}

\begin{pf}
1. \textitt{Consider} $\psi$
\textitt{so that} $0\leq\psi\leq1-1/z_{n}$ \textitt{for any sequence} $0<z_{n}<\mathrm{o}%
(n^{1/2})$. Rewrite the process of interest as
%
\begin{equation}
\bigl(\hat{c}_{\psi}^{0}\bigr)^{k}\mathsf{f}
\bigl(\hat{c}_{\psi}^{0}\bigr)-(c_{\psi}%
)^{k}
\mathsf{f}(c_{\psi})=\bigl\{\bigl(\hat{c}_{\psi}^{0}
\bigr)^{k}-(c_{\psi})^{k}%
\bigr\}
\mathsf{f}(c_{\psi})+\bigl(\hat{c}_{\psi}^{0}
\bigr)^{k}\mathsf{f}\bigl(\hat{c}_{\psi} 
^{0}
\bigr)\biggl\{1-\frac{\mathsf{f}(c_{\psi})}{\mathsf{f}(\hat{c}_{\psi}^{0})}\biggr\}. \label{ptdhat0hatrewrite}%
\end{equation}
The first term is zero for $k=0$. For $k=1$, $n^{1/2}(\hat{c}_{\psi}%
^{0}-c_{\psi})\mathsf{f}(c_{\psi})=-\widehat{\mathbb
{G}}_{n}^{1,0}(0,c_{\psi
})+\mathrm{o}_{\mathsf{P}}(1)$ by Lemma~\ref{tdhat0}(a) using Assumption~\ref{asFS}(i)(b).
This in turn is tight due to Lemma~\ref{lGhat0}. Overall,
the first term is $\mathrm{O}_{\mathsf{P}}(n^{-1/2})$. For the second term,
note that $(\hat{c}_{\psi}^{0})^{k}\mathsf{f}(\hat{c}_{\psi}^{0})$ is bounded
uniformly in $0\leq\psi\leq1$ due to Assumption~\ref{asFS}(i)(a) with $q=1$,
while $1-\mathsf{f}(c_{\psi})/\mathsf{f}(\hat{c}_{\psi}^{0})$ vanishes by
Lemma~\ref{lf2f}(a) using Assumption~\ref{asFS}(i)(b).

2. \textitt{Consider} $\psi$ \textitt{so that} $\psi_{n}\leq\psi\leq1$
\textitt{for any sequence} $\psi_{n}\rightarrow1$. Assumption~\ref{asFS}(i)(a)
and the continuity of $\mathsf{f}$ implies that $(c_{\psi})^{k}\mathsf
{f}%
(c_{\psi})$ is continuous and convergent for $\psi\rightarrow1$, and
hence for
$c_{\psi}\rightarrow G^{-1}(1)$. Rewrite
\[
\hat{c}_{\psi_{n}}^{0}=\mathsf{G}^{-1}\bigl\{
\mathsf{G}\bigl(\hat{c}_{\psi_{n}}%
^{0}\bigr)\bigr\}=
\mathsf{G}^{-1}\bigl[\psi_{n}+\bigl\{\mathsf{G}\bigl(
\hat{c}_{\psi
_{n}}^{0}\bigr)-\psi _{n}\bigr\}\bigr]
\geq\mathsf{G}^{-1}(\psi_{n}-g_{n}),
\]
where $g_{n}=\sup_{0\leq\psi\leq1}\{\mathsf{G}(\hat{c}_{\psi}^{0}%
)-\psi\}=\mathrm{O}_{\mathsf{P}}(n^{-1/2})$ due to Lemmas \ref{lGhat0},
\ref{tdhat0}(c) using Assumption~\ref{asFS}(i)(b). By the
continuity of
$\mathsf{G}^{-1}$, $\hat{c}_{\psi_{n}}^{0}\rightarrow\mathsf
{G}^{-1}(1)$ in
probability and therefore $(\hat{c}_{\psi}^{0})^{k}\mathsf{f}(\hat
{c}_{\psi
}^{0})$ and $(c_{\psi})^{k}\mathsf{f}(c_{\psi})$ converge to the same
limit in
probability, and their difference vanishes.
\end{pf}

\subsection{A one-step result for the least squares estimator}

A one-step result for the least squares estimator now follows. Equation
(\ref{betahatG}) represents the one-step least squares estimator $\hat
{\beta
}^{(m+1)}$ in terms of $\widehat{\mathsf{G}}_{n}^{g,p}$. That
expression has
the random quantities $\hat{b}^{(m)}$ and $\sigma^{-1}\hat{z}^{(m)}$ as
arguments. Replacing these by a deterministic quantity $b$ and the residual
$\hat{c}_{\psi}^{b}$ defined in (\ref{chatb}) gives the following asymptotic
uniform linearization result.

\begin{lemma}
\label{tbhat}Let $c_{\psi}=\mathsf{G}^{-1}(\psi)$ and%
%
\begin{equation}
\rho_{\psi}=2c_{\psi}\mathsf{f}(c_{\psi})/\psi.
\label{rho}%
\end{equation}
Suppose Assumption~\ref{asFS}\textup{(i)(a), (i)(b), (ii)} holds, while $\psi_{0}>0$ and
$\eta\leq1/4$, but with $q_{0}=1+2^{r+1}$ only. Then:
\begin{longlist}[(a)]
\item[(a)] $\sup_{0\leq\psi\leq1,|b|\leq n^{1/4-\eta}B}|n^{1/2}\widehat{\mathsf{G}}_{n}%
^{x,1}(b,\hat{c}_{\psi}^{b})-\mathbb{G}_{n}^{x,1}(0,c_{\psi})-2c_{\psi
}\mathsf{f}(c_{\psi})\Sigma_{n}b|=\mathrm{o}_{\mathsf{P}}(1)$;

\item[(b)]
$\sup_{0\leq\psi\leq1,|b|\leq n^{1/4-\eta}B}|n^{1/2}\{\widehat{\mathsf
{G}}%
_{n}^{xx,0}(b,\hat{c}_{\psi}^{b})-\Sigma_{n}\psi\}|=\mathrm{O}_{\mathsf
{P}%
}(1)$;

\item[(c)] $\sup_{\psi_{0}\leq\psi\leq1,|b|\leq n^{1/4-\eta}%
B}|\{\widehat{\mathsf{G}}_{n}^{xx,0}(b,\hat{c}_{\psi}^{b})\}^{-1}%
n^{1/2}\widehat{\mathsf{G}}_{n}^{x,1}(b,\hat{c}_{\psi}^{b})-(\psi\Sigma
_{n})^{-1}\mathbb{G}_{n}^{x,1}(0,c_{\psi})-\break \rho_{\psi}b|=\mathrm{o}%
_{\mathsf{P}}(1)$.
\end{longlist}
\end{lemma}

\begin{pf}
(a) The inequality of Lemma~\ref{tdhatalg}
implies that%
%
\begin{equation}
\sup_{0\leq\psi\leq1}\sup_{|b|\leq n^{1/4-\eta}B}n^{1/2-\kappa}\bigl|\hat
{c}_{\psi
}^{b}-\hat{c}_{\psi}^{0}\bigr|=
\mathrm{O}_{\mathsf{P}}\bigl(n^{1/4-\eta}\bigr), \label{ptbhatorderstat}%
\end{equation}
since $\max_{1\leq i\leq n}|x_{in}|=\mathrm{O}_{\mathsf{P}}(n^{\kappa-1/2})$
by Assumption~\ref{asFS}(ii)(b), where $0\leq\kappa<\eta\leq1/4$. By
definition
\[
n^{1/2}\widehat{\mathsf{G}}_{n}^{x,1}
\bigl(b,c_{\psi}+n^{\kappa
-1/2}d\bigr)=\mathbb{G}%
_{n}^{x,1}
\bigl(b,c_{\psi}+n^{\kappa-1/2}d\bigr)+n^{1/2}\overline{
\mathsf{G}}_{n}%
^{x,1}\bigl(b,c_{\psi}+n^{\kappa-1/2}d
\bigr).
\]
Lemma~\ref{lG}(a), (b), using Assumption~\ref{asFS}(i)(a), (ii)(b), (ii)(c)
along with
the definitions $g_{in}=n^{1/2}x_{in}$ and $\Sigma_{n}=%
{ \sum_{i=1}^{n}}
x_{in}x_{in}^{\prime}$ gives, uniformly in $|b|,|d|\leq n^{1/4-\eta}B$ and
$0\leq\psi\leq1$,
\[
n^{1/2}\widehat{\mathsf{G}}_{n}^{x,1}
\bigl(b,c_{\psi}+n^{\kappa
-1/2}d\bigr)=\mathbb{G}%
_{n}^{x,1}(0,c_{\psi})+n^{1/2}
\overline{\mathsf{G}}_{n}^{x,1}(0,c_{\psi
})+2c_{\psi}
\mathsf{f}(c_{\psi})\Sigma_{n}b+\mathrm{o}_{\mathsf{P}}(1).
\]
Note that $\overline{\mathsf{G}}_{n}^{x,1}(0,c_{\psi})=0$ due to the symmetry
of $\mathsf{f}$. Replace $c_{\psi}$ by $\hat{c}_{\psi}^{0}$ and $d$ by
$n^{1/2-\kappa}(\hat{c}_{\psi}^{b}-\hat{c}_{\psi}^{0})$, which is
$\mathrm{O}_{\mathsf{P}}(n^{1/4-\eta})$ due to (\ref{ptbhatorderstat}).
Thus, it holds on a set with large probability that%
%
\begin{equation}
n^{1/2}\widehat{\mathsf{G}}_{n}^{x,1}\bigl(b,
\hat{c}_{\psi}^{b}\bigr)=\mathbb{G}%
_{n}^{x,1}
\bigl(0,\hat{c}_{\psi}^{0}\bigr)+2\hat{c}_{\psi}^{0}
\mathsf{f}\bigl(\hat {c}_{\psi
}^{0}\bigr)
\Sigma_{n}b+\mathrm{o}_{\mathsf{P}}(1), \label{ptbhataexpand}%
\end{equation}
uniformly in $|b|\leq n^{1/4-\eta}B$ and $0\leq\psi\leq1$. The two
terms are
analysed in turn.

\textitt{First term}. Let $a_{\psi}=n^{1/2}\{\mathsf{G}(\hat{c}_{\psi}%
^{0})-\psi\}$. Theorem~\ref{tdhat0}(c) using Assumption~\ref{asFS}(i)(b)
shows that $a_{\psi}=-\mathbb{G}_{n}(c_{\psi})+\mathrm{o}_{\mathsf{P}}(1)$
uniformly in $0\leq\psi\leq1$, which in turn is tight due to Lemma~\ref{lGhat0}. Expand
%
\begin{equation}
\hat{c}_{\psi}^{0}=\mathsf{G}^{-1}\bigl\{
\mathsf{G}\bigl(\hat{c}_{\psi}^{0}%
\bigr)\bigr
\}=c_{\mathsf{G}(\hat{c}_{\psi}^{0})}=c_{\psi+n^{-1/2}a_{\psi}}. \label{ptbhatchat}%
\end{equation}
Lemma~\ref{lG}(c) using Assumption~\ref{asFS}(i)(a), (ii)(b), (ii)(c) shows
$\mathbb{G}_{n}^{x,1}(0,\hat{c}_{\psi}^{0})=\mathbb
{G}_{n}^{x,1}(0,c_{\psi
})+\mathrm{o}_{\mathsf{P}}(1)$.

\textitt{Second term.} Use that $\hat{c}_{\psi}^{0}\mathsf{f}(\hat
{c}_{\psi
}^{0})=c_{\psi}\mathsf{f}(c_{\psi})+\mathrm{o}_{\mathsf{P}}(1)$
uniformly in
$\psi$ by Lemma~\ref{tdhat0hat} using Assumption~\ref{asFS}(i)(a), (i)(b).

(b) An expansion as in (\ref{ptbhataexpand}) gives
\[
\widehat{\mathsf{G}}_{n}^{xx,0}\bigl(b,\hat{c}_{\psi}^{b}
\bigr)=n^{-1/2}\mathbb{G} 
_{n}^{xx,0}
\bigl(0,\hat{c}_{\psi}^{0}\bigr)+\overline{\mathsf
{G}}_{n}^{xx,0}\bigl(0,\hat {c}_{\psi}^{0}
\bigr)+2\mathsf{f}\bigl(\hat{c}_{\psi}^{0}\bigr)
\Sigma_{n}\bigl(\hat{c}_{\psi} 
^{b}-
\hat{c}_{\psi}^{0}\bigr)+\mathrm{o}_{\mathsf{P}}
\bigl(n^{-1/2}\bigr),
\]
uniformly in $b,\psi$. The three terms are analysed in turn.

\textitt{First term.} This is $n^{-1/2}\mathbb{G}_{n}^{xx,0}(0,\hat
{c}_{\psi
}^{0})=n^{-1/2}\mathbb{G}_{n}^{xx,0}(0,c_{\psi})+\mathrm{o}_{\mathsf{P}%
}(n^{-1/2})$ by an argument as for the first term of (\ref{ptbhataexpand}).

\textitt{Second term.} Use that $\Sigma_{n}=%
{ \sum_{i=1}^{n}}
x_{in}x_{in}^{\prime}$ is tight by Assumption~\ref{asFS}(ii)(a), while
$\mathsf{G}(\hat{c}_{\psi}^{0})=\psi+\mathrm{O}_{\mathsf{P}}(n^{-1/2})$
uniformly in $\psi$ by Lemma~\ref{lGhat0}, \ref{tdhat0}(c) using
Assumption~\ref{asFS}(i)(b). Thus,
\[
\overline{\mathsf{G}}_{n}^{xx,0}\bigl(0,
\hat{c}_{\psi}^{0}\bigr)=\frac{1}{n}%
{
\sum_{i=1}^{n}} 
nx_{in}x_{in}^{\prime}
\mathsf{E}_{i-1}1_{(|\varepsilon_{i}|\leq\sigma
\hat
{c}_{\psi}^{0})}=\Sigma_{n}\mathsf{G}\bigl(
\hat{c}_{\psi}^{0}\bigr)=\Sigma_{n}%
\psi+\mathrm{O}_{\mathsf{P}}\bigl(n^{-1/2}\bigr).
\]

\textitt{Third term.} This is $\mathrm{o}_{\mathsf{P}}(n^{-1/2})$ since
$\mathsf{f}(\hat{c}_{\psi}^{0})(\hat{c}_{\psi}^{b}-\hat{c}_{\psi}%
^{0})=\mathrm{o}_{\mathsf{P}}(n^{-1/2})$ uniformly in $\psi,b$ by Lemma~\ref{tdhatf} using
Assumption~\ref{asFS}(i)(a), (ii)(a), (ii)(b), while $\Sigma_{n}$
is tight by Assumption~\ref{asFS}(ii)(a).

(c) Combine (a), (b). The denominator from (b)
satisfies
\[
\widehat{\mathsf{G}}_{n}^{xx,0}\bigl(b,\hat{c}_{\psi}^{b}
\bigr)=\psi\Sigma _{n}\bigl\{1+\mathrm{o}_{\mathsf{P}}(1)\bigr\},
\]
for $\psi\geq\psi_{0}>0$ since $\Sigma_{n}\rightarrow\Sigma$ in distribution
where $\Sigma>0$ a.s. by Assumption~\ref{asFS}(ii)(a). Combine with the
expression for the numerator in (a).
\end{pf}

For the variance estimator, expansions of the same type are needed.

\begin{lemma}
\label{lshat}Suppose Assumption~\ref{asFS}\textup{(i)(a), (i)(b), (i)(d), (ii)} holds while
$\psi_{0}>0$ and $\eta\leq1/4$. Then:
\begin{longlist}[(a)]
\item[(a)] $\sup_{\psi_{0}\leq
\psi
\leq1,|b|\leq n^{1/4-\eta}B}|\{\widehat{\mathsf{G}}_{n}^{x,1}(b,\hat
{c}_{\psi
}^{b})\}^{\prime}\{\widehat{\mathsf{G}}_{n}^{xx,0}(b,\hat{c}_{\psi}%
^{b})\}^{-1}\{\widehat{\mathsf{G}}_{n}^{x,1}(b,\hat{c}_{\psi}^{b}%
)\}|=\mathrm{O}_{\mathsf{P}}(n^{-1/2-2\eta})$;

\item[(b)] $\sup_{\psi
_{0}\leq\psi\leq n/(n+1),|b|\leq n^{1/4-\eta}B}|n^{1/2}\{\widehat
{\mathsf{G}%
}_{n}^{1,2}(b,\hat{c}_{\psi}^{b})-\tau_{\psi}\sigma^{2}\}-\mathbb{G}_{n}
^{1,2}(0,c_{\psi})+\sigma^{2}c_{\psi}^{2}\mathbb{G}_{n}^{1,0}(c_{\psi
})|=\mathrm{o}_{\mathsf{P}}(1)$.
\end{longlist}
\end{lemma}

\begin{pf}
(a) Lemma~\ref{tbhat}(a), (c) using Assumption~\ref{asFS}(i)(a), (i)(b), (ii) shows%
%
\begin{eqnarray}
n^{1/2}\widehat{\mathsf{G}}_{n}^{x,1}\bigl(b,
\hat{c}_{\psi}^{b}\bigr) & =&\mathbb {G}%
_{n}^{x,1}(0,c_{\psi})+2c_{\psi}
\mathsf{f}(c_{\psi})\Sigma_{n}b+\mathrm {o}%
_{\mathsf{P}}(1),\label{ptshatexpandGx,1}
\\
\bigl\{\widehat{\mathsf{G}}_{n}^{xx,0}\bigl(b,
\hat{c}_{\psi}^{b}\bigr)\bigr\}^{-1}n^{1/2}%
\widehat{\mathsf{G}}_{n}^{x,1}\bigl(b,\hat{c}_{\psi}^{b}
\bigr) & =&(\Sigma_{n}\psi )^{-1}\mathbb{G}_{n}^{x,1}(0,c_{\psi})-
\rho_{\psi}b+\mathrm{o}_{\mathsf
{P}%
}(1),\label{ptshatexpandGxx,0}%
\end{eqnarray}
uniformly in $|b|\leq n^{1/4-\eta}B$, $\psi_{0}\leq\psi\leq1$ for $\eta
\leq1/4$. Because $\mathbb{G}_{n}^{x,1}(0,c_{\psi})$ is tight by Lemma~\ref{lG}(c) using
Assumption~\ref{asFS}(i)(a), (ii)(b), (ii)(c), since $\Sigma
_{n}\rightarrow\Sigma$ in distribution where $\Sigma>0$ a.s. by Assumption~\ref{asFS}(ii)(a) and since $|b|\leq n^{1/4-\eta}B$, then both
$\widehat{\mathsf{G}}_{n}^{x,1}(b,\hat{c}_{\psi}^{b})$, see
(\ref{ptshatexpandGx,1}), and $\{\widehat{\mathsf
{G}}_{n}^{xx,0}(b,\hat
{c}_{\psi}^{b})\}^{-1}\widehat{\mathsf{G}}_{n}^{x,1}(b,\hat{c}_{\psi}^{b})$,
see (\ref{ptshatexpandGxx,0}), are $\mathrm{O}_{\mathsf
{P}}(n^{-1/4-\eta
})$. Thus, their product is $\mathrm{O}_{\mathsf{P}}(n^{-1/2-2\eta})$
as desired.

(b) The argument relates to that of the proof of Lemma~\ref{tbhat}.

1. \textitt{Expansion.} By definition
\[
n^{1/2}\widehat{\mathsf{G}}_{n}^{1,2}
\bigl(b,c_{\psi}+n^{\kappa
-1/2}d\bigr)=\mathbb{G}%
_{n}^{1,2}
\bigl(b,c_{\psi}+n^{\kappa-1/2}d\bigr)+n^{1/2}\overline{
\mathsf{G}}_{n}%
^{1,2}\bigl(b,c_{\psi}+n^{\kappa-1/2}d
\bigr).
\]
Apply Lemma~\ref{lG}(a), (b) using Assumption~\ref{asFS}(i)(a), (ii)(b), (ii)(c)
to get
\[
n^{1/2}\widehat{\mathsf{G}}_{n}^{1,2}
\bigl(b,c_{\psi}+n^{\kappa
-1/2}d\bigr)=\mathbb{G}%
_{n}^{1,2}(0,c_{\psi})+n^{1/2}
\overline{\mathsf{G}}_{n}^{1,2}(0,c_{\psi
})+2\sigma
c_{\psi}^{2}\mathsf{f}(c_{\psi})n^{\kappa}
\sigma d+\mathrm{o} 
_{\mathsf{P}}(1),
\]
uniformly in $|b|,|d|\leq n^{1/4-\eta}B$, $0\leq\psi\leq1$. Combine the first
two terms to get
\[
n^{1/2}\widehat{\mathsf{G}}_{n}^{1,2}
\bigl(b,c_{\psi}+n^{\kappa
-1/2}d\bigr)=n^{1/2}%
\widehat{\mathsf{G}}_{n}^{1,2}(0,c_{\psi})+2
\sigma^{2}c_{\psi}^{2}%
\mathsf{f}(c_{\psi})n^{\kappa}d+\mathrm{o}_{\mathsf{P}}(1).
\]
Replace $c_{\psi}$ by $\hat{c}_{\psi}^{0}$. Since $n^{1/2-\kappa}(\hat
{c}_{\psi}^{b}-\hat{c}_{\psi}^{0})=\mathrm{O}_{\mathsf{P}}(n^{1/4-\eta})$
uniformly in $0\leq\psi\leq1$, $|b|\leq n^{1/4-\eta}B$ by
(\ref{ptbhatorderstat}), we can replace $n^{\kappa}d$ by
$n^{1/2}(\hat
{c}_{\psi}^{b}-\hat{c}_{\psi}^{0})$ on a set with large probability. Subtract
$n^{1/2}\tau_{\psi}\sigma^{2}$ on both sides. Add and subtract $n^{1/2}%
\tau_{\mathsf{G}(\hat{c}_{\psi}^{0})}\sigma^{2}$ on the right.
Altogether we
get%
%
\begin{eqnarray}\label{ptshatexpand}
&&n^{1/2}\bigl\{\widehat{\mathsf{G}}_{n}^{1,2}
\bigl(b,\hat{c}_{\psi}^{b}\bigr)-\tau_{\psi
}
\sigma^{2}\bigr\}\nonumber\\
&&\quad=n^{1/2}\bigl\{\widehat{\mathsf{G}}_{n}^{1,2}
\bigl(0,\hat{c}_{\psi}%
^{0}\bigr)-
\sigma^{2}\tau_{\mathsf{G}(\hat{c}_{\psi}^{0})}\bigr\}
\\
&&\qquad{}+2\sigma^{2}\bigl(\hat{c}_{\psi}^{0}
\bigr)^{2}\mathsf{f}\bigl(\hat{c}_{\psi
}^{0}
\bigr)n^{1/2}%
\bigl(\hat{c}_{\psi}^{b}-
\hat{c}_{\psi}^{0}\bigr)+\sigma^{2}n^{1/2}
\{\tau_{\mathsf
{G}%
(\hat{c}_{\psi}^{0})}-\tau_{\psi}\}+\mathrm{o}_{\mathsf{P}}(1),
\nonumber
\end{eqnarray}
uniformly in $0\leq\psi\leq1$, $|b|\leq n^{1/4-\eta}B$. The three terms are
analysed in turn.

2. \textitt{First term of} (\ref{ptshatexpand}). Since
$\overline{\mathsf{G}}_{n}^{1,2}(0,c)=\sigma^{2}\tau_{\mathsf{G}(c)}$, the
first term equals $\mathbb{G}_{n}^{1,2}(0,\hat{c}_{\psi}^{0})$. Lemmas
\ref{lGhat0}, \ref{tdhat0}(c) show that $\hat{c}_{\psi}^{0}=c_{\psi
+n^{-1/2}\phi}$ where $\phi=n^{1/2}\{\mathsf{G}(\hat{c}_{\psi}^{0}%
)-\psi\}=\mathbb{G}_{n}(c_{\psi})+\mathrm{o}_{\mathsf{P}}(1)$ is tight.
Tightness of $\mathbb{G}_{n}^{1,2}$ was established in Lemma~\ref
{lG}(c) under the Assumption~\ref{asFS}(i)(a), (ii)(b), (ii)(c), then implies
that the first term equals
$\mathbb{G}_{n}^{1,2}(0,c_{\psi})+\mathrm{o}_{\mathsf{P}}(1)$ uniformly in
$0\leq\psi\leq1$.

3. \textitt{The order of} $\hat{c}_{\psi}^{0}$ is $\mathrm{o}_{\mathsf{P}
}(n^{\nu/2})$ for some $\nu<\eta-\kappa\leq1/4$. The reason is that
$\hat
{c}_{\psi}^{0}\leq\max_{i\leq n}|\varepsilon_{i}|$, that $\mathsf{E}%
|\varepsilon_{i}|^{q}<\infty$ for some $q>2/(\eta-\kappa)$ by Assumption~\ref{asFS}(i)(a), so that $q(\eta-\kappa)/2>1+\epsilon$ for some
$\epsilon>0$. Thus, Boole's and Markov's inequalities imply that
$\mathsf{P}(\max_{i}|\varepsilon_{i}|>Cn^{\nu/2})\leq\sum_{i=1}^{n}%
\mathsf{P}(|\varepsilon_{i}|>Cn^{\nu/2})\leq n(Cn^{\nu/2})^{-q}\mathsf
{E}%
|\varepsilon_{i}|^{q}$ vanishes if $\nu=(\eta-\kappa)/(1+\epsilon)$.

4. \textitt{The order of} $c_{\psi}^{2}$ is $\mathrm{o}(n^{1/4-2\lambda
})$ for
some $\lambda>0$ when $\psi\leq1-n^{-1}$. Because $\mathsf
{E}|\varepsilon
_{i}|^{q}<\infty$ for some $q>8$ by Assumption~\ref{asFS}(i)(a),
$\mathsf{P}(|\varepsilon_{i}|>\sigma c_{\psi})\leq c_{\psi}^{-q}%
\mathsf{E}(|\varepsilon_{i}/\sigma|^{q})$ by the Markov inequality. Thus,
$c_{\psi}^{2}=\mathrm{O}\{(1-\psi)^{-2/q}\}$. In particular, for $\psi
\leq1-n^{-1}$, $c_{\psi}^{2}=\mathrm{O}(n^{2/q})=\mathrm
{o}(n^{1/4-2\lambda})$
for $1/4-2\lambda>2/q$ so that $\lambda<(q-8)/(8q)$.

5. \textitt{Second term of} (\ref{ptshatexpand}) vanishes. Indeed,
$\mathsf{f}(\hat{c}_{\psi}^{0})n^{1/2}(\hat{c}_{\psi}^{b}-\hat{c}_{\psi}
^{0})=\mathrm{o}_{\mathsf{P}}(n^{-\omega})$ for all $\omega<\eta-\kappa$
uniformly in $0\leq\psi\leq1$, $b\leq n^{1/4-\eta}B$ by Lemma~\ref{tdhatf}
using Assumption~\ref{asFS}(i)(a), (ii)(b), (ii)(c). By item 3 then $(\hat
{c}_{\psi
}^{0})^{2}=\mathrm{o}_{\mathsf{P}}(n^{\nu})$ for some $\nu<\eta-\kappa$
and an
$\omega$ exists so $\nu<\omega$.

6. \textitt{Third term of} (\ref{ptshatexpand}). We will argue that%
%
\begin{equation}
n^{1/2}(\tau_{\psi+n^{-1/2}\hat{\phi}}-\tau_{\psi})-c_{\psi}^{2}
\hat {\phi }=\mathrm{o}_{\mathsf{P}}(1), \label{ptshatexpand3}%
\end{equation}
for $\psi_{0}\leq\psi\leq n/(n+1)$ and $\hat{\phi}=-n^{1/2}\{\mathsf{G}%
(\hat{c}_{\psi}^{0})-\psi\}$. This suffices since Lemmas \ref{lGhat0},
\ref{tdhat0}(c) using Assumption~\ref{asFS}(c) show
%
\begin{equation}
\hat{\phi}=-\mathbb{G}_{n}(c_{\psi})+\mathrm{o}_{\mathsf{P}}
\bigl(n^{\zeta-1/4}\bigr), \label{display1}%
\end{equation}
for all $\zeta>0$ while item 4 shows $c_{\psi}^{2}=\mathrm
{o}(n^{1/4-2\lambda
})$ for some $\lambda>0$. This implies
\[
n^{1/2}\{\tau_{\mathsf{G}(\hat{c}_{\psi}^{0})}-\tau_{\psi}\}+c_{\psi}%
^{2}
\mathbb{G}_{n}(c_{\psi})=\mathrm{o}\bigl(n^{\zeta-2\lambda}
\bigr)+\mathrm{o}%
_{\mathsf{P}}(1)=\mathrm{o}_{\mathsf{P}}(1),
\]
as desired. To prove (\ref{ptshatexpand3}), write
\[
\mathcal{S}_{3}=n^{1/2}(\tau_{\psi+n^{-1/2}\phi}-
\tau_{\psi})-c_{\psi
}^{2}%
\phi=n^{1/2}\int_{c_{\psi}}^{c_{\psi+n^{-1/2}\phi}}
\bigl(u^{2}-c_{\psi}%
^{2}\bigr)2
\mathsf{f}(u)\,\mathrm{d}u.
\]
Change variable $y=\mathsf{G}(u)$, $\mathrm{d}y=2\mathsf{f}(u)\,\mathrm{d}u$, and Taylor
expand to
get
\[
\mathcal{S}_{3}=n^{1/2}\int_{\psi}^{\psi+n^{-1/2}\phi}
\bigl(c_{y}^{2}-c_{\psi
}%
^{2}
\bigr)\,\mathrm{d}y=\phi\bigl(c_{\psi^{\ast}}^{2}-c_{\psi}^{2}
\bigr),
\]
for some $\psi^{\ast}$ so $|\psi^{\ast}-\psi|\leq\phi$. Rewrite this,
for some
$\upsilon>0$ yet to be chosen,
\begin{eqnarray*}
\mathcal{S}_{3}&=&\bigl\{\psi(1-\psi)\bigr\}^{-2\upsilon}\biggl\{
\frac{\psi(1-\psi)}%
{\mathsf{f}(c_{\psi})}\biggr\}(c_{\psi^{\ast}}+c_{\psi})\\
&&{}\times  \biggl[
\frac{\phi}%
{\{\psi(1-\psi)\}^{1/2-\upsilon}} \biggr] \biggl[ \frac{\mathsf{f}(c_{\psi
})n^{1/2}(c_{\psi^{\ast}}-c_{\psi})}{\{\psi(1-\psi)\}^{1/2-\upsilon
}} \biggr] n^{-1/2}.
\end{eqnarray*}
The first component is
%
\begin{equation}
\bigl\{\psi(1-\psi)\bigr\}^{-2\upsilon}=\mathrm{O}\bigl(n^{2\upsilon}
\bigr), \label{display2}%
\end{equation}
for $\psi_{0}\leq\psi\leq n/(n+1)$. The second component is $\psi
(1-\psi)/\mathsf{f}(c_{\psi})=\mathrm{O}(c_{\psi})$ by Assumption~\ref{asFS}(i)(d).
Since $c_{\psi}=\mathrm{o}(n^{1/8-\lambda})$ for some
$\lambda>0$ due to item 4, then
\[
\mathcal{S}_{3}=(c_{\psi^{\ast}}+c_{\psi}) \biggl[
\frac{\phi}{\{\psi
(1-\psi)\}^{1/2-\upsilon}} \biggr] \biggl[ \frac{\mathsf{f}(c_{\psi}%
)n^{1/2}(c_{\psi^{\ast}}-c_{\psi})}{\{\psi(1-\psi)\}^{1/2-\upsilon
}} \biggr]
\mathrm{O}_{\mathsf{P}}\bigl(n^{2\upsilon+1/8-\lambda-1/2}\bigr).
\]
Evaluate this expression for $\phi$ replaced by $\hat{\phi}$. The first term
is $c_{\psi}^{\ast}+c_{\psi}\leq\hat{c}_{\psi}^{0}+2c_{\psi}=\mathrm{o}%
_{\mathsf{P}}(n^{1/8})$ due to items 3, 4 so that
\[
\mathcal{S}_{3}= \biggl[ \frac{\hat{\phi}}{\{\psi(1-\psi)\}^{1/2-\upsilon}
} \biggr] \biggl[
\frac{\mathsf{f}(c_{\psi})n^{1/2}(c_{\psi^{\ast
}}-c_{\psi}%
)}{\{\psi(1-\psi)\}^{1/2-\upsilon}} \biggr] \mathrm{o}_{\mathsf{P}%
}\bigl(n^{2\upsilon-\lambda-1/4}
\bigr).
\]
The first component is $\{\mathbb{G}_{n}(c_{\psi})+\mathrm{o}_{\mathsf
{P}%
}(n^{\zeta-1/4})\}/\{\psi(1-\psi)\}^{1/2-\upsilon}$ by (\ref{display1}). The
first normalised summand is $\mathrm{O}_{\mathsf{P}}(1)$ by the H\"{o}lder
continuity of Lemma~\ref{tBB}. The second summand is $\mathrm
{o}_{\mathsf{P}%
}(n^{\zeta-1/4})\mathrm{o}_{\mathsf{P}}(n^{1/2-\upsilon})$ for $\psi_{0}
\leq\psi\leq n/(n+1)$ as in (\ref{display2}). Thus,
\[
\mathcal{S}_{3}=\frac{\mathsf{f}(c_{\psi})n^{1/2}(c_{\psi^{\ast
}}-c_{\psi}%
)}{\{\psi(1-\psi)\}^{1/2-\upsilon}}\mathrm{o}_{\mathsf{P}}
\bigl(n^{\upsilon
+\zeta-\lambda}\bigr).
\]
For the first component note $|c_{\psi}^{\ast}-c_{\psi}|\leq|c_{\psi}%
+n^{-1/2}\phi-c_{\psi}|=|c_{\hat{\psi}}-c_{\psi}|$. Lemma~\ref{tOReilly}(b)
using Assumption~\ref{asFS}(i)(b) then implies that a sequence of Brownian
bridges $\mathbb{B}_{n}$ exists so that the first component is bounded by
$\mathrm{o}_{\mathsf{P}}(1)+|\mathbb{B}_{n}(\psi)|/\{\psi(1-\psi
)\}^{1/2-\upsilon}$ uniformly in $\psi_{0}\leq\psi\leq n/(n+1)$. This
in turn
is $\mathrm{O}_{\mathsf{P}}(1)$ by the H\"{o}lder continuity of Lemma~\ref{tBB}. Overall it follows that $\mathcal{S}_{3}=\mathrm{o}_{\mathsf
{P}%
}(n^{\upsilon+\zeta-\lambda})=\mathrm{o}_{\mathsf{P}}(1)$ since we can choose
$\upsilon+\zeta<\lambda$.
\end{pf}

\subsection{The forward plot of least squares estimators}

The forward plot of least squares estimators is now considered. The one-step
result in Lemma~\ref{tbhat} implies that the Forward Search iteration
can be
viewed as a fixed point problem. Indeed, the one-step result in Lemma~\ref{tbhat} implies an autoregressive relation between the one-step updated
estimation error $\hat{b}^{(m+1)}$ and the previous estimation error
$\hat
{b}^{(m)}$. That is,%
%
\begin{equation}
\hat{b}^{(m+1)}=\rho_{\psi}\hat{b}^{(m)}+(\psi
\Sigma_{n})^{-1}\mathbb{G} 
_{n}^{x,1}(0,c_{\psi})+e_{\psi}
\bigl(\hat{b}^{(m)}\bigr), \label{ARu}%
\end{equation}
for $\psi=m/n+\mathrm{o}(1)$, an ``autoregressive coefficient'' $\rho
_{\psi}$
defined in (\ref{rho}) and a vanishing remainder term $e_{\psi}$. This
autoregressive representation generalises Johansen and Nielsen \cite{JN2010}, Theorem~5.2,
which was concerned with a location-scale model, a fixed $\psi\sim
m/n$, and convergent initial estimators, $\hat{b}^{(m)}=\mathrm{O}(1)$.

It is first established that $\rho_{\psi}$ has nice properties for unimodal
densities $\mathsf{f}$.

\begin{lemma}
\label{trho}Suppose Assumption~\ref{asFS}\textup{(i)(a), (i)(c)} holds. Then $\rho
_{\psi
}=2c_{\psi}\mathsf{f}(c_{\psi})/\psi$ satisfies:
\begin{longlist}[(a)]
\item[(a)] $\rho_{\psi}>0$
for all $\psi>0$;

\item[(b)] $\sup_{\psi_{0}\leq\psi<1}\rho_{\psi}<1$ for
all $\psi_{0}>0$.
\end{longlist}
\end{lemma}

\begin{pf}
(a) holds because $\mathsf{f}(c_{\psi
})>0$ for
$0<\psi<1$.

(b) If the conclusion were incorrect, there would exist a sequence
$\psi
_{n}$ so that $\rho_{\psi_{n}}\rightarrow1$ for $n\rightarrow\infty$. Let
$\psi^{\dag}$ be a limit point. We consider the cases where $\psi^{\dag}<1$
and $\psi^{\dag}=1$.

Suppose $\psi^{\dag}<1$. Then $\rho_{\psi^{\dag}}=1$, which implies
$2c_{\psi^{\dag}}\mathsf{f}(c_{\psi^{\dag}})=\psi^{\dag}$. Since $\psi
^{\dag
}=2\int_{0}^{c_{\psi^{\dag}}}\mathsf{f}(x)\,\mathrm{d}x$ it holds that $\int_{0}%
^{c_{\psi^{\dag}}}\{\mathsf{f}(x)-\mathsf{f}(c_{\psi^{\dag}})\}\,\mathrm{d}x=0$. This
contradicts Assumption~\ref{asFS}(i)(c).

Suppose $\psi^{\dag}=1$. Because $\psi_{n}\rightarrow1$, it must hold
in this
case that $c_{\psi_{n}}\mathsf{f}(c_{\psi_{n}})\rightarrow1$ for
$n\rightarrow0$. This contradicts that $c\mathsf{f}(c)\rightarrow0$ for
$c\rightarrow\infty$ by Assumption~\ref{asFS}(i)(a).
\end{pf}

The next result investigates the forward estimator $\hat{\beta}^{(m+1)}$.
There are two results: first, the Forward Search preserves the order of the
initial estimator, and second, by infinite iteration a slowly converging
initial estimator can be improved to consistency at a standard rate.
The proof
of this result is related to that of Johansen and Nielsen \cite{JN2013}, Theorem~3.3.

\setcounter{theorem}{12}
\begin{theorem}
\label{tbhatmax}Suppose Assumption~\ref{asFS}\textup{(i)(a)--(i)(c), (ii), (iii)} holds, but
with $q_{0}=1+2^{r+1}$ only. Then, for all $\psi_{1}>\psi_{0}>0$ and
$m_{0}/n=\psi_{0}+\mathrm{o}(1)$, the estimator $\hat{\beta}_{\psi}$
satisfies:
\begin{longlist}[(a)]
\item[(a)] $\sup_{\psi_{0}\leq\psi\leq1}|N^{-1}(\hat{\beta
}_{\psi
}-\beta)|=\mathrm{O}_{\mathsf{P}}(n^{1/4-\eta})$;

\item[(b)] $\sup_{\psi
_{1}\leq\psi\leq1}|N^{-1}(\hat{\beta}_{\psi}-\beta)|=\mathrm{O}_{\mathsf
{P}%
}(1)$.
\end{longlist}
\end{theorem}

\begin{pf}
Due to the embedding (\ref
{embedding}), it
suffices to evaluate $N^{-1}(\hat{\beta}_{\psi}-\beta)$ at the grid points
$\psi=m/n$. Introduce notation $K_{\psi}^{n}=\Sigma_{n}^{-1}\mathbb{G}%
_{n}^{x,1}(0,c_{\psi})$.

(a) Solve the autoregressive equation (\ref{ARu}) recursively to get
\[
\hat{b}^{(m+1)}=%
{ \sum
_{k=m_{0}}^{m}} 
\Biggl(%
{
\prod_{\ell=k+1}^{m}} 
\rho_{\ell/n}\Biggr)\biggl\{\frac{n}{k}K_{k/n}^{n}+e_{k/n}
\bigl(\hat{b}^{(k)}\bigr)\biggr\}+\Biggl(%
{
\prod_{k=m_{0}}^{m}} 
\rho_{k/n}\Biggr)\hat{b}^{(m_{0})},
\]
with the convention that an empty product equals unity. Lemma~\ref{trho}
using Assumption~\ref{asFS}(i)(a), (i)(c) shows that $\rho_{\psi}\leq\rho
_{0}$ for
some $\rho_{0}<1$ for $\psi\geq\psi_{0}$, and therefore $%
{ \sum_{k=m_{0}}^{m}}
\rho_{0}^{m-k}\leq\sum_{k=0}^{\infty}\rho_{0}^{k}=C$. This gives the
bound%
%
\begin{equation}
\bigl|\hat{b}^{(m+1)}\bigr|\leq C\Bigl\{\sup_{\psi_{0}\leq\psi\leq1}\bigl|
\psi^{-1}K_{\psi}%
^{n}\bigr|+\max
_{m_{0}\leq k\leq m}\bigl|e_{k/n}\bigl(\hat{b}^{(k)}\bigr)\bigr|
\Bigr\}+\rho _{0}^{m-m_{0}+1}\bigl|\hat{b}^{(m_{0})}\bigr|.\label{ptbhatmaxbound}%
\end{equation}
In this expression, the process $\psi^{-1}K_{\psi}^{n}$ in $D[\psi
_{0},1]$ for
$\psi_{0}>0$, is tight by Lemma~\ref{lG}(c) using Assumption~\ref
{asFS}%
(i)(a), (ii)(b), (ii)(c). Therefore, for any $\epsilon>0$ we first choose $B$ so large
that $\mathsf{P}(C\sup_{\psi_{0}\leq\psi\leq1}|\psi^{-1}K_{\psi
}^{n}|\geq
B)\leq\epsilon/3$ for all $n$. The initial estimator is $\hat{b}^{(m_{0}
)}=\mathrm{O}_{\mathsf{P}}(n^{1/4-\eta})$ by Assumption~\ref{asFS}(iii),
and we next choose $B$ so large that \textsf{P}$(|\hat{b}^{(m_{0})}|\geq
Bn^{1/4-\eta})\leq\epsilon/3$ for all $n$. Finally, by Lemma~\ref
{tbhat}%
(c), $\sup_{\psi_{0}\leq\psi\leq1}\sup_{|b|\leq3n^{1/4-\eta}B}|e_{\psi
}(b)|=\mathrm{o}_{\mathsf{P}}(1)$, using Assumption~\ref{asFS}(i)(a), (i)(b), (ii).
Thus, there is an $n_{0}$ such that
\[
\mathsf{P}\Bigl(C\sup_{\psi_{0}\leq\psi\leq1}\sup_{|b|\leq3n^{1/4-\eta
}B}\bigl|e_{\psi
}(b)\bigr|
\geq B\Bigr)\leq\epsilon/3,
\]
for $n\geq n_{0}$. This implies that the set
\[
\mathcal{A}_{n}=\Bigl(C\sup_{\psi_{0}\leq\psi\leq1}\bigl|
\psi^{-1}K_{\psi
}^{n}\bigr|\leq B\Bigr)\cap\Bigl(C\sup
_{\psi_{0}\leq\psi\leq1}\sup_{|b|\leq3n^{1/4-\eta}B}\bigl|e_{\psi
}(b)\bigr|\leq B
\Bigr)\cap\bigl(\bigl|\hat{b}^{(m_{0})}\bigr|\leq n^{1/4-\eta}B\bigr)
\]
has probability larger than $1-\epsilon$. An induction over $m$ is now
used to
prove that
\[
\max_{m_{0}\leq k\leq m}\bigl|\hat{b}^{(k)}\bigr|\leq3n^{1/4-\eta}B\qquad
\mbox{for }m=m_{0},\ldots,n,
\]
on the set $\mathcal{A}_{n}$, which implies the desired result. For $m=m_{0}$,
the initial estimator satisfies $|\hat{b}^{(m_{0})}|\leq n^{1/4-\eta}B$
on the
set $\mathcal{A}_{n}$. Suppose the result holds for some $m$. This implies
that
%
\begin{equation}
C\sup_{\psi_{0}\leq\psi\leq1}\max_{m_{0}\leq k\leq m}\bigl|e_{\psi
}
\bigl(\hat {b}^{(k)}\bigr)\bigr|\leq B\label{ptbhatmaxe}%
\end{equation}
on the set $\mathcal{A}_{n}$. Thus, the bound (\ref{ptbhatmaxbound})
becomes
\[
\bigl|\hat{b}^{(m+1)}\bigr|\leq B+B+n^{1/4-\eta}B\leq3n^{1/4-\eta}B,
\]
because $n^{1/4-\eta}\geq1$ for $\eta\leq1/4$. Thus, the result holds for
$m+1$, which completes the induction.

(b) Consider next (\ref{ptbhatmaxbound}) for $\psi_{1}n\leq m\leq n$.
Here, $%
{ \sum_{k=0}^{n}}
\rho_{0}^{k}\leq C$, the first term is $\sup_{\psi_{1}\leq\psi\leq
1}|\psi
^{-1}K_{\psi}^{n}|=\mathrm{O}_{\mathsf{P}}(1)$ due to tightness, while the
second, as remarked above, is $\sup_{\psi_{1}\leq\psi\leq1}\max_{m_{0}%
\leq k<n}|e_{\psi}(\hat{b}^{(k)})|=\mathrm{o}_{\mathsf{P}}(1)$. Because
$\rho_{0}^{m-m_{0}}\leq\rho_{0}^{\mathrm{int}(\psi_{1}n)-\operatorname
{int}(\psi
_{0}n)}$ declines exponentially, $\rho_{0}^{m-m_{0}}<n^{-1/4}$ for
large $n$
and therefore the last term is $\max_{m\geq m_{1}}\rho
_{0}^{m-m_{0}}|\hat
{b}^{(m_{0})}|=\mathrm{O}_{\mathsf{P}}(n^{-1/4+1/4-\eta})=\mathrm{o}%
_{\mathsf{P}}(1)$, which proves (b).
\end{pf}

\subsection{Proofs of main Theorems \texorpdfstring{\protect\ref{txi}}{3.1}--\texorpdfstring{\protect\ref{tK}}{3.7}}

Lemmas \ref{tdhat0}, \ref{tdhatalg} are now combined to show that the
forward residuals scaled with a known variance, $\sigma^{-1}\hat
{z}_{\psi}$,
have the same Bahadur representation as the quantile process for the
innovations $\sigma^{-1}\varepsilon_{i}$. This is the main theorem
stated with
slightly weaker conditions.

\begin{remark}
The proof below of Theorem~\ref{txi} only requires Assumption~\ref
{asFS}%
(i)(a)--(i)(c), (ii), (iii) with $q_{0}=1+2^{r+1}$.
\end{remark}

\begin{pf*}{Proof of Theorem~\ref{txi}}
It is first argued that the forward plot
of the
estimators is bounded in the sense that for all $\epsilon>0$ a $B>0$
exists so
that the set
\[
\mathcal{C}_{n}=\Bigl(\sup_{\psi_{0}\leq\psi\leq1}\bigl|N^{-1}(
\hat{\beta}_{\psi}%
-\beta)\bigr|\leq n^{1/4-\eta}B\Bigr)
\]
has $\mathsf{P}(\mathcal{C}_{n})\geq1-\epsilon$. This follows from Lemma~\ref{tbhatmax} using Assumption~\ref{asFS}(i)(a)--(i)(c), (ii), (iii).
Now, on
$\mathcal{C}_{n}$ it holds that $\sigma^{-1}\hat{z}_{\psi}=\hat{c}_{\psi
}%
^{b}$, see (\ref{xihat}), for some $|b|\leq n^{1/4-\eta}B$. Thus it suffices
to show that
\[
\sup_{\psi_{0}\leq\psi\leq n/(n+1)}\sup_{|b|\leq n^{1/4-\eta}B}\bigl|\mathbb
{C}%
_{\psi}^{b}\bigr|=\mathrm{o}_{\mathsf{P}}(1)
\qquad\mbox{for }\mathbb {C}_{\psi
}^{b}=2\mathsf{f}(c_{\psi})n^{1/2}
\bigl(\hat{c}_{\psi}^{b}-c_{\psi}\bigr)+\mathbb
{G}%
_{n}^{1,0}(c_{\psi}).
\]
Now, write $(\hat{c}_{\psi}^{b}-c_{\psi})=(\hat{c}_{\psi}^{0}-c_{\psi}%
)+(\hat{c}_{\psi}^{b}-\hat{c}_{\psi}^{0})$, so that
\[
\mathbb{C}_{\psi}^{b}=\bigl\{2\mathsf{f}(c_{\psi})n^{1/2}
\bigl(\hat{c}_{\psi}%
^{0}-c_{\psi}\bigr)+
\mathbb{G}_{n}^{1,0}(c_{\psi})\bigr\}+2
\frac{\mathsf
{f}(c_{\psi}%
)}{\mathsf{f}(\hat{c}_{\psi}^{0})}n^{1/2}\mathsf{f}\bigl(\hat{c}_{\psi}^{0}%
\bigr) \bigl(\hat{c}_{\psi}^{b}-\hat{c}_{\psi}^{0}
\bigr).
\]
The first term is $\mathrm{o}_{\mathsf{P}}(n^{\zeta-1/4})$ for all
$\zeta>0$
uniformly in $0\leq\psi\leq1$ by Lemma~\ref{tdhat0}(a) using Assumption~\ref{asFS}(i)(b). In the second term, the ratio $\mathsf{f}(c_{\psi
})/\mathsf{f}(\hat{c}_{\psi}^{0})$ is $\mathrm{O}_{\mathsf{P}}(1)$ uniformly
in $0\leq\psi\leq n/(n+1)$ by Lemma~\ref{lf2f} using Assumption~\ref{asFS}(i)(a), (i)(b), while $n^{1/2}\mathsf{f}(\hat{c}_{\psi}^{0})(\hat
{c}_{\psi}^{b}-\hat{c}_{\psi}^{0})=\mathrm{o}_{\mathsf{P}}(n^{-\omega})$
uniformly in $0\leq\psi\leq1$ by Lemma~\ref{tdhatf} using Assumption~\ref{asFS}(i)(a), (ii)(b), (ii)(c).
Combining the first statement with Lemma~\ref{tdhat0}(a) gives the second statement.
\end{pf*}

\begin{remark}
The proof below of Theorem~\ref{tshat} only requires Assumption~\ref{asFS}(i)(a), (i)(b), (i)(d), (ii).
\end{remark}

\begin{pf*}{Proof of Theorem~\ref{tshat}}
The above theory for $\sigma^{-1}\hat
{z}_{\psi
}$ involves the population variance $\sigma^{2}$. The result gives an
asymptotic expansion for $\hat{\sigma}_{\psi,\mathrm{cor}}^{2}$, recalling, from
(\ref{sigma^m}), (\ref{shatcorrect}), (\ref{Ghat}) that%
%
\begin{eqnarray}\label{sigmaexpand}
&& n^{1/2}\bigl(\hat{\sigma}_{\psi,\mathrm{cor}}^{2}-
\sigma^{2}\bigr)
\nonumber
\\[-8pt]
\\[-8pt]
\nonumber
&&\qquad =\frac{1}{\tau_{\psi}}n^{1/2}\bigl[\bigl\{\widehat{\mathsf{G}}_{n}^{1,2}
\bigl(\hat{b} 
,\hat{c}_{\psi}^{\hat{b}}\bigr)-
\tau_{\psi}\sigma^{2}\bigr\}-\bigl\{\widehat{\mathsf
{G}}%
_{n}^{x,1}\bigl(\hat{b},
\hat{c}_{\psi}^{\hat{b}}\bigr)\bigr\}^{\prime}\bigl\{\widehat
{\mathsf{G}%
}_{n}^{xx,0}\bigl(\hat{b},
\hat{c}_{\psi}^{\hat{b}}\bigr)\bigr\}^{-1}\bigl\{\widehat{
\mathsf {G}%
}_{n}^{x,1}\bigl(\hat{b},
\hat{c}_{\psi}^{\hat{b}}\bigr)\bigr\}\bigr].\quad
\end{eqnarray}
Compare also the definitions in (\ref{GG}), (\ref{HH}) with (\ref
{Gtilde}) to
see
%
\begin{equation}
\mathbb{G}_{n}(c_{\psi})=\mathbb{G}_{n}^{1,0}(0,c_{\psi}),\qquad
\tau _{\psi
}\mathbb{L}_{n}(c_{\psi})=
\sigma^{-2}\mathbb{G}_{n}^{1,2}(0,c_{\psi
})-c_{\psi
}^{2}
\mathbb{G}_{n}^{1,0}(0,c_{\psi}).
\end{equation}

Lemma~\ref{lshat} using Assumption~\ref{asFS}(i)(a), (i)(b), (i)(d), (ii) shows
the first
term in (\ref{sigmaexpand}) equals the leading term $\mathbb
{L}_{n}(c_{\psi
})+\mathrm{o}_{\mathsf{P}}(1)$ uniformly in $\psi_{0}\leq\psi\leq n/(n+1)$
while the second term in~(\ref{sigmaexpand}) vanishes.
\end{pf*}

\begin{pf*}{Proof of Theorem~\ref{txisigma}}
Note the identity
\[
\frac{\hat{z}_{\psi}}{\hat{\sigma}_{\psi,\mathrm{cor}}}-c_{\psi}=\frac{\hat
{z}_{\psi
}/\sigma-c_{\psi}}{\hat{\sigma}_{\psi,\mathrm{cor}}/\sigma}-c_{\psi}
\frac{\hat
{\sigma
}_{\psi,\mathrm{cor}}^{2}-\sigma^{2}}{\hat{\sigma}_{\psi,\mathrm{cor}}(\hat{\sigma}_{\psi
,\mathrm{cor}}+\sigma)}.
\]
Multiply this by $2\mathsf{f}(c_{\psi})n^{1/2}$. Use that $2\mathsf{f}%
(c_{\psi})n^{1/2}(\hat{z}_{\psi}/\sigma-c_{\psi})$ and $n^{1/2}(\hat
{\sigma
}_{\psi,\mathrm{cor}}^{2}/\sigma^{2}-1)$ have the leading terms $-\mathbb{G}%
_{n}(c_{\psi})$ and $\mathbb{L}_{n}(c_{\psi})$, respectively, due to Theorems
\ref{txi}, \ref{tshat}. In particular $\hat{\sigma}_{\psi},_{\mathrm{cor}}$ is
consistent for $\sigma$.
\end{pf*}

\begin{pf*}{Proof of Theorem~\ref{tdeletion}}
We first show that $\hat{d}^{(m)}
\leq\hat{z}^{(m)}$ and then we find an upper bound for $\hat{z}%
^{(m)}-\hat{d}^{(m)}$, and finally show that the difference is small.

1. \textitt{Inequality} $\hat{d}^{(m)}\leq\hat{z}^{(m)}$.
Indeed, if
$S^{(m)}$ is the ranks of $\hat{\xi}_{(1)}^{(m)},\ldots,\hat{\xi}
_{(m)}^{(m)}$ then $\hat{d}^{(m)}=\hat{z}^{(m)}$. If $S^{(m)}$ does
not have this form, then its complement must include one of the ranks of
$\hat{\xi}_{(1)}^{(m)},\ldots,\hat{\xi}_{(m)}^{(m)}$, for
instance that
of $i^{\dag}$. In that situation $\hat{d}^{(m)}\leq\hat{\xi
}_{i^{\dag
}}^{(m)}\leq\hat{\xi}_{(m)}^{(m)}\leq\hat{\xi}_{(m+1)}^{(m)}%
=\hat{z}^{(m)}$.

2. \textitt{The set} $S^{(m)}$ \textitt{consists} of the ranks of
$\hat{\xi
}_{(1)}^{(m-1)},\ldots,\hat{\xi}_{(m)}^{(m-1)}$. It follows that for all
$i\notin S^{(m)}$ then $\hat{\xi}_{i}^{(m-1)}\geq\hat{\xi}%
_{(m+1)}^{(m-1)}\geq\hat{\xi}_{(m)}^{(m-1)}=\hat{z}^{(m-1)}$.

3. \textitt{Inequality for deletion residual.} The absolute residual for
observation $i$ based on the set $S^{(m)}$, $\xi_{i}^{(m-1)}$ in step $m-1$,
satisfies
\begin{eqnarray*}
\hat{\xi}_{i}^{(m-1)}&=&\bigl|y_{i}-x_{i}^{\prime}
\hat{\beta}^{(m-1)}\bigr|\leq \bigl|y_{i}-x_{i}^{\prime}
\hat{\beta}^{(m)}\bigr|+\bigl|x_{i}^{\prime}\bigl(\hat {
\beta }^{(m)}-\hat{\beta}^{(m-1)}\bigr)\bigr|
\\
&\leq&\hat{\xi}_{i}^{(m)}+\max_{1\leq i\leq n}\bigl|N^{\prime}x_{i}\bigr|\bigl|N^{-1}%
\bigl(\hat{\beta}^{(m)}-\hat{\beta}^{(m-1)}\bigr)\bigr|.
\end{eqnarray*}
For $i\notin S^{(m)}$, we have from item 2 that $\xi_{i}^{(m-1)}\geq
\hat{\xi
}_{(m)}^{(m-1)}=\hat{z}^{(m-1)}$ and $\hat{d}^{(m)}=\min_{i\notin
S^{(m)}}\hat{\xi}_{i}^{(m)}$ giving
\[
\hat{z}^{(m-1)}\leq\hat{d}^{(m)}+\max_{1\leq i\leq n}\bigl|N^{\prime}%
x_{i}\bigr|\bigl|N^{-1}
\bigl(\hat{\beta}^{(m)}-\hat{\beta}^{(m-1)}\bigr)\bigr|,
\]
and therefore, using $\hat{d}^{(m)}\leq\hat{z}^{(m)}$ we find
%
\begin{equation}
0\leq\hat{z}^{(m)}-\hat{d}^{(m)}\leq
\hat{z}^{(m)}-\hat {z}%
^{(m-1)}+\bigl|N^{-1}
\bigl(\hat{\beta}^{(m)}-\hat{\beta}^{(m-1)}\bigr)\bigr|\max
_{i}|Nx_{i}|. \label{ptdeletionineq1}%
\end{equation}

4. \textitt{Embed in the interval} $[0,1]$ \textitt{using} $\psi=m/n$. The
asymptotic expansion for $\hat{z}^{(m)}$ in Theorem~\ref{txi} combined
with the tightness of $\mathbb{G}_{n}$ in Lemma~\ref{tbhatmax} shows
\[
\sup_{\psi_{0}\leq\psi\leq n/(n+1)}\bigl|2\mathsf{f}(c_{\psi}) (\hat
{z}_{\psi
}-\hat{z}_{\psi-1/n})\bigr|=\mathrm{o}_{\mathsf{P}}
\bigl(n^{-1/2}\bigr),
\]
while the asymptotic result for $\hat{\beta}^{(m)}$ in Lemma~\ref{tbhatmax}
shows
\[
\sup_{\psi_{1}\leq\psi\leq n/(n+1)}\bigl|N^{-1}\bigl(\hat{
\beta}^{(m)}%
-\hat{\beta}^{(m-1)}\bigr)\bigr|=
\mathrm{o}_{\mathsf{P}}\bigl(n^{-1/2}\bigr).
\]

5. \textitt{Combine.} The bound (\ref{ptdeletionineq1}) and the triangle
inequality give
\begin{eqnarray*}
0 & \leq&2\mathsf{f}(c_{m/n}) \bigl(\hat{z}^{(m)}-
\hat{d}^{(m)}\bigr)
\\
& \leq&2\mathsf{f}(c_{m/n})\bigl|\hat{z}^{(m)}-
\hat{z}^{(m-1)}%
\bigr|+2\mathsf{f}(c_{m/n})|N^{-1}
\bigl(\hat{\beta}^{(m)}-\hat{\beta}%
^{(m-1)}
\bigr)|\max_{i}|Nx_{i}|.
\end{eqnarray*}
The bounds in item 4, combined with the condition $\max_{i}|Nx_{i}%
|=\mathrm{O}_{\mathsf{P}}(n^{\kappa-1/2})$ for some $\kappa<\eta\leq
1/4$ by
Assumption~\ref{asFS}(ii)(b), give a further bound
\[
0\leq2\mathsf{f}(c_{m/n}) \bigl(\hat{z}^{(m)}-
\hat{d}^{(m)}\bigr)\leq \mathrm{o}_{\mathsf{P}}
\bigl(n^{-1/2}\bigr)+\mathrm{o}_{\mathsf{P}}\bigl(n^{-1/2}%
n^{\kappa-1/2}
\bigr)=\mathrm{o}_{\mathsf{P}}\bigl(n^{-1/2}\bigr),
\]
as desired.
\end{pf*}

\begin{pf*}{Proof of Theorem~\ref{tbhatexpand}}
Lemma~\ref{tbhat}(c) using Assumption~\ref{asFS}(i)(a), (i)(b), (ii) shows
\[
b^{\dag}=\bigl\{\widehat{\mathsf{G}}_{n}^{xx,0}
\bigl(b,\hat{c}_{\psi}^{b}\bigr)\bigr\}^{-1}%
\bigl\{n^{1/2}\widehat{\mathsf{G}}_{n}^{x,1}\bigl(b,
\hat{c}_{\psi}^{b}\bigr)\bigr\}=(\Sigma _{n}
\psi)^{-1}\mathbb{G}_{n}^{x,1}(0,c_{\psi})+
\rho_{\psi}b+\mathrm{o}%
_{\mathsf{P}}(1),
\]
uniformly in $|b|\leq n^{1/4-\eta}B$, $\psi_{0}\leq\psi\leq1$. Lemma~\ref{tbhatmax}(b) using Assumption~\ref{asFS}(i)(a)--(i)(c), (ii), (iii)
shows that
$N^{-1}(\hat{\beta}_{\psi}-\beta)$ is uniformly bounded for $\psi\geq
\psi
_{1}$. Thus, on a set with large probability both $b^{\dag}$ and $b$
can be
replaced by $N^{-1}(\hat{\beta}_{\psi}-\beta)+\mathrm{o}_{\mathsf{P}}(1)$.
Lemma~\ref{trho} using Assumption~\ref{asFS}(i)(a), (i)(c) shows that $\rho
_{\psi}\leq\rho_{0}<1$ for $\psi\geq\psi_{0}$. Thus, it holds
\[
N^{-1}(\hat{\beta}_{\psi}-\beta)=\frac{1}{1-\rho_{\psi}}(
\Sigma_{n}\psi )^{-1}\mathbb{G}_{n}^{x,1}(0,c_{\psi})+
\mathrm{o}_{\mathsf{P}}(1).
\]
Insert $\rho_{\psi}=2c_{\psi}\mathsf{f}(c_{\psi})/\psi$ and $\mathbb{K}%
_{n}(c_{\psi})=\mathbb{G}_{n}^{x,1}(0,c_{\psi})$ to get the desired expansion.
\end{pf*}

\begin{pf*}{Proof of Theorems \ref{tGhat0} and \ref{tK}}
Tightness follows from Lemma~\ref{lG}(c), and convergence of finite dimensional distributions
follows from the central limit theorem for martingale differences, see
Helland \cite{Helland}, Theorem~3.2b, using Assumption \ref
{asFS}(ii)(c).
\end{pf*}

\section{A result on order statistics of $\mathsf{t}$-distributed
variables}
\label{st-order}

\begin{theorem}
\label{torderedt}Let $v_{1},\ldots,v_{n}$ be independent absolute
$\mathsf{t}_{m\mbox{-}\dim x}$ distributed. Consider the $(m+1)^{\prime}$st smallest
order statistic $\hat{v}_{(m+1)}^{(m)}$. Suppose $\dim x$ is fixed while
$m\sim\psi n$ for some $0<\psi<1$. Let $\varphi$ be the standard normal
density. Then, as $n\rightarrow\infty$,
\[
2\varphi(c_{m/n})n^{1/2}\bigl(\hat{v}_{(m+1)}^{(m)}-c_{m/n}
\bigr)\overset{\mathsf {D}%
} {\rightarrow}\mathsf{N}\bigl\{0,\psi(1-
\psi)\bigr\}.
\]
\end{theorem}

\begin{pf*}{Sketch of the proof of Theorem~\ref{torderedt}}
Let $\hat
{v}_{(m+1)}^{(m)}$ be the $(m+1)^{\prime}$st quantile of a sample of $n$
scaled, absolute $\mathsf{t}_{m-\dim x}$ variables. To get a handle on the
asymptotic distribution of $\hat{v}_{(m+1)}^{(m)}$ consider first the
$(m+1)^{\prime}$st smallest order statistic, $\hat{w}_{(m+1)}$ say,
from $n$
draws of absolute standard normal variables with distribution function
$2\Phi(y)-1$. This satisfies
\[
2\varphi(c_{m/n})n^{1/2}(\hat{w}_{(m+1)}-c_{m/n})
\overset{\mathsf{D}%
} {\rightarrow}\mathsf{N}\bigl\{0,\psi(1-\psi)\bigr
\},
\]
for $m\sim\psi n$ and $c_{\psi}=\mathsf{G}^{-1}(\psi)$ due to Lemmas
\ref{lGhat0}, \ref{tdhat0}(a). For the $\mathsf{t}_{m-\dim x}$ order
statistic $\hat{v}_{(m+1)}^{(m)}$ it is useful to Edgeworth expand
$\mathsf{P}(\mathsf{t}_{m-\dim x}\leq y)=2\{\Phi(y)+\mathrm{O}(n^{-1})\}-1$,
for $m\sim\psi n$, which indicates that the same asymptotic distribution
arises as in the normal case. A more formal argument will keep track of the
remainder terms. The starting point could be the expression for $\mathsf
{P}%
(\hat{v}_{(m+1)}^{(m)}\leq y)$ in terms of the distribution of an
$\mathsf{F}$
variate as given in Guenther \cite{Guentherequation3}, equation (3). This can be
expanded using
the approximation to the $\log\mathsf{F}$ distribution by Aroian \cite{Aroian}, Section~15. These considerations lead to the result.
\end{pf*}
\end{appendix}

\section*{Acknowledgements}
The first author is grateful to CREATES -- Center for
Research in Econometric Analysis of Time Series (DNRF78), funded by the Danish
National Research Foundation. We thank two anonymous referees for many
constructive suggestions for
improvement of the manuscript, and Xiyu Jiao and James Duffy for a careful
reading of the final manuscript.


%




\printhistory
\end{document}